\documentclass[11pt, reqno]{amsart}
\PassOptionsToPackage{dvipsnames}{xcolor}
\usepackage{eucal}
\usepackage{lmodern}
\usepackage[dvipsnames]{xcolor}
\usepackage[hidelinks]{hyperref}
\usepackage{comment}
\usepackage{color}
\usepackage{amsmath,amssymb,graphicx}
\usepackage{tikz,tikz-cd}
\usepackage{xparse}
\usepackage{fancyhdr}
\usepackage{caption}
\usepackage{longtable}
\usepackage{mathrsfs}

\DeclareRobustCommand{\g}{\ensuremath{\mathrm{g}}}
\DeclareRobustCommand{\TT}{\ensuremath{\mathsf{T}^2}}
\DeclareRobustCommand{\dd}{\ensuremath{\text{d}}}

\newcommand{\Vca}{\mathcal{V}}
\newcommand{\Hca}{\mathcal{H}}
\newcommand{\Rm}{\operatorname{Rm}}
\newcommand{\R}{\mathbb{R}}
\newcommand{\Lie}{\operatorname{Lie}}
\newcommand{\vol}{\operatorname{vol}}

\newcommand{\etau}{\hbox{\teneur\char28}}

\newtheorem{maintheorem}{Theorem}

\newtheorem{theorem}{Theorem}[section]
\newtheorem{proposition}[theorem]{Proposition} 
\newtheorem{lemma}[theorem]{Lemma}
\newtheorem{corollary}[theorem]{Corollary}

\theoremstyle{definition}
\newtheorem{remark}[theorem]{Remark} 
\newtheorem{definition}[theorem]{Definition} 
\newtheorem{notation}[theorem]{Notation}

\theoremstyle{remark}
\newtheorem{example}[theorem]{Example}

\title{$\mathsf{T}^2$-invariant Einstein metrics on closed $4$-manifolds}

\author{Kyle Broder}
\address{School of Mathematics and Physics, The University of Queensland, St Lucia QLD 4072, Australia}
\email{k.broder@uq.edu.au}

\author{Ramiro A. Lafuente}
\address{School of Mathematics and Physics, The University of Queensland, St Lucia QLD 4072, Australia}
\email{r.lafuente@uq.edu.au}

\thanks{Both authors are supported by the Australian Research Council project DP240101772. The second-named author is also supported by the Australian Research Council project FT240100361. The first-named author acknowledges  support from OpenAI's Academic Researcher Scholarship. }

\makeatletter
\@ifundefined{subjclassname@2020}{%
  \@namedef{subjclassname@2020}{\textup{2020} Mathematics Subject Classification}}{}
\makeatother

\font\teneur=eurm10

\makeatletter
\renewcommand{\section}{%
  \@startsection{section}{1}{\z@}%
    {.7\linespacing\@plus\linespacing}%
    {.5\linespacing}%
    {\normalfont\bfseries\centering}%
}
\makeatother

\usetikzlibrary{arrows.meta}

\hypersetup{pdftitle={T2-invariant Einstein metrics on closed 4-manifolds},
pdfauthor={Kyle Broder and Ramiro A. Lafuente},
pdfkeywords={Einstein metrics, torus actions, four-manifolds, Weyl curvature, conical metrics}}

\begin{document}
\raggedbottom

\begin{abstract}
We construct $\mathsf{T}^2$-invariant Einstein metrics with positive Ricci curvature on $\sharp_m(S^2\times S^2)$ for every integer $m>1$. These provide  infinite family of Einstein metrics on closed $4$-manifolds with $\lambda>0$, which are in addition  not  Hermitian. We further establish the uniqueness of the $\TT$-invariant Einstein metrics on $S^4$, $\mathbb{CP}^2$, $S^2\times S^2$, and $\mathbb{CP}^2\mathbin{\sharp}\overline{\mathbb{CP}}^2$.
\end{abstract}

\maketitle

\begingroup
\hypersetup{linkcolor=red}
\color{red}
\renewcommand{\baselinestretch}{1}

\tableofcontents
\endgroup

\section{Introduction}

\noindent A Riemannian metric $\g$ is \emph{Einstein} if $\text{Ric}(\g)=\lambda\g$ for some constant $\lambda\in\mathbb{R}$. Determining which manifolds admit Einstein metrics is a basic problem in Riemannian geometry (see, e.g., \cite{besse1987einstein}). Previously, Einstein metrics with $\lambda>0$ on closed simply connected $4$-manifolds were known on only eleven diffeomorphism types: $S^4$, $S^2 \times S^2$, and $\mathbb{CP}^2 \sharp k \overline{\mathbb{CP}^2}$ for $0 \leq k \leq 8$. The examples were the standard symmetric metrics on $S^4$, $S^2 \times S^2$, and $\mathbb{CP}^2$, the Page metric \cite{Page1978CompactRotatingInstanton} for $k=1$, the Chen--LeBrun--Weber metric \cite{ChenLeBrunWeber2008ConformallyKahlerEinstein} for $k=2$, and the K\"ahler--Einstein metrics \cite{Tian1990CalabiConjectureSurfaces} for $3 \leq k \leq 8$. In particular, all previously known examples are locally Hermitian. 

The first main theorem produces an infinite family of diffeomorphism types admitting Einstein metrics with $\lambda>0$.

\begin{maintheorem}\label{thm:intro-existence}
For each $m \in \mathbb{N}$, there exists a $\TT$-invariant Einstein metric with $\lambda>0$ on $\sharp_m(S^2\times S^2)$.
\end{maintheorem}  

Since $b_2^+=m$, LeBrun's classification \cite{LeBrun2012EinsteinHermitian} shows that, for $m>1$, these Einstein metrics are not  Hermitian. 

The idea for the proof of Theorem~\ref{thm:intro-existence} begins by realising $\sharp_m(S^2 \times S^2)$ as a degree-$(m+1)$ equivariant branched covering of $S^4$, branched along a principal orbit $\mathcal{T}$. We show that if $S^4 \setminus \mathcal{T}$ admits an Einstein metric with \emph{conical singularties along $\mathcal{T}$} and cone angle $2\pi/(m+1)$, then its pullback extends to a smooth $\TT$-invariant Einstein metric on $\sharp_m (S^2 \times S^2)$ (Proposition~\ref{prop:removability}). 

More generally, we prove existence of toric conical Einstein metrics for any cone angle $2\pi \beta$, $\beta \in (0,1/2]$, by using a mod-2 degree argument for the Einstein equation in this setting, vieweing it as a family of Fredholm operators of index zero between Banach spaces. This is inspired in part by earlier work of Anderson \cite{AndersonAH}, and is in turn  based on a classical framework due to Smale \cite{Smale1965Sard}.

The use of conical metrics to produce (smooth) Einstein metrics is well-2 \cite{Donaldson2012Cone, CDS1, Tian2015KStability} (cf~\cite{JMR}), as is the use of branched coverings \cite{FinePremoselli2020, HamenstadtJackel2024NegativelyCurvedEinstein, GuenanciaHamenstadt2025KahlerEinstein}. While previous conical arguments involved deforming the cone angle $2\pi \beta \to 2\pi$, our approach is to deform the angle \emph{downwards} from $\pi$. In the precise setting of $4$-manifolds, a particularly relevant reference is \cite{AtiyahLeBrun}. In this regard, we mention that while our precise definition  of conical metric, given in $\S$\ref{sec:space_conical} (Definition \ref{def:X_beta}), is more general, under the Einstein condition, regularity improves and after a radial reparametrisation one recovers the edge-cone notion from \cite{AtiyahLeBrun} (see the proof of Lemma \ref{lem:global-bounds}). Other related works on Einstein $4$-manifolds with symmetry include \cite{CP02,LS25}. 

To implement the mod-two degree argument, we consider an Ansatz  that reduces the toric conical Einstein equation to three scalar equations on the orbit space (Proposition~\ref{prop:reduced-equations}). We also impose an additional $\mathbb Z_2 \times \mathbb Z_2$  symmetry for technical reasons (see $\S$\ref{sec:branched-cover}). 
Near $\mathcal{T}$, $\TT$-invariance removes derivatives along $\mathcal{T}$ from Mazzeo's edge operators \cite{Mazzeo}, ensuring that to study these conical metrics we need only invoke Melrose's $b$-calculus \cite{Melrose1993}. Beyond the linear Fredholm theory (Theorem \ref{thm:scalar-index-zero}), a key ingredient for the degree theory to work is a compactness theorem (Theorem \ref{thm:compact-zero-set}), showing that  solutions cannot escape to infinity as we vary the cone angle. This is obtained  by performing a blow-up analysis and keeping track of good cordinate systems. Finally, we must also compute the degree at some known metric. This is achieved by showing that the standard Einstein metric on $S^2 \times S^2$ is non-degenerate (in the sense that it has no toric Einstein deformations), and unique under the symmetry assumptions.

In this direction, our next main  result establishes  uniqueness for toric Einstein metrics on $S^2\times S^2$ and other spaces, thus extending previous uniqueness results in the cohomogeneity-one setting by the first-named author \cite{BroderCohom1}:

\begin{maintheorem}\label{thm:intro-uniqueness}
Every Einstein metric admitting an effective isometric $\TT$-action on
\[
  S^4,\qquad \mathbb{CP}^2,\qquad S^2\times S^2,\qquad \mathbb{CP}^2\mathbin{\sharp}\overline{\mathbb{CP}}^2
\]
is, up to scaling and isometry, respectively the round metric, the Fubini--Study metric, the equal-radius round product, or the Page metric.
\end{maintheorem}

The proof of Theorem~\ref{thm:intro-uniqueness} uses a series of  subtle maximum principles to show that a certain distinguished Weyl eigenvalue must be negative (Proposition \ref{prop:weyl-sign}), or else the manifold is half-conformally flat.
For $S^4$ and $\mathbb{CP}^2$, the boundary behavior of the trace-normalized Gram map and a non-degeneracy criterion for planar elliptic equations force Weyl components to vanish; for $S^2\times S^2$ and $\mathbb{CP}^2\sharp\overline{\mathbb{CP}}^2$, we combine Gauss--Bonnet with Wu's maximum principle argument \cite{Wu2021} to obtain repeated Weyl eigenvalues. This allows us to then apply the results of Derdzi\'nski \cite{Derdzinski1983SelfDualKahler} and LeBrun \cite{LeBrun2012EinsteinHermitian} to identify the metrics.

Uniqueness for $S^4$, $\mathbb{CP}^2$ and $S^2\times S^2$ under the additional assumption of non-negative curvature was shown by Liu in \cite{Liu2025ToricEinstein4Manifolds}.

The Weyl curvature estimates also give the following refinement of the Hitchin--Thorpe inequality in the toric case. This allows us to rule out the existence of $\TT$-invariant Einstein metrics on several infinite families in the Orlik--Raymond classification \cite{OrlikRaymond1}. 

\begin{maintheorem}\label{thm:Hitchin-Thorpe}
A closed simply connected $4$-manifold admitting a $\TT$-invariant Einstein metric satisfies 
\[
 3|\etau(M)|\leq\chi(M).
\]
\end{maintheorem}

The rest of the article is organised as follows. Section \ref{sec:general-equations} contains preliminaries on $\TT$-invariant metrics, the Weyl curvature estimates, and the proof of Theorem \ref{thm:Hitchin-Thorpe}. In Section \ref{sec:uniq} we prove Theorem \ref{thm:intro-uniqueness}. The remaining sections \ref{sec:strategy}--\ref{sec:half-angle} are devoted to proving Theorem \ref{thm:intro-existence}. In Section \ref{sec:strategy} we present the general strategy for the proof, as well as the details on the mod-two degree theory approach. Section \ref{sec:conical} sets up our Ansatz for conical toric Einstein metrics, and the corresponding function spaces, in the language of Melrose's b-calculus. Section \ref{sec:fredholm} establishes the Fredholm property and computes the index of the Einstein operator. In Section \ref{sec:compactness} we prove the compactness for solutions of the conical toric Einstein equation with varying angles. The degree at $\beta=1/2$, corresponding to the standard Einstein metric on $S^2 \times S^2$, is computed in Section \ref{sec:half-angle}, and with this we also finish the proof of Theorem \ref{thm:intro-existence}.

\subsection{Related work of Chen--Dong, Li--Sun and Ozuch}
The results of this paper have been obtained independently by Li--Sun \cite{LS26} and Ozuch \cite{Ozu26}. Both papers establish the existence of $\TT$-invariant Einstein metrics on $\sharp_m(S^2 \times S^2)$ for every $m>1$. In addition, Li--Sun also prove the uniqueness statements of Theorem~\ref{thm:intro-uniqueness}, as well as  Theorem~\ref{thm:Hitchin-Thorpe}, showing among other things, that in fact  the refined Hitchin-Thorpe inequality is also sufficient for admitting toric Einstein metrics. Their construction combines local gluing with a global virtual counting argument. Ozuch constructs his two infinite families by varying the cone angle along a principal torus orbit as in our approach. His proof uses a continuity method, whereas we prove existence by a mod-two degree argument.  Also, uniqueness for $S^4$ has been recently obtained by Chen--Dong \cite{ChenDong2026}. 

\subsection{AI Declaration}
The existence results of this article were obtained with substantial use of ChatGPT 5.6 Sol in July-August 2026. In particular, we have used the model not only for checking and literature review, but also for exploring possible strategies, elaborate detailed arguments, and producing new useful identities, among other things. The subsequent write-up and proof-checking were aided by the use of ChatGPT 6 Astra and 6.1 Sol. Since the aforementioned papers appeared on the arXiv (October 6--8, 2026), we have not tried to correct or improve the mathematics in our paper in any way by the use of AI, and have focused solely on presentation issues.

\subsection{Acknowledgments}
The first-named author thanks Christoph Böhm and Wolfgang Ziller for raising questions that motivated the pursuit of the results in this paper. The second-named author thanks Song Sun for useful conversations on toric Einstein $4$-manifolds, and for his wonderful lectures on the topic delivered at the AMSI Winter School in Brisbane in July 2026. Both authors would also like to thank him and Tristan Ozuch for their positive and encouraging comments in regards to this submission and its significant overlap with their work.

\section{\texorpdfstring{$\TT$-invariant Einstein metrics on closed $4$-manifolds}{T2-actions on closed Einstein 4-manifolds}}\label{sec:general-equations}

\subsection{\texorpdfstring{Effective $\TT$-actions on closed $4$-manifolds}{Effective T2-actions on closed 4-manifolds}}\label{sec:orbit-topology}
\noindent Let $(M,\g)$ be a closed, simply-connected oriented $4$-manifold with $\operatorname{Ric}(\g)=\lambda\g$, admitting an effective, isometric $\TT$-action. We may assume $\lambda > 0$: indeed, for a nonzero Killing field $X$, the Bochner identity gives $\lambda\int_M|X|^2\dd v_\g=\int_M|\nabla X|^2\dd v_\g$, and equality would yield a nonzero parallel $1$-form, contradicting\textbf{} $H^1(M;\mathbb{R})=0$.

\phantomsection\label{notloc:orbit-space}
By \cite[Lemmas~5.1--5.2]{OrlikRaymond1}, the orbit space $\Sigma:=M/\TT$ is a topological closed disk. Its interior $\Sigma^\circ$ represents those orbits in $M$ where the action is free. The boundary $\partial\Sigma$ has $N:=\chi(M)$  vertices, corresponding to fixed points of the action, and $N$ open edges $E_1,\ldots,E_N$, indexed cyclically, corresponding to circle orbits in $M$. Under the quotient map $\pi:M\to\Sigma$, each closed edge has preimage $\mathcal D_i:=\pi^{-1}(\overline E_i)$, a smoothly embedded $2$-sphere, which we call a \emph{characteristic sphere}. Said sphere is a component of the fixed-point set of the circle stabilizer over $E_i$, and is therefore totally geodesic. 
We denote by $\mathcal D_i^2\in\mathbb Z$ the self-intersection number of $\mathcal D_i$ in $M$.
 The isotropy subgroup of $\mathcal{D}_i$ is a circle generated by a primitive vector $v_i\in\mathbb{Z}^2$ (a vector with coprime coordinates), unique up to sign. The generators at each vertex form a lattice basis: $|\det(v_i,v_{i+1})|=1$.

The Orlik--Raymond classification \cite{OrlikRaymond1} implies that the diffeomorphism type of such an $M$ is either $S^4$,  $\sharp_m(S^2\times S^2)$ for some $m\geq1$, or  $\sharp_p\mathbb{CP}^2\mathbin{\sharp}_q\overline{\mathbb{CP}}^2$, where $p,q\geq0$ and $p+q>0$. Since $\lambda>0$, the Hitchin--Thorpe inequality is strict: $3|\etau|<2\chi$; see \cite[pp.~161--162]{besse1987einstein}. On the last family, this reads $3|p-q|<2(2+p+q)$. Theorem~\ref{thm:toric-hitchin-thorpe} will improve this bound using the characteristic spheres.

\begin{example}
Two $\TT$-actions are equivalent if a diffeomorphism intertwines them after an automorphism of $\TT$. Up to this equivalence, $S^4$ and $\mathbb{CP}^2$ each admit a unique effective smooth $\TT$-action, whereas $S^2\times S^2$ and $\mathbb{CP}^2\mathbin{\sharp}\overline{\mathbb{CP}}^2$ each admit countably infinitely many; see \cite[\S\S~5.4--5.6]{OrlikRaymond1}. The latter actions are represented by the standard actions on the Hirzebruch surfaces $F_n=\mathbb P(\mathcal O\oplus\mathcal O_{\mathbb{CP}^1}(n))$, with $n\geq0$ even and odd, respectively. The following diagrams show the actions on $S^4$, $\mathbb{CP}^2$, $F_0$, and $F_1$; each edge is labelled by its stabilizer generator, up to sign.
\begin{center}
\resizebox{\linewidth}{!}{%
\begin{tikzpicture}[
    scale=0.85,
    boundary/.style={black, line width=0.55pt},
    normal/.style={->, black, line width=0.45pt},
    fixed/.style={circle, fill=black, inner sep=1.05pt},
    rodlabel/.style={font=\scriptsize, inner sep=1pt},
    title/.style={font=\small, anchor=base}
]

\begin{scope}
    \node[title] at (0,1.60) {$S^4$};

    \filldraw[boundary, fill=gray!18]
        (0,-0.70)
        .. controls (-0.77,-0.315) and (-0.77,0.315) .. (0,0.70)
        .. controls (0.77,0.315) and (0.77,-0.315) .. cycle;

    \foreach \p in {(0,-0.70),(0,0.70)}
        \node[fixed] at \p {};

    \node[rodlabel, left=3pt] at (-0.5775,0)
        {$[e_1]$};

    \node[rodlabel, right=3pt] at (0.5775,0)
        {$[e_2]$};
\end{scope}

\begin{scope}[shift={(4.15,0)}]
    \node[title] at (0,1.60) {$\mathbb{CP}^2$};

    \filldraw[boundary, fill=gray!18]
        (-0.70,-0.70) -- (0.70,-0.70)
        -- (-0.70,0.70) -- cycle;

    \foreach \p in {(-0.70,-0.70),(0.70,-0.70),(-0.70,0.70)}
        \node[fixed] at \p {};

    \draw[normal] (0,-0.70) -- (0,-0.33);
    \node[rodlabel, below=4pt] at (0,-0.70)
        {$[e_2]$};

    \draw[normal] (0,0) -- (-0.26,-0.26);
    \node[rodlabel, above right=3pt] at (0,0)
        {$[e_1+e_2]$};

    \draw[normal] (-0.70,0) -- (-0.33,0);
    \node[rodlabel, left=3pt] at (-0.70,0)
        {$[e_1]$};
\end{scope}

\begin{scope}[shift={(7.95,0)}]
    \node[title] at (0,1.60) {$S^2\times S^2$};

    \filldraw[boundary, fill=gray!18]
        (-0.70,-0.70) rectangle (0.70,0.70);

    \foreach \p in {(-0.70,-0.70),(0.70,-0.70),
                    (0.70,0.70),(-0.70,0.70)}
        \node[fixed] at \p {};

    \draw[normal] (0,-0.70) -- (0,-0.33);
    \node[rodlabel, below=4pt] at (0,-0.70)
        {$[e_2]$};

    \draw[normal] (0.70,0) -- (0.33,0);
    \node[rodlabel, right=3pt] at (0.70,0)
        {$[e_1]$};

    \draw[normal] (0,0.70) -- (0,0.33);
    \node[rodlabel, above=4pt] at (0,0.70)
        {$[e_2]$};

    \draw[normal] (-0.70,0) -- (-0.33,0);
    \node[rodlabel, left=3pt] at (-0.70,0)
        {$[e_1]$};
\end{scope}

\begin{scope}[shift={(12.65,0)}]
    \node[title] at (0,1.60)
        {$\mathbb{CP}^2\,\#\,\overline{\mathbb{CP}}^2$};

    \filldraw[boundary, fill=gray!18]
        (-1.05,-0.70) -- (1.05,-0.70)
        -- (-0.35,0.70) -- (-1.05,0.70) -- cycle;

    \foreach \p in {(-1.05,-0.70),(1.05,-0.70),
                    (-0.35,0.70),(-1.05,0.70)}
        \node[fixed] at \p {};

    \draw[normal] (-0.28,-0.70) -- (-0.28,-0.33);
    \node[rodlabel, below=4pt] at (-0.28,-0.70)
        {$[e_2]$};

    \draw[normal] (0.35,0) -- (0.09,-0.26);
    \node[rodlabel, above right=3pt] at (0.35,0)
        {$[e_1+e_2]$};

    \draw[normal] (-0.70,0.70) -- (-0.70,0.33);
    \node[rodlabel, above=4pt] at (-0.70,0.70)
        {$[e_2]$};

    \draw[normal] (-1.05,0) -- (-0.68,0);
    \node[rodlabel, left=3pt] at (-1.05,0)
        {$[e_1]$};
\end{scope}

\end{tikzpicture}%
}
\end{center}
\end{example}

\subsection{\texorpdfstring{The curvature operator of $\TT$-invariant Einstein metrics}{The curvature operator of T2-invariant Einstein metrics}}

\phantomsection\label{notloc:distributions}
\noindent By \cite[Theorem~1.2]{Liu2025ToricEinstein4Manifolds}, the distribution orthogonal to the principal orbits is integrable (i.e.~ the action is \emph{polar}). On $M^\circ:=\pi^{-1}(\Sigma^\circ)$, let $\Vca$ be the tangent distribution to the orbits and let $\Hca:=\Vca^\perp$. Choose circle generators $K_1,K_2 \in \mathfrak{t} = \Lie(\TT) \simeq \R^2$ with period $2\pi$. A local horizontal section gives
\begin{equation}\label{eqn:polar-metric}
 \g=\hat\g_{\alpha\beta}(x)\dd x^\alpha\dd x^\beta
      +G_{ij}(x)\dd\vartheta^i\dd\vartheta^j,
 \qquad \vartheta^i\in\mathbb R/2\pi\mathbb Z,
\textbf{}\end{equation}
where $\hat\g$ is the quotient metric on $\Sigma^\circ$ (the unique metric such that $\pi : M^\circ \to \Sigma^\circ$ is a Riemannian submersion) and $G_{ij}=\g(K_i,K_j)$. The local isometry $(x,\vartheta)\mapsto(x,2\vartheta_0-\vartheta)$ fixes the horizontal section $\vartheta=\vartheta_0$, so the section is totally geodesic. It also forces every curvature component with an odd number of vertical indices to vanish. Thus $\operatorname{Ric}(\Hca,\Vca)=0$, and the curvature operator preserves the plane spanned by the horizontal and vertical area forms.

\begin{notation}
\phantomsection\label{notloc:curvature-conventions}
Write $K$ for the Gauss curvature of $\hat\g$, $\hat\nabla$ for its Levi--Civita connection, and $\Delta_{\hat\g}=\operatorname{div}_{\hat\g}\hat\nabla$. Unless  otherwise specified, derivatives, inner products, and norms on $\Sigma^\circ$ are taken with respect to $\hat\g$. Our curvature convention is $\Rm (U,V)=\nabla_U\nabla_V-\nabla_V\nabla_U-\nabla_{[U,V]}$, and the curvature operator is $\mathcal R : \Lambda^2 T^ *M \to \Lambda^2 T^*M$,
\[
 \langle\mathcal R(U^\flat\wedge V^\flat),Z^\flat\wedge W^\flat\rangle
 =\g(\Rm(U,V)W,Z).
\]
On the oriented four-manifold $(M,\g)$, the Hodge star $\star$ gives the decomposition $\Lambda^2T^*M=\Lambda^2_+\oplus\Lambda^2_-$ into self-dual and anti-self-dual $2$-forms. We denote by $\mathcal W^\pm : \Lambda^2_\pm \to \Lambda^2_\pm$ the corresponding Weyl curvature endomorphisms. Since $\g$ is Einstein,
\[
    \mathcal R|_{\Lambda^2_\pm} =\mathcal W^\pm+\frac{\lambda}{3}\operatorname{Id}.
\]
Norms of $\mathcal W^\pm$ are always  Hilbert--Schmidt norms.
\end{notation}

\subsubsection{Reduction of the toric Einstein equations}

\phantomsection\label{notloc:gram-decomposition}
Set
\[
 \rho:=\sqrt{\det G},\qquad f:=\log\rho,\qquad \Phi:=\rho^{-1}G.
\]
Each principal orbit has area $4\pi^2\rho$.
\phantomsection\label{notloc:hyperbolic-plane}
Since $\Phi$ is positive definite with determinant one, it defines a map $\Phi:\Sigma^\circ\to \mathsf{SL}_2(\mathbb R) / \mathsf{SO}(2) \simeq \mathbb{H}^2$ to the hyperbolic plane, with the identification
\[
 \mathbb{H}^2=\{Z\in\operatorname{Sym}_2(\mathbb{R}):Z>0,\ \det Z=1\},
 \qquad
 (\g_{\mathbb{H}^2})_Z(A,B)=\tfrac12\operatorname{tr}(Z^{-1}AZ^{-1}B).
\]
We use the orientation
$\omega_{\mathbb{H}^2}=-Z_{11}^{-1}\dd Z_{11}\wedge\dd Z_{12}$.
We will also use the Klein disk compactification of $\mathbb H^2$: the map $Z\mapsto Z/\operatorname{tr}Z$ identifies $\mathbb H^2$ with the positive-definite symmetric matrices of trace one, whose closure is $\overline{\mathbb H}^2:=\{Y\in\operatorname{Sym}_2(\mathbb R):Y\geq0,\ \operatorname{tr}Y=1\}$. Its ideal boundary $\partial_\infty\mathbb H^2$ consists precisely of the rank-one matrices. In this model, hyperbolic geodesics are Euclidean line segments. 
For each nonzero $v\in\mathfrak t$,  we let $P_v\in\partial_\infty\mathbb H^2$ denote the unique rank-one positive semi-definite symmetric matrix with $\operatorname{tr}P_v=1$ and $\ker P_v=\langle v\rangle$. In particular, $P_v=P_{cv}$ for every $c\neq0$.

For a function $v\in C^\infty(\Sigma^\circ)$ and the map $\Phi$, we consider the weighted Laplacian and tension field:
\begin{equation}\label{eqn:weighted-operators}
 \mathbf L_fv:=\Delta v+\langle\dd f,\dd v\rangle
     =\rho^{-1}\operatorname{div}(\rho\hat\nabla v),
 \qquad
 \etau_f(\Phi):=\operatorname{tr}_{\hat\g}\nabla\dd\Phi
                         +\dd\Phi(\hat\nabla f).
\end{equation}
For invariant functions $w = \pi^* v$ on $M$, the total space Laplacian is $ \Delta_{\g}(w) = \pi^* \mathbf L_f (v)$. 

The Ricci formulas in~\cite[\S~4.2, equation~(4.6), p.~501]{Lott2010}, applied to~\eqref{eqn:polar-metric} with $G=e^f\Phi$, give the following reduced system.

\begin{proposition}\label{prop:quotient-equations}
A metric $\g$ as in~\eqref{eqn:polar-metric} is Einstein with constant
$\lambda$ if and only if, on $\Sigma^\circ$,
\begin{align}
 \Delta f+|\dd f|^2&=-2\lambda,
       \qquad \qquad \etau_f(\Phi)=0,\label{eqn:general-vertical}\\
 \hat\nabla^2f&=(K-\lambda)\hat\g
       -\tfrac12\dd f\otimes\dd f-\tfrac12\Phi^*\g_{\mathbb{H}^2}.
       \label{eqn:general-horizontal}
\end{align}
In this case, the following also hold:
\begin{equation}\label{eqn:general-trace}
 \Delta\rho=-2\lambda\rho,
 \qquad 4K=|\dd\Phi|^2-|\dd f|^2.
\end{equation}
\end{proposition}

\begin{remark}\label{rem:diagonal-reduction}
If $K_1\perp K_2$ on $\Sigma^\circ$, then $\Phi=\operatorname{diag}(e^u,e^{-u})$,
where $u=\log(|K_1|/|K_2|)$ is the arclength parameter along a geodesic of
$\mathbb{H}^2$. Consequently,
$\Phi^*\g_{\mathbb{H}^2}=\dd u\otimes\dd u$ and
$\etau_f(\Phi)=0$ is equivalent to $\mathbf L_fu=0$.
Proposition~\ref{prop:quotient-equations} is local: its equations also
hold for two commuting Killing fields with integrable common
orthogonal distribution, wherever they are independent.
Periodicity is used only in the orbit-area normalization.
\end{remark}

\subsubsection{Smoothness conditions at singular orbits}
At the singular orbits, corresponding to points on the edges or vertices of $\partial\Sigma$, there are a number of conditions to be satisfied for a metric on the regular locus $M^\circ$ to extend to a smooth metric on $M$.   
These are contained in the following lemma, which is essentially \cite[Lemmas~2.1--2.2]{Liu2025ToricEinstein4Manifolds} in the particular case of a polar action:

\begin{lemma}\cite{Liu2025ToricEinstein4Manifolds}\label{lem:slice-regularity}
Near an open edge $E \subset \partial \Sigma$, choose a primitive generator $K_1$ that vanishes there and complete it to a lattice basis $K_1,K_2$. With respect to the inward distance $r$ and edge arclength $s$, we may write
\begin{equation}\label{eqn:edge-slice-expansion}
 \hat\g=\dd r^2+h(r,s)^2\dd s^2,\qquad
 G=\begin{pmatrix}r^2A&r^2B\\r^2B&C\end{pmatrix},\qquad
 \rho=rD,
\end{equation}
where the coefficients are smooth and even in $r$, with
$D^2=AC-r^2B^2>0$ and $h(0,s)=A(0,s)=1$.
In particular, $\partial_\nu\rho=|K_2|$ along the edge, where $\nu$ is the inward unit normal.

At a vertex, the adjacent stabilizers give a lattice basis and equivariant normal coordinates give a section $r_1,r_2\geq0$ in which
\begin{equation}\label{eqn:vertex-slice-expansion}
\hat\g=h_{11}\dd r_1^2
+2r_1r_2h_{12}\dd r_1\dd r_2
+h_{22}\dd r_2^2, \qquad
 G=\begin{pmatrix}
 r_1^2A_{11}&r_1^2r_2^2A_{12}\\
 r_1^2r_2^2A_{12}&r_2^2A_{22}
 \end{pmatrix},\quad
 \rho=r_1r_2D,
\end{equation}
where all coefficients $h_{ij}$ and $A_{ij}$ are smooth
and even in each variable, with
\[
A_{11}(0,0)=A_{22}(0,0)=1,\qquad
h_{11}(0,r_2)=A_{11}(0,r_2),\qquad
h_{22}(r_1,0)=A_{22}(r_1,0).
\]
Here, 
\[
D=\sqrt{A_{11}A_{22}-r_1^2r_2^2A_{12}^2}>0.
\]
In particular, the edges are geodesic, the
corner angles are $\pi/2$, and $K$ extends
smoothly and evenly to the reflected sections.
\end{lemma}
\begin{proof}
The expansions and smoothness conditions follow directly from
\cite[Lemmas~2.1--2.2]{Liu2025ToricEinstein4Manifolds}. The identity
$\partial_\nu\rho=|K_2|$ follows from
$\rho=rD$ and $D(0,s)=\sqrt{C(0,s)}$.
The remaining conclusions follow from the
reflection symmetries of the local sections.
\end{proof}


\begin{remark}
Since the $\TT$-action is polar, by \cite{HLO} there exists a complete immersed section $S$, meeting every orbit orthogonally. Geometrically, one can think of $S$ as obtained by reflecting $\Sigma$ across its boundary edges. For example, when $M=S^4$, the orbit space $\Sigma$ is a bigon with angles $\pi/2$, and one may take $S\simeq S^2$ to be embedded. The polar group $\Pi$, consisting of the elements of $\TT$ preserving $S$, is then isomorphic to $\mathbb Z_2\times\mathbb Z_2$, and $\Sigma\simeq S/\Pi$.
\end{remark}

We say that a function $u$ on $\Sigma^\circ$ \emph{extends smoothly and evenly to the reflected sections} if, near each open edge, it extends smoothly to signed normal coordinates $(r,s)$, with $u(-r,s)=u(r,s)$, and near each vertex, it extends smoothly to signed coordinates $(r_1,r_2)$, with $u(-r_1,r_2)=u(r_1,r_2)=u(r_1,-r_2)$. This is a local condition, independent of the choice of a global section.

\subsubsection{Eigenvalues of the Weyl curvature}
\phantomsection\label{notloc:distinguished-forms}
We now endow $\Sigma^\circ$ with an orientation so that its orientation followed by that of $(K_1,K_2)$ agrees with the orientation of $M$. Let $\star$ denote the Hodge star operator of $(M,g)$, and let $\dd \mu_{\hat g}$ be the Riemannian volume form of $\hat g$. Set
\[
 \Phi^*\omega_{\mathbb{H}^2}
   =\operatorname{Jac}_\Phi\,\dd \vol_{\hat\g},\qquad
 \omega_\Hca:=\pi^*\dd\vol_{\hat\g},\qquad
 \omega_\Vca:=\star\omega_\Hca
\]
For $\varepsilon\in\{+1,-1\}$, put
$\Omega_\varepsilon:=2^{-1/2}(\omega_\Hca +\varepsilon\omega_\Vca)$.

\begin{lemma}\label{lem:distinguished-weyl}
For $\varepsilon\in\{+1,-1\}$, the unit $2$-form $\Omega_\varepsilon$ is an eigenform of $\mathcal W^\varepsilon$, with eigenvalue
\begin{equation}\label{eqn:distinguished-weyl}
 \mu_\varepsilon=K+\tfrac\varepsilon2\operatorname{Jac}_\Phi-\tfrac\lambda3.
\end{equation}
Both $\operatorname{Jac}_\Phi$ and $\mu_\varepsilon$ extend smoothly and evenly to the reflected sections. Moreover,
\begin{equation}\label{eqn:weyl-Q-relation}
 Q_\varepsilon:=|\dd\Phi|^2+2\varepsilon\operatorname{Jac}_\Phi
       =4\mu_\varepsilon+\tfrac{4\lambda}3+|\dd f|^2\geq0.
\end{equation}
\end{lemma}
\begin{proof}
On the plane spanned by $\omega_\Hca,\omega_\Vca$, the curvature operator has both diagonal entries equal to $K$: horizontal sections are totally geodesic, and $[\mathcal{R},\star]=0$  by the Einstein condition. To compute its off-diagonal entry $c$, write $\Rm_{\alpha\beta ij}:=\g(\Rm(\partial_\alpha,\partial_\beta)K_j,K_i)$. The vertical connection gives, with $S_\alpha:=G^{-1}\partial_\alpha G$ and $U_\alpha:=\Phi^{-1}\partial_\alpha\Phi$ for this calculation,
\[
 (\operatorname R_{\alpha\beta ij})_{ij}
       =-\tfrac14G[S_\alpha,S_\beta]
       =-\tfrac\rho4\Phi[U_\alpha,U_\beta].
\]
Choose a constant determinant-one vertical frame in which $\Phi=\operatorname{Id}$ at the point, and a positively-oriented orthonormal horizontal frame. Write the diagonal entries of $\partial_\alpha\Phi$ as $u_\alpha,-u_\alpha$ and its off-diagonal entry as $v_\alpha$. Then $\operatorname{Jac}_\Phi=-(u_1v_2-u_2v_1)$, while the commutator has $(1,2)$-entry $2(u_1v_2-u_2v_1)$. Dividing by the area density $\rho$ of the vertical frame gives $c=\tfrac12\operatorname{Jac}_\Phi$. The eigenforms of this block are $\Omega_\varepsilon$. Subtracting $\lambda/3$ gives~\eqref{eqn:distinguished-weyl}.

To establish regularity at the singular orbits, observe that
\[
\Phi^*\omega_{\mathbb H^2}=-\dd\left(\frac{G_{11}}{\rho}\right)\wedge\dd\left(\frac{G_{12}}{G_{11}}\right).
\]
Near an open edge, the two factors are $rA/D$ and $B/A$, so the wedge product of their differentials is smooth and even in $r$. At a vertex, they are $r_1a/r_2$ and $r_2^2b$, where $a=A_{11}/D$ and $b=A_{12}/A_{11}$ are smooth and even in both variables. Direct differentiation shows that
\[
\dd(r_1a/r_2)\wedge\dd(r_2^2b)=\big[(a+r_1\partial_1a)(2b+r_2\partial_2b)-r_1(r_2\partial_2a-a)\partial_1b\big]\,\dd r_1\wedge\dd r_2,
\]
whose coefficient is again smooth and even in both variables. Since the quotient volume density has the same regularity, so does $\operatorname{Jac}_\Phi$. The assertion for $\mu_\varepsilon$ now follows from Lemma~\ref{lem:slice-regularity}.


The equality in~\eqref{eqn:weyl-Q-relation} follows from~\eqref{eqn:general-trace}. If $J$ is positive rotation by $\pi/2$ on $\mathbb{H}^2$, then
\[
 Q_\varepsilon
 =|\dd\Phi(e_1)-\varepsilon J\dd\Phi(e_2)|^2
\]
for a positively-oriented orthonormal horizontal frame, proving nonnegativity.
\end{proof}

\begin{remark}
Geometrically, $Q_\varepsilon$ measures the failure of $\Phi$ to be conformal with the corresponding orientation. Equivalently, for the canonical almost-K\"ahler structure $\omega_\varepsilon=\rho^{-1}(\omega_{\mathcal H}+\varepsilon\omega_{\mathcal V})$ on $(M^\circ,\widetilde\g=\rho^{-1}\g)$, $Q_\varepsilon$ measures its failure from being K\"ahler:
\[
|\widetilde\nabla\omega_\varepsilon|_{\widetilde\g}^2=\rho Q_\varepsilon.
\]
\end{remark}

The following Weyl curvature estimate in Proposition \ref{prop:weyl-sign} plays a central role in the paper. It underlies the improved Hitchin--Thorpe inequality (Theorem~\ref{thm:Hitchin-Thorpe}), the uniqueness results in Section~\ref{sec:uniq}, and analogous estimates in the conical setting used in the compactness argument for Theorem~\ref{thm:intro-existence}. We begin with two differential identities for the weighted distinguished Weyl eigenvalues.

\begin{lemma}\label{lem:weyl-differential-identities}
For $\varepsilon\in\{\pm1\}$, set $F_\varepsilon:=\rho^{3/2}\mu_\varepsilon$. \label{notloc:weighted-weyl} Then, on $\Sigma^\circ$,
\begin{align}
\Delta F_\varepsilon&=\left(\lambda+\tfrac34|\dd f|^2\right)F_\varepsilon+2\rho^{3/2}|\mathcal W^\varepsilon|^2,\label{eqn:weyl-smooth-scalar}\\
|\dd F_\varepsilon|^2&=\tfrac12\rho^3Q_\varepsilon\left(|\mathcal W^\varepsilon|^2-\tfrac32\mu_\varepsilon^2\right).\label{eqn:weyl-gap-gradient}
\end{align}
\end{lemma}
\begin{proof}
Fix $\varepsilon$ and write $F=F_\varepsilon$, $\mu=\mu_\varepsilon$, and $Q=Q_\varepsilon$. Set $\tilde\g=\rho^{-1}\g$, write $\tilde s$ for its scalar curvature, and put
\[
 \omega=\rho^{-1}(\omega_\Hca+\varepsilon\omega_\Vca)
       =\rho^{-1}\pi^*\dd v_{\hat\g}
          +\varepsilon\dd\vartheta^1\wedge\dd\vartheta^2.
\]
The form $\omega$ is closed and has $\tilde\g$-norm $\sqrt2$. With the orientation determined by $\omega$, the metric $\tilde\g$ is almost K\"ahler. Put $A=\rho^{1/2}\mathcal W^\varepsilon_{\tilde\g}$. Conformal covariance of the Einstein Bianchi identity gives $\delta_{\tilde\g}A=0$, and $A\omega=F\omega$. Set $P=\omega\otimes\omega$. Since $\langle A,\tilde\nabla P\rangle=0$ and $\langle\tilde\nabla A,P\rangle=2\dd F$, the product rule and the Weitzenb\"ock identities in \cite[\S~2, equation~(2) and Lemma~1]{LeBrun2015HarmonicForms} give
\begin{align*}
 -2\Delta_{\tilde\g}F
 &=\langle\tilde\nabla^*\tilde\nabla A,P\rangle
       -\langle A,\tilde\nabla^*\tilde\nabla P\rangle\\
 &=\bigl(-\tilde sF+12\rho\mu F-4\rho^{-1/2}|A|^2\bigr)
       -\bigl(12\rho\mu F-2\tilde sF\bigr)\\
 &=\tilde sF-4\rho^{-1/2}|A|^2.
\end{align*}
The minus sign in the product rule follows by differentiating $\langle A,\tilde\nabla P\rangle=0$, which gives $\langle\tilde\nabla A,\tilde\nabla P\rangle =\langle A,\tilde\nabla^*\tilde\nabla P\rangle$. 
The scalar curvature formula and~\eqref{eqn:general-vertical} give $\tilde s=-\rho(2\lambda+\tfrac32|\dd f|^2)$. The conformal Laplacian formula and torus invariance give $\Delta_{\tilde\g}F=\rho\Delta F$; since $|A|^2=\rho^3|\mathcal W^\varepsilon|^2$, this proves~\eqref{eqn:weyl-smooth-scalar}. The harmonic-form identity also gives $|\tilde\nabla\omega|^2=4\rho\mu-2\tilde s/3=\rho Q$.

At a point, complete $\omega$ to an orthogonal frame $\omega,\phi,\psi$ of self-dual forms of norm $\sqrt2$ in the chosen orientation, diagonalizing $A|_{\omega^\perp}$ there, with eigenvalues $a,b$. Write $\tilde\nabla\omega=\alpha\otimes\phi+\beta\otimes\psi$ at that point. Closure of $\omega$ and the projection of $\delta_{\tilde\g}A=0$ onto $\omega$ give
\[
 \alpha\wedge\phi+\beta\wedge\psi=0,
 \qquad
 \dd F\wedge\omega=(a-b)\alpha\wedge\phi.
\]
For the second identity, regard $A$ as a self-dual bundle-valued $2$-form, so that $\delta A=0$ is equivalent to $\dd^{\tilde\nabla}A=0$, and apply the product rule to $\dd(A\omega)=\dd(F\omega)$. Wedge multiplication by any of $\omega,\phi,\psi$ is an isometry on $1$-forms. Thus $|\alpha|=|\beta|$ and $|\tilde\nabla\omega|^2=4|\alpha|^2$. Trace-freeness gives $a+b=-F$, so
\[
 |\dd F|_{\tilde\g}^2
 =\tfrac14(a-b)^2|\tilde\nabla\omega|^2
 =\tfrac14(2|A|^2-3F^2)\rho Q.
\]
Since $|\dd F|_{\tilde\g}^2=\rho|\dd F|^2$ and $|A|^2=\rho^3|\mathcal W^\varepsilon|^2$, this proves~\eqref{eqn:weyl-gap-gradient}, including at $Q=0$.
\end{proof}

\begin{proposition}\label{prop:weyl-sign}
For $\varepsilon\in\{+1,-1\}$, either $\mu_\varepsilon<0$ on the closed disk $\Sigma$, or $\mathcal W^\varepsilon\equiv0$ on $M$. In particular,
\begin{equation}\label{eqn:estimateKJacPhi}
 K+\tfrac12|\operatorname{Jac}_\Phi|\leq\tfrac\lambda3.
\end{equation}
Also, if $Q_\varepsilon>0$ throughout $\Sigma^\circ$, then $\mathcal W^\varepsilon\equiv0$.
\end{proposition}
\begin{proof}
Fix $\varepsilon\in\{\pm1\}$, and write $\mu=\mu_\varepsilon$, $Q=Q_\varepsilon$ and $F=F_\varepsilon$. By Lemmas~\ref{lem:slice-regularity} and~\ref{lem:distinguished-weyl}, $F$ extends continuously by zero to $\partial\Sigma$. The first identity in Lemma~\ref{lem:weyl-differential-identities} rules out a positive interior maximum of $F$, so $F\leq0$. The strong maximum principle gives $F<0$ or $F\equiv0$ in $\Sigma^\circ$. In the latter case, the same identity implies $\mathcal W^\varepsilon\equiv0$.

To include the boundary in the strict alternative, expand~\eqref{eqn:weyl-smooth-scalar} using~\eqref{eqn:general-vertical}:
\begin{equation}\label{eqn:weyl-full-block}
\mathbf L_{3f}\mu = \Delta\mu + 3\langle\dd f,\dd\mu\rangle =4\lambda\mu+2|\mathcal W^\varepsilon|^2.
\end{equation}
Hence $w=-\mu\geq0$ satisfies $(\mathbf L_{3f}-4\lambda)w\leq0$. At a corner, Lemma~\ref{lem:slice-regularity} gives $\rho=r_1r_2D$, with $D$ smooth, positive and even. The diagonal terms of $\mathbf L_{3f}$ are $\hat\g^{jj}(\partial_{r_j}^2+3r_j^{-1}\partial_{r_j})$; the mixed second-order coefficient is $r_1r_2$ times a smooth even function, and the remaining first-order coefficients are $r_j$ times such functions. Consequently, replacing $r_j$ by $|x_j|$, $x_j\in\mathbb R^4$, yields a differential operator with smooth coefficients, since $\partial_{r_j}^2+3r_j^{-1}\partial_{r_j}$ becomes the Euclidean Laplacian on $\mathbb R^4$, and the factors $r_j$ in the other terms absorb the radial derivatives. The lifted principal part is the Euclidean Laplacian at the corner, hence uniformly elliptic nearby. The smooth even function $w$ also lifts smoothly. At an open edge, the same argument applies using one radial variable and the tangential coordinate. A boundary zero is therefore an interior minimum of a nonnegative solution of a uniformly elliptic inequality. The strong maximum principle~\cite[Theorem~3.5]{GilbargTrudinger} forces $w=0$ near that point and then throughout the connected interior. This proves the strict alternative on $\Sigma$, and~\eqref{eqn:distinguished-weyl} gives the curvature bound.

Finally, suppose $Q>0$ throughout $\Sigma^\circ$. If $F<0$, then $F$ attains a negative minimum at an interior point. There the gradient identity in Lemma~\ref{lem:weyl-differential-identities} gives $|\mathcal W^\varepsilon|^2=\tfrac32\mu^2$. Substituting into the first identity and using~\eqref{eqn:weyl-Q-relation}, we obtain $\Delta F=\tfrac34QF<0$, contradicting minimality. Hence $F\equiv0$ and $\mathcal W^\varepsilon\equiv0$.
\end{proof}

\begin{remark}\label{rem:weyl-applications}
For the standard Einstein metric  on $S^2\times S^2$, one has $\mu_\pm=-\lambda/3$ and $Q_+=Q_-=|\dd u|^2$, which vanishes at the equatorial torus. Since $\mathcal W^\pm\not\equiv0$, the hypothesis $Q_\varepsilon>0$ in Proposition~\ref{prop:weyl-sign} cannot be omitted.

More generally, if $K=\lambda/3$ at an interior point of the orbit space of a toric Einstein four-manifold, \eqref{eqn:estimateKJacPhi} forces $\operatorname{Jac}_\Phi=0$ there. Both distinguished eigenvalues vanish, so Proposition~\ref{prop:weyl-sign} gives $\mathcal W^+=\mathcal W^-=0$, and hence constant sectional curvature $\lambda/3$.
\end{remark}

\subsection{\texorpdfstring{An improved Hitchin--Thorpe inequality for $\TT$-invariant Einstein metrics}{A Hitchin-Thorpe inequality for T2-invariant Einstein metrics}}

\begin{theorem}\label{thm:toric-hitchin-thorpe}
A closed, simply-connected Einstein $4$-manifold admitting an effective isometric $\TT$-action satisfies
\[
3|\etau(M)|\leq\chi(M).
\]
\end{theorem}

\begin{proof}
Let $\kappa_i$ and $\kappa_i^\perp$ denote the tangent and normal curvatures of $\mathcal D_i$. By Gauss--Bonnet and the Chern--Weil formula for the Euler class of its normal bundle, respectively,
\[
\int_{\mathcal D_i}\kappa_i\,\dd A=4\pi,\qquad \int_{\mathcal D_i}\kappa_i^\perp\,\dd A=2\pi\mathcal D_i^2.
\]

At a point of $\mathcal D_i$, choose an oriented orthonormal frame $e_1,e_2,e_3,e_4$, with $e_1,e_2$ tangent and $e_3,e_4$ normal to $\mathcal D_i$. Writing $e^{ab}=e^a\wedge e^b$, an orthonormal basis of $\Lambda^2_\pm$ is
\[
\eta_0^\pm=\tfrac1{\sqrt2}(e^{12}\pm e^{34}),\qquad \eta_1^\pm=\tfrac1{\sqrt2}(e^{13}\mp e^{24}),\qquad \eta_2^\pm=\tfrac1{\sqrt2}(e^{14}\pm e^{23}).
\]
Since $\mathcal D_i$ is totally geodesic and $\g$ is Einstein, the Gauss and Ricci equations give
\[
\langle\mathcal R e^{12},e^{12}\rangle=\langle\mathcal R e^{34},e^{34}\rangle=\kappa_i,\qquad \langle\mathcal R e^{12},e^{34}\rangle=\kappa_i^\perp.
\]
The stabilizer circle fixes $\eta_0^\pm$ and rotates the plane spanned by $\eta_1^\pm,\eta_2^\pm$. Since $\mathcal W^\pm$ is symmetric, trace-free and commutes with these rotations, its matrix in this basis is
\[
\mathcal W^\pm=\operatorname{diag}(a_\pm,-a_\pm/2,-a_\pm/2),\qquad a_\pm=\kappa_i\pm\kappa_i^\perp-\lambda/3.
\]

In the adapted frames of Lemma~\ref{lem:distinguished-weyl}, the limiting horizontal and vertical planes each contain one tangent and one normal direction to $\mathcal D_i$. Thus $\Omega_\pm$ belongs to the repeated eigenspace, and
\begin{equation}\label{eqn:weyl-edge-limit}
\mu_\pm=\frac{\lambda}{6}-\frac12(\kappa_i\pm\kappa_i^\perp)\qquad\text{along }\mathcal D_i.
\end{equation}
Proposition~\ref{prop:weyl-sign} therefore yields $\kappa_i\pm\kappa_i^\perp\geq\lambda/3$. Integrating gives
\[
\frac\lambda3\operatorname{Area}(\mathcal D_i)\leq4\pi-2\pi|\mathcal D_i^2|.
\]
Since $\lambda>0$, we conclude that $|\mathcal D_i^2|<2$, and hence $\mathcal D_i^2\in\{-1,0,1\}$.

If $N:=\chi(M)=2$, then $b_2(M)=0$ by simple connectedness, so $\etau(M)=0$. Otherwise $N\geq3$, and the quotient is a simple $N$-gon. The stable tangent splitting of \cite[Theorem~3.16 and its proof]{BR01}, together with the Hirzebruch signature theorem, gives $3\etau(M)=\sum_i\mathcal D_i^2$. Consequently,
\[
3|\etau(M)|=\left|\sum_i\mathcal D_i^2\right|\leq\sum_i|\mathcal D_i^2|\leq N=\chi(M),
\]
as claimed.
\end{proof}

\begin{corollary}\label{cor:connected-sum-obstruction}
If $\sharp_p\mathbb{CP}^2\mathbin{\sharp}_q\overline{\mathbb{CP}}^2$ admits a $\TT$-invariant Einstein metric, then $p\leq2q+1$ and $q\leq2p+1$.
\end{corollary}

\section{Uniqueness results}\label{sec:uniq}

Our goal in this section is to establish the uniqueness results in Theorem~\ref{thm:intro-uniqueness}. For $S^4$ and $\mathbb{CP}^2$, we exploit the trace-normalised Gram map $\Psi$, its boundary extension, and a non-degeneracy criterion for solutions of planar elliptic equations with monotone boundary values. For $S^2\times S^2$ and $\mathbb{CP}^2\mathbin{\sharp}\overline{\mathbb{CP}}^2$, we combine Gauss--Bonnet with Wu's maximum principle \cite{Wu2021} to obtain repeated Weyl eigenvalues, and then apply the classification results of Derdzi\'nski \cite{Derdzinski1983SelfDualKahler} and LeBrun \cite{LeBrun2012EinsteinHermitian} on Einstein Hermitian four-manifolds.

Fix a lattice basis for the torus action and retain the notation $G$, $\rho=\sqrt{\det G}$, $f=\log\rho$, and $\Phi=\rho^{-1}G$. \phantomsection\label{notloc:resolved-vertices}
Let $\Sigma^{\mathrm v}$ be the disk obtained by replacing each vertex of $\Sigma$ by its interval of inward directions. In order to understand the boundary behavior of $\Phi : \Sigma^\circ \to \mathbb H^2$, it is more natural to work with the Klein disk compactification $\overline{\mathbb{H}}^2$ (see $\S$\ref{notloc:gram-decomposition}). The corresponding normalisation of the Gram matrix $G$ satisfies the following equation:

\begin{lemma}\label{lem:normalized-gram-map}
Set $\Psi:=G/\operatorname{tr}G$. Then
\begin{equation}\label{eq:normalized-gram-equation}
 \Delta_{\hat\g}\Psi+
 \left\langle2\dd\log(\operatorname{tr}G)-\dd f,\dd\Psi\right\rangle_{\hat\g}=0
\end{equation}
entrywise on $\Sigma^\circ$. The map $\Psi$ extends continuously to $\Sigma^{\mathrm v}$, taking the constant value $P_{v_i}$ along each original edge $E_i$. Each added interval maps homeomorphically onto the segment joining the matrices associated with its adjacent edges.
\end{lemma}
\begin{proof}
Choose a local $\hat\g$-orthonormal frame $e_1,e_2$, and set $B_j:=\Phi^{-1}\dd\Phi(e_j)$. Differentiating $\det\Phi=1$ gives $\operatorname{tr}B_j=0$. A trace-free $2\times2$ matrix satisfies $B_j^2=\tfrac12\operatorname{tr}(B_j^2)\operatorname{Id}$. The weighted harmonic-map equation~\eqref{eqn:general-vertical} therefore gives
\[
 \Phi^{-1}\mathbf L_f\Phi=\sum_{j=1}^2B_j^2=|\dd\Phi|^2\operatorname{Id}.
\]
Taking traces shows that $\operatorname{tr}\Phi$ satisfies the same scalar equation as each entry of $\Phi$. Substituting $\Phi=(\operatorname{tr}\Phi)\Psi$ and applying the product rule gives $\mathbf L_f\Psi+2\langle\dd\log(\operatorname{tr}\Phi),\dd\Psi\rangle=0$. Since $\operatorname{tr}\Phi=e^{-f}\operatorname{tr}G$, this is \eqref{eq:normalized-gram-equation}.

Along an open edge $E$, the matrix $G$ has rank one and kernel $\langle v_E\rangle$. Its trace is positive, so $\Psi=P_{v_E}$. At a vertex, take equivariant geodesic normal coordinates and use the adjacent isotropy generators as a lattice basis. If $r_1,r_2$ are the slice radii and $R^2=r_1^2+r_2^2$, the metric differs from its value at the vertex by $O(R^2)$, while the circle generators are linear rotation fields. Their Gram matrix is therefore $\operatorname{diag}(r_1^2,r_2^2)+O(R^4)$. Returning to the fixed lattice basis gives $G=r_1^2Q_1+r_2^2Q_2+O(R^4)$, where $Q_1,Q_2$ are fixed rank-one nonnegative matrices with distinct kernels. Consequently,
\[
 \Psi=\frac{r_1^2Q_1+r_2^2Q_2}
 {r_1^2\operatorname{tr}Q_1+r_2^2\operatorname{tr}Q_2}+O(R^2).
\]
The $O(R^2)$ estimate is uniform in the direction of approach. The leading matrix is a convex combination of $Q_1/\operatorname{tr}Q_1$ and $Q_2/\operatorname{tr}Q_2$. Its first coefficient increases strictly from zero to one as $r_1/r_2$ increases from zero to infinity. Thus the added interval traverses the segment once, and its endpoints agree with the adjacent edge values.
\end{proof}

The following lemma is a non-degeneracy criterion for solutions of planar elliptic equations: under suitable monotonicity of the boundary values, we can rule out interior critical points. Applied  to the extension of $\Psi$  to $\Sigma^{\mathrm v}$, it  yields a function without critical points on the bigon orbit space of $S^4$, and an immersion into $\mathbb H^2$ in the case of $\mathbb{CP}^2$.

\begin{lemma}[Alessandrini--Sigalotti]\label{lem:boundary-monotonicity}\label{lem:two-boundary-sign-arcs}
Let $D\subset\mathbb{R}^2$ be the bounded region enclosed by a simple closed curve. Let $h$ be a smooth Riemannian metric and $X$ a smooth vector field on $D$. Suppose that $u\in\mathcal C^0(\overline D)\cap\mathcal C^\infty(D)$ is nonconstant and satisfies $\Delta_hu+X(u)=0$. If $p_-,p_+\in\partial D$ are distinct and $u$ is nondecreasing along each boundary arc traversed from $p_-$ to $p_+$, then $\dd u\neq0$ throughout $D$.
\end{lemma}
This is the $N=1$ case of \cite[Theorem~2.3, p.~253, and Proposition~2.7, p.~255]{AlessandriniSigalotti2001}.
\begin{proof}
Suppose that $\dd u(p)=0$, and set $t=u(p)$. The strong maximum principle gives $\min_{\partial D}u<t<\max_{\partial D}u$. Boundary monotonicity therefore makes $B:=\{q\in\partial D:u(q)>t\}$ a nonempty connected arc. Set $U:=\{q\in D:u(q)>t\}$. Every component of $U$ has closure meeting $B$; otherwise $u$ would attain a maximum greater than $t$ in the interior. Each point of $B$ has a neighborhood whose intersection with $D$ is connected and contained in $U$. The component containing that intersection is locally constant along $B$, so $U$ is connected.

The Hartman--Wintner theorem gives at least four local branches of $u=t$ through $p$, separating sectors on which $u-t$ has alternating signs; see \cite[Theorem~(H.--W.) and Remark~1.2, p.~231]{Alessandrini1987}. Since $U$ is connected, a path in $U$ joins two distinct positive sectors. A portion of that path joins distinct sectors outside a small disk about $p$. Removing loops and joining the endpoints to $p$ within the sectors produces a simple closed curve. It bounds a disk compactly contained in $D$, with $u\geq t$ on its boundary and, by the alternating signs near $p$, $u<t$ somewhere inside. Thus $u$ attains an interior minimum, contradicting the strong maximum principle.
\end{proof}

\subsection{\texorpdfstring{The round metric on $S^4$}{The round metric on S4}}\label{sec:round-s4}
\begin{theorem}\label{thm:round-s4-uniqueness}
Every Einstein metric on $S^4$ admitting an effective isometric $\TT$-action is homothetic to the round metric.
\end{theorem}

\begin{proof}
Since $\chi(S^4)=2$, the orbit disk has two edges, whose primitive isotropy generators $K_1, K_2$ form a lattice basis. By Lemma~\ref{lem:normalized-gram-map}, the boundary image of $\Psi$ on the vertex-resolved disk $\Sigma^{\mathrm v}$ is the segment joining $\operatorname{diag}(1,0)$ and $\operatorname{diag}(0,1)$ in $\overline{\mathbb{H}}^2$. In particular, $\Psi_{12}$ vanishes on $\partial\Sigma^{\mathrm v}$, so the maximum principle applied to~\eqref{eq:normalized-gram-equation} gives $\Psi_{12}=0$ throughout $\Sigma^\circ$ (cf.~\cite[Proposition~5]{Liu2025ToricEinstein4Manifolds}).

Write $\Psi=\operatorname{diag}(h,1-h)$, where $0<h<1$. The function $h$ takes the values $0$ and $1$ on the two original edges, while at each added interval (of inward directions at vertices) it traverses $[0,1]$ monotonically. Thus the boundary hypothesis of Lemma~\ref{lem:boundary-monotonicity} is satisfied, giving $\dd h\neq0$ throughout $\Sigma^\circ$.

Set $u:=\log(|K_1|/|K_2|)=\tfrac12\log(h/(1-h))$. Since $\Phi=\operatorname{diag}(e^u,e^{-u})$, $|\dd\Phi|^2=|\dd u|^2$ and Remark~\ref{rem:diagonal-reduction} gives $\operatorname{Jac}_\Phi=0$, from which
\[
Q_+=Q_-=|\dd u|^2=\frac{|\dd h|^2}{4h^2(1-h)^2}>0.
\]
Proposition~\ref{prop:weyl-sign} therefore yields $\mathcal W^\pm=0$. Hence $\g$ has constant sectional curvature $\lambda/3>0$. By completeness and simple connectedness, $(S^4,\g)$ is a round sphere of radius $\sqrt{3/\lambda}$.
\end{proof}

\subsection{\texorpdfstring{The Fubini--Study metric on $\mathbb{CP}^2$}{The Fubini-Study metric on CP2}}\label{sec:fubini-study}
\begin{theorem}\label{thm:fubini-study-uniqueness}
Every Einstein metric on $\mathbb{CP}^2$ admitting an effective isometric $\TT$-action is homothetic to the Fubini--Study metric.
\end{theorem}
\begin{proof}
Since $\chi(\mathbb{CP}^2)=3$, the orbit disk is a triangle. Every pair of edges meets at a vertex, so their isotropy lines are distinct. The corresponding matrices $P_{v_1},P_{v_2},P_{v_3}$ are distinct points on $\partial\overline{\mathbb H}^2$, and hence are the vertices of a nondegenerate triangle in the affine plane of trace-one symmetric matrices. By Lemma~\ref{lem:normalized-gram-map}, the boundary map $\Psi$ traverses the boundary of this triangle once, with the original edges mapped to its vertices.

We claim that $\Psi$ is an immersion on $\Sigma^\circ$. Indeed, suppose that $\dd\Psi_p$ has rank at most one at some point $p\in\Sigma^\circ$. Then there exists a nonconstant affine function $L$ on the plane of trace-one symmetric matrices such that $\dd(L\circ\Psi)_p=0$. Choose vertices of the triangle where $L$ attains its minimum and maximum. Since $L$ is affine, $L\circ\Psi$ is nondecreasing along each of the two boundary arcs joining these vertices. Moreover, $L\circ\Psi$ is nonconstant and satisfies the scalar equation~\eqref{eq:normalized-gram-equation}, since that equation holds entrywise and has no zeroth-order term. Lemma~\ref{lem:boundary-monotonicity} therefore gives $\dd(L\circ\Psi)\neq0$ throughout $\Sigma^\circ$, a contradiction. Since $\Phi=\Psi/\sqrt{\det\Psi}$, $\Phi$ is also an immersion.

In particular, ${\rm Jac}_\Phi$ does not change its sign $\varepsilon$ on $\Sigma^\circ$.  Consequently $Q_\varepsilon=|\dd\Phi|^2+2|\operatorname{Jac}_\Phi|>0$, and Proposition~\ref{prop:weyl-sign} gives $\mathcal W^\varepsilon=0$. Since $\lambda>0$, Hitchin's classification identifies $\g$ up to homothety with the Fubini--Study metric; see \cite{Hitchin1981KahlerianTwistorSpaces} and \cite[Remark~7, p.~429]{Derdzinski1983SelfDualKahler}.
\end{proof}

\subsection{\texorpdfstring{The product metric on $S^2\times S^2$ and the Page metric on $\mathbb{CP}^2\mathbin{\sharp}\overline{\mathbb{CP}}^2$}{The product metric and the Page metric}}
\label{sec:product-page}

For these cases, we need the following auxiliary result:

\begin{lemma}[Local Weyl rigidity]\label{lem:local-wu}
Let $(N^4,h)$ be an oriented Riemannian four-manifold, fix $\varepsilon\in \{\pm1\}$, and assume that $\delta_h \mathcal W^\varepsilon=0$. Let $V\subset N$ be a relatively compact open set such that $\det {\mathcal W}^\varepsilon>0$ on $\overline V$. If ${\mathcal W}^\varepsilon$ has repeated eigenvalues on $\partial V$, then it also has repeated eigenvalues on all of $V$.
\end{lemma}

\begin{proof}
Recall that the discriminant of a trace-free symmetric $3\times3$ matrix is
\[
\mathscr D(W)=\prod_{i<j}(\lambda_i-\lambda_j)^2=\tfrac12|W|^6-27(\det W)^2.
\]
The pointwise differential inequalities in the proof of \cite[Theorem~1.2]{Wu2021} show that, for every sufficiently large integer $k$, the function
\[
F_k=|{\mathcal W}^\varepsilon|^{1/3}\left(\frac{2\mathscr D({\mathcal W}^\varepsilon)}{|{\mathcal W}^\varepsilon|^6}\right)^k
\]
is subharmonic on $V$. Indeed, the choice of $k$ depends only on positive lower bounds for $(\det {\mathcal W}^\varepsilon)/|{\mathcal W}^\varepsilon|^3$ and $(\det {\mathcal W}^\varepsilon)/|{\mathcal W}^\varepsilon|^2$, and a lower bound for the scalar curvature, all of which exist on $\overline V$. On every component of $V$ with nonempty boundary, the maximum principle gives $F_k\equiv0$. A component with empty boundary is compact, so \cite[Theorem~1.2]{Wu2021} applies directly. Hence $\mathscr D(\mathcal W^\varepsilon)\equiv0$ on $V$.
\end{proof}

\begin{theorem}\label{thm:product-page-uniqueness}
An Einstein metric admitting an effective isometric $\TT$-action on $S^2\times S^2$ or $\mathbb{CP}^2\mathbin{\sharp}\overline{\mathbb{CP}}^2$ is, up to scaling and isometry, respectively the equal-radius round product or the Page metric.
\end{theorem}
\begin{proof}
The idea is to show that the metrics are \emph{algebraically special} (namely, $\mathcal{W}^\varepsilon$ has a repeated eigenvalue), and then to apply classification results in Hermitian and K\"ahler geometry.

For both manifolds $M$ in the statement, $\chi(M)=4$. By Lemma~\ref{lem:slice-regularity}, $\Sigma$ is a geodesic quadrilateral with right angles, hence has zero total curvature by Gauss--Bonnet. Proposition~\ref{prop:weyl-sign} gives $\mu_\varepsilon<0$ on the closed disk $\Sigma$ for both signs $\varepsilon\in\{\pm1\}$: otherwise $\mathcal W^\varepsilon=0$, and Hitchin's classification, as used in the proof of Theorem~\ref{thm:fubini-study-uniqueness}, would give $\chi(M)\in\{2,3\}$. To measure the separation of the remaining eigenvalues $\nu_1,\nu_2$ of $\mathcal W^\varepsilon$, and motivated by Derdzi\'nski's conformally K\"ahler rescaling \cite[Proposition~5(i)]{Derdzinski1983SelfDualKahler}, define
\[
\mathcal G_\varepsilon:=\frac{(\nu_1-\nu_2)^2}{9\mu_\varepsilon^2}=\frac{2|\mathcal W^\varepsilon|^2}{9\mu_\varepsilon^2}-\frac13,\qquad \hat\g_\varepsilon:=(-\mu_\varepsilon)^{2/3}\hat\g.
\]
Notice that $\mathcal G_\varepsilon=1$ precisely when $\mu_\varepsilon$ is a repeated eigenvalue, while the conformal exponent $2/3$ produces a useful cancellation in the Gauss curvature formula for $\hat\g_\varepsilon$.
 Clearly, $\mathcal G_\varepsilon$ is smooth and nonnegative, with $\mathcal G_\varepsilon=1$ on the boundary by the proof of Theorem~\ref{thm:toric-hitchin-thorpe}. The conformal factor is smooth and even by Lemma~\ref{lem:distinguished-weyl}, so $\hat\g_\varepsilon$ still has geodesic edges and right angles. Gauss--Bonnet gives zero total curvature for each $\hat\g_\varepsilon$.

Write $z_\varepsilon=\log(-\mu_\varepsilon)$. Since $F_{\varepsilon}=\rho^{3/2}\mu_\varepsilon$, equation~\eqref{eqn:weyl-gap-gradient} and the definition of $\mathcal G_\varepsilon$ give
\[
\mathcal G_\varepsilon Q_\varepsilon=\left|\dd f+\tfrac23\dd z_\varepsilon\right|^2.
\]
The logarithmic form of~\eqref{eqn:weyl-full-block} is
\[
\Delta z_\varepsilon=4\lambda+3\mu_\varepsilon(1+3\mathcal G_\varepsilon)-|\dd z_\varepsilon|^2-3\langle\dd f,\dd z_\varepsilon\rangle.
\]
Combining these identities with~\eqref{eqn:weyl-Q-relation} yields
\[
\Delta z_\varepsilon=(4-3\mathcal G_\varepsilon)\lambda+3\mu_\varepsilon+\tfrac94(1-\mathcal G_\varepsilon)|\dd f|^2.
\]
The conformal curvature formula and~\eqref{eqn:distinguished-weyl} therefore give
\begin{equation}\label{eqn:four-edge-conformal-curvature}
 K_{\hat\g_\varepsilon}\dd v_{\hat\g_\varepsilon}
 =\left[(\mathcal G_\varepsilon-1)
 \left(\lambda+\tfrac34|\dd f|^2\right)
 -\tfrac\varepsilon2\operatorname{Jac}_\Phi\right]\dd v_{\hat\g}.
\end{equation}
Lemma~\ref{lem:distinguished-weyl} gives a smooth extension of $\operatorname{Jac}_\Phi$. The remaining terms are integrable by Lemma~\ref{lem:slice-regularity}: evenness and $\mathcal G_\varepsilon=1$ on every edge give $\mathcal G_\varepsilon-1=O(r^2)$ there and $O(r_1^2r_2^2)$ at a corner. Adding the two zero total-curvature identities cancels the Jacobians and yields
\begin{equation}\label{eqn:four-edge-total-curvature}
 0=\int_\Sigma(\mathcal G_++\mathcal G_--2)
       \left(\lambda+\frac34|\dd f|^2\right)\dd v_{\hat\g}.
\end{equation}

It remains to show that $\mathcal G_\varepsilon\leq1$, or equivalently, that $\mu_\varepsilon$ is the smallest eigenvalue of $\mathcal W^\varepsilon$ everywhere on $\Sigma$. Otherwise, let $U$ be a connected component of $\{\mathcal G_\varepsilon>1\}$ and set $V=\pi^{-1}(U)$. On $\overline U$, the eigenvalues satisfy $\nu_1\leq\mu_\varepsilon<0<\nu_2$, so $\det\mathcal W^\varepsilon>0$ on $\overline V$. Moreover, $\mathcal G_\varepsilon=1$ along $\partial U$, hence $\mathcal W^\varepsilon$ has repeated eigenvalues on $\partial V$. Lemma~\ref{lem:local-wu} then gives repeated eigenvalues throughout $V$, contradicting $\mathcal G_\varepsilon>1$ on $U$. Thus $\mathcal G_\varepsilon\leq1$.

Equation~\eqref{eqn:four-edge-total-curvature} now forces $\mathcal G_+=\mathcal G_-=1$, from which  $\mathcal W^\varepsilon$ has eigenvalues $\mu_\varepsilon,\mu_\varepsilon,-2\mu_\varepsilon$ everywhere. Derdzi\'nski's result \cite[Proposition~5(i)]{Derdzinski1983SelfDualKahler} then makes $(-\mu_\varepsilon)^{2/3}\g$ locally K\"ahler. The local K\"ahler forms have norm $\sqrt2$ with respect to $(-\mu_\varepsilon)^{2/3}\g$ and span the eigenline of $-2\mu_\varepsilon$. They therefore differ by signs on overlaps. Since $H^1(M;\mathbb{Z}_2)=0$, they glue to a global K\"ahler form in each orientation. In particular, $\g$ is Hermitian. In the non-K\"ahler case, LeBrun's classification \cite[Theorem~A]{LeBrun2012EinsteinHermitian} gives the Page metric: its other non-K\"ahler possibility, the Chen--LeBrun--Weber metric, has Euler characteristic five.

In the K\"ahler--Einstein case, $\lambda>0$ makes the complex surface del Pezzo with $c_1^2=2\chi+3\etau=8$. The classification of del Pezzo surfaces gives either $\mathbb{CP}^1\times\mathbb{CP}^1$ or  the one-point blow-up of $\mathbb{CP}^2$; see \cite[proof of Theorem~B]{LeBrun2015HarmonicForms}. The one-point blow-up admits no constant-scalar-curvature K\"ahler metric \cite[Remark~7]{Derdzinski1983SelfDualKahler}. On $\mathbb{CP}^1\times\mathbb{CP}^1$, Bando--Mabuchi uniqueness \cite[Theorem~A(ii)]{BandoMabuchi1987} identifies $\g$, up to scale and isometry, with the equal-radius round product.
\end{proof}

\section{Existence strategy}\label{sec:strategy}


In this section, we introduce our framework for proving Theorem \ref{thm:intro-existence}. As mentioned above, the idea is to construct $\TT$-invariant Einstein metrics on $S^4\setminus\mathcal T$, where $\mathcal T$ is a given orbit, with conical singularities along $\mathcal T$ of cone angle $2\pi/d$, $d \in \mathbb N$. That is, the metrics satisfy $\operatorname{Ric}(\g)=\lambda\g$ away from $\mathcal T$, and near $\mathcal T$ they are asymptotic to
\[
a_0(\dd r^2+\beta^2r^2\dd\vartheta^2)+b_0(\dd\vartheta_1^2+\dd\vartheta_2^2),
\qquad a_0,b_0>0,
\]
where $(r,\vartheta)$ are polar coordinates normal to $\mathcal T$. The precise asymptotic regularity will be specified in detail in Section~5.

In $\S$\ref{sec:branched-cover} we explain how one may view $\sharp_{d-1}(S^2\times S^2)$ as a $\TT$-equivariant ramified $d$-fold cover of $S^4$, ramified along $\mathcal T$. The  pullback of the conical Einstein metrics will be later shown to  extend smoothly across the ramification torus. 

To prove existence of the conical Einstein metrics, we vary the cone angle continuously downwards from $\beta = 1/2$ (corresponding to $S^2 \times S^2$), and use a mod-two degree argument for the toric conical Einstein equation, the details of which are contained in $\S$\ref{sec:localized-degree}.

\subsection{Branched covers and conformal coordinates}\label{sec:branched-cover}

\phantomsection\label{notloc:s4-action}
Fix the Einstein constant $\lambda>0$. Let $\TT$ act on $S^4\subset\mathbb C^2\times\mathbb{R}$ by
$(e^{i\vartheta_1},e^{i\vartheta_2})\cdot(z_1,z_2,t)=(e^{i\vartheta_1}z_1,e^{i\vartheta_2}z_2,t)$, and let $\pi:S^4\to\Sigma$ be the quotient map. The orbit disk has vertices $p_\pm:=\pi(0,0,\pm1)$. The Clifford torus $\mathcal T=\{t=0,\ |z_1|=|z_2|=1/\sqrt2\}$ is a principal orbit; write $c:=\pi(\mathcal T)$. The commuting involutions $\sigma(z_1,z_2,t)=(z_2,z_1,t)$ and $\tau(z_1,z_2,t)=(z_1,z_2,-t)$ preserve $\mathcal T$. Set 
\[
 \Gamma:=\langle\sigma,\tau\rangle\simeq\mathbb{Z}_2\times\mathbb{Z}_2, \qquad \mathsf G:=\TT\rtimes\Gamma,
\]
with $\sigma$ exchanging the circle factors and $\tau$ acting trivially on $\TT$.


\phantomsection\label{notloc:branched-cover}
For each integer $m\geq1$, put $d=m+1$. Consider a regular $2d$-gon $P$ whose successive edges are labelled alternately by the two coordinate circle subgroups of $\TT$. Since adjacent labels generate $\TT$, the Orlik--Raymond reconstruction theorem~\cite[\S4.2]{OrlikRaymond1} gives a smooth closed $4$-manifold $\widetilde X_m$ with an effective $\TT$-action and orbit space $P$. 
The resulting manifold is diffeomorphic to $\sharp_{d-1}(S^2\times S^2)$. Indeed, the characteristic spheres all have self-intersection zero, since their two neighbouring isotropy subgroups coincide. Thus the intersection form is even, excluding projective-plane summands in the classification. The number of summands follows from $\chi(\widetilde X_m)=2d$.

The rotation $R$ of $P$ by a $2\pi/d$-angle preserves the characteristic labels and therefore lifts to an order-$d$ diffeomorphism of $\widetilde X_m$ commuting with $\TT$. Indeed, in the standard reconstruction $\widetilde X_m=(P\times\TT)/{\sim}$, the lift is induced by $(q,t)\mapsto(Rq,t)$. The resulting $\mathbb Z_d$-action fixes precisely the principal orbit over the centre of $P$ and is free elsewhere. Its quotient orbit space $P/\mathbb Z_d$ is a bigon with characteristic labels $e_1,e_2$, and hence the quotient manifold is equivariantly diffeomorphic to $S^4$ with the standard $\TT$-action, again by \cite{OrlikRaymond1}. We choose this identification to send the central orbit to the Clifford torus $\mathcal T$. Near the fixed torus, the cyclic action is rotation in the normal plane, so the quotient map $p_m:\widetilde X_m\to S^4$ has the local form
\[
p_m(y,w)=(y,w^d),\qquad y\in\mathcal T.
\]
Thus $p_m$ is a cyclic $d$-fold branched covering, with ramification locus the central torus and branch locus $\mathcal T$.

\subsubsection{Doubling the orbit space}

In later sections it will be useful to work with a fixed surface obtained by reflecting the orbit disk across its boundary. We use both the four-copy global section $S_\Sigma^2$, which provides the setting for the scalar Fredholm analysis in Section~\ref{sec:fredholm}, and the two-copy conformal double $\widehat\Sigma$, which provides the normalized conformal coordinates.

\phantomsection\label{notloc:double}
The sphere $S_{\Sigma}^2 =S^4\cap\{z_1,z_2\in\mathbb{R}\}$ has coordinates $(r_1,r_2,t)$. Identify $\Sigma$ with its quadrant $r_1,r_2\geq0$, equipped with the round conformal structure. Thus $S_\Sigma^2$ consists of four reflected copies of $\Sigma$. Let $\widehat{\Sigma}$ be the conformal double across $\partial\Sigma$, let $j$ interchange its two copies, and denote the lifts of $c$ by $c_\pm$. Reflection defines the complex structure along each open edge; at a vertex, the square of a holomorphic quadrant coordinate defines a coordinate on the double. Both $j$ and $\sigma$ are antiholomorphic.

\begin{lemma}\label{lem:fixed-marking}
For every conformal structure on $\widehat{\Sigma}$ making $j$ and $\sigma$ antiholomorphic, there is a unique biholomorphism $Z:\widehat{\Sigma}\to\mathbb{CP}^1$ with
\begin{equation}\label{eq:fixed-marking}
 Z(p_-)=0,\qquad Z(p_+)=\infty,\qquad Z(c_+)=i,\qquad Z(c_-)=-i.
\end{equation}
\end{lemma}
\begin{proof}
Uniformization \cite{Forster81} gives a unique $Z$ with the first three prescribed values. Since $j$ and $\sigma$ are commuting antiholomorphic involutions, $j\sigma$ is a holomorphic involution fixing $p_\pm$ and exchanging $c_\pm$. In particular, it is nontrivial. In the $Z$-coordinate it must therefore be $Z\mapsto-Z$, which gives the fourth value.
\end{proof}

\begin{remark}
Alternatively, one could work directly with the four-copy sphere $S_\Sigma^2$, which identifies with the global section for the polar action in $M$. Indeed, the involution $(r_1,r_2,t)\mapsto(-r_1,-r_2,t)$ identifies opposite quadrants, and its quotient is precisely the conformal double $\widehat\Sigma$. The quotient map is a two-fold branched covering, ramified at $p_\pm$. Thus the conformal normalization above could equivalently be carried out on $S_\Sigma^2$.
\end{remark}

\subsection{Localized mod-two degree}\label{sec:localized-degree}

Our existence argument uses a mod-two degree for the Einstein equation. Roughly speaking, this counts the number of solutions modulo two after a small generic perturbation. The count is unchanged under continuous deformations, provided solutions cannot escape to infinity. We recall the version of this principle needed here, based on ~\cite[\S3]{Smale1965Sard} and \cite[\S\S2,4]{Milnor1965}, including a proof for the reader's convenience.

For the application in this paper, $s$ will be the Einstein operator, $x$ the metric, $t$ the cone angle, and $P_{t,x}$ the projection onto  the contracted-Bianchi subspace.

\begin{proposition}\label{prop:localized-degree}
Let $E,F$ be real Banach spaces, $U\subset E$ open, and $P_x:F\to F$, $x\in U$, a smooth family of bounded projections. Suppose $s:U\to F$ is smooth and satisfies $s(x)\in P_xF$, for all $x\in U$. If the zero set $Z=s^{-1}(0)$ is compact and
\[
\mathrm D_xs:E\longrightarrow P_xF
\]
is Fredholm of index $0$ at every $x\in Z$, then $s$ has a well-defined localized degree $\deg_2(s)\in\mathbb Z_2$, with the following properties:
\begin{enumerate}
    \item \emph{(Count)} If all zeros are nondegenerate, then $\deg_2(s)=\#Z\pmod 2$. In particular, if $\deg_2(s) \neq 0$ then $Z \neq \varnothing$.
    \item \emph{(Homotopy invariance)} For a smooth family $(P_{t,x} : F \to F,s_t)$, $t\in[0,1]$, satisfying the same Fredholm condition at every zero,  $\deg_2(s_t)$ is independent of $t$ provided the total zero set $\{(t,x)\in[0,1]\times U:s_t(x)=0\}$ is compact.
\end{enumerate}
\end{proposition}

\begin{proof}
The idea is to reduce the equation to a finite-dimensional problem near its zeros. Since the linearizations are Fredholm, their cokernels are finite-dimensional. Compactness of $Z$ therefore allows us to choose a finite-dimensional subspace $V\subset F$ such that
\[
\operatorname{im}\mathrm D_xs+P_xV=P_xF
\qquad\text{for all }x\in Z.
\]
Consequently, the implicit function theorem shows that the solutions of the perturbed equation
\[
s(x)=P_xv,\qquad v\in V,
\]
form a smooth manifold $M$ of dimension $\dim V$ near $Z\times\{0\}$. Let $p:M\to V$ denote the projection onto the perturbation parameter.

Choose a compact neighborhood $N\subset M$ of $Z\times\{0\}$, with smooth boundary disjoint from $p^{-1}(0)$. For every sufficiently small regular value $v\in V$, the fiber $N\cap p^{-1}(v)$ is finite. We define
\[
\deg_2(s):=\#(N\cap p^{-1}(v))\pmod 2.
\]
The parity is independent of $v$: a generic path between two regular values produces a compact one-dimensional cobordism between their fibers, and the boundary of a compact one-manifold has even cardinality. Shrinking $N$ does not change the count, while independence of $V$ follows by comparing two choices with their sum. Thus the degree depends only on the original section near its zero set.

For a smooth family, compactness of the total zero set permits a common finite-dimensional reduction. A generic small perturbation then gives a compact one-manifold whose boundary consists of the perturbed solutions at $t=0$ and $t=1$. Their numbers agree modulo two, proving homotopy invariance.
\end{proof}

\section{Conical metrics and weighted spaces}\label{sec:conical}

The purpose of this section is to turn the geometric problem of finding conical Einstein metrics into a well-defined nonlinear elliptic problem on fixed Banach spaces, and to establish the regularity needed to recover smooth Einstein metrics on branched covers.

More precisely, we develop a framework for the conical Einstein equation using weighted H\"older spaces on the real blow-up of the orbit disk, based on Melrose's $b$-calculus \cite[Chapters~2, 4]{Melrose1993}. We  prove that, within a diagonal toric Ansatz (Definition \ref{def:X_beta}), the  Einstein equation is equivalent to a reduced system of three scalar equations (Proposition \ref{prop:reduced-equations}), using the Bianchi identity to recover the missing horizontal equations. Finally, we establish regularity of solutions (Lemma~\ref{lem:cone-regularity}), including smooth extension to branched covers when the cone angle is  $2\pi/d$, $d\in \mathbb N$ (Proposition~\ref{prop:removability}).


Throughout, we fix $0 < \beta_0<\frac{1}{2}$, an integer $k \geq 2$, and $0 < \alpha<1$. Set $I := \left[ \beta_0, \frac{1}{2}  \right]$, and choose $0 < \delta < \min \{ 1, 2 \beta_0 \}$.

\subsection{The fixed blow-up and model operators}
\phantomsection\label{notloc:conformal-coordinate}
Choose a conformal coordinate $z$ centered at $c$ with $\sigma(z)=-\overline z$ and $\tau(z)=\overline z$. For the round structure, stereographic projection from the antipode of $c$, with suitable axes, gives such a coordinate.

\begin{definition}
\phantomsection\label{notloc:blow-up}
The \emph{real blow-up} $\Sigma_b:=[\Sigma;c]$ replaces $c$ by its circle of directions $F_c:=(T_c\Sigma\setminus\{0\})/\mathbb{R}_{>0}$. It is a compact smooth surface with corners (topologically an annulus), retaining the original corners of $\Sigma$, and acquiring an additional smooth boundary component $F_c$. Its smooth structure agrees with that of $\Sigma$ off $F_c$ and is defined near $F_c$ by $z=xe^{i\vartheta}$, with $x\geq0$ and $\vartheta\in\mathbb{R}/2\pi\mathbb{Z}$. The blow-down map collapses the boundary circle $F_c$ to $c$ and is the identity elsewhere.
\end{definition} 

Extend $x$ to a $\Gamma$-invariant function in $\mathcal{C}^{\infty}(\Sigma_b)$, positive on $\Sigma_b \setminus F_c$ and equal to $1$ near $\partial \Sigma$.

\phantomsection\label{notloc:collar}
Choose $0<\varepsilon<1$ so that $\{x\leq\varepsilon\}$ lies in the coordinate neighborhood of $F_c$. Set $\Omega=[0,\varepsilon]\times S^1$ and $\Omega^\circ=\Omega\setminus F_c$. Smooth functions on $\Omega$ are smooth up to both boundary circles.

\begin{definition}
\phantomsection\label{notloc:b-geometry}
A \emph{$b$-vector field} is a smooth vector field tangent to $F_c$, equivalently one satisfying $Vx\in x\mathcal C^\infty(\Omega)$. These are the sections of the rank-two bundle ${}^b\mathcal T_\Omega$, with frame $x\partial_x,\partial_\vartheta$ and dual co-frame $\dd x/x,\dd\vartheta$. Set
\begin{equation*}
 \g_b=\left(\frac{\dd x}{x}\right)^2+\dd\vartheta^2,
 \qquad \g_\beta=x^{2\beta}\g_b,\qquad 0<\beta<1.
\end{equation*}
The \emph{cylinder metric} $\g_b$ is smooth  on ${}^b\mathcal T_\Omega$; $\g_\beta$ is a flat cone metric of angle $2\pi\beta$ on $\Omega^\circ$.
\end{definition}

\begin{remark}\label{rem:radial-convention}
For $r_\beta = x^\beta/\beta$, we have $\g_\beta=\dd r_\beta^2+\beta^2r_\beta^2\dd\vartheta^2$ and $x\partial_x=\beta \, r_\beta\partial_{r_\beta}$. The cylinder metric $\g_b$ is complete toward $F_c$, whereas the cone metric places $F_c$ at finite distance. For nonintegral $\beta$, $r_\beta$ is not a smooth boundary defining function in the fixed smooth structure. The smooth structure and weighted spaces are defined using $x$; we also use $r_\beta$ to express the cone metric and its coframe.
\end{remark}

We use the notation $H=O_b(x^a)$ to mean that $H$ and all its derivatives with respect to the $b$-vector fields $x\partial_x$ and $\partial_\vartheta$ are $O(x^a)$.

\begin{definition}
\phantomsection\label{notloc:b-operators}
A \emph{$b$-differential operator} of order $k$ has the form
\begin{equation*}
 L=\sum_{a+b\leq k}a_{ab}(x,\vartheta)(xD_x)^aD_\vartheta^b,
 \qquad a_{ab}\in\mathcal C^\infty(\Omega;\mathbb C),
 \qquad D_x=-i\partial_x,\quad D_\vartheta=-i\partial_\vartheta.
\end{equation*}
Its \emph{$b$-principal symbol} at $\xi\,\dd x/x+\eta\,\dd\vartheta$ is $\sum_{a+b=k}a_{ab}\xi^a\eta^b$. It is \emph{$b$-elliptic} if this polynomial is nonzero on every nonzero $b$-covector. 

Replacing each coefficient $a_{ab}(x,\vartheta)$ by its boundary value $a_{ab}(0,\vartheta)$ gives the `frozen-coefficients' operator $L_0$, which is translation-invariant in the cylindrical coordinate $t=-\log x$. Its \emph{indicial family}, with powers measured in $x$, is 
\begin{equation*}
 I(L,s)\varphi:=x^{-s}L_0(x^s\varphi)
 =\sum_{a+b\leq k}a_{ab}(0,\vartheta)(-is)^aD_\vartheta^b\varphi,
 \qquad \varphi=\varphi(\vartheta).
\end{equation*}
For $b$-elliptic $L$, an \emph{indicial root} is a complex number $s$ for which $I(L,s):H^k(S^1)\to L^2(S^1)$ is not invertible. A real weight is \emph{admissible} if no root has that real part.
\end{definition}

A $b$-elliptic operator is not uniformly elliptic up to $x=0$ in the usual smooth structure. Its indicial roots describe the possible leading powers $x^s$ of solutions near this boundary, and hence determine the weights at which elliptic estimates may fail.

\begin{example}
\phantomsection\label{notloc:b-laplacian}
The \emph{$b$-Laplacian} is the Laplace-Beltrami operator of $g_b$, $$\Delta_b \, = \, (x\partial_x)^2+\partial_\vartheta^2.$$ 
It has principal symbol $-(\xi^2+\eta^2)$ and indicial family $s^2+\partial_\vartheta^2$. On $e^{in\vartheta}$, the latter is multiplication by $s^2-n^2$. The indicial exponents measured in $x$ are therefore $\mathbb{Z}$, and those measured in $r=x^\beta/\beta$ are $\beta^{-1}\mathbb{Z}$. The solutions associated with the exponent $0$ are $1$ and $\log x$.
\end{example}

\begin{example}\label{ex:quotient-b-operator}
Consider the product conical metric $\g=\g_\beta+\g_{\mathcal T}$ on $\Omega^\circ\times\mathcal T$, where $\g_{\mathcal T}$ is translation-invariant. Its rescaled Laplacian is
\[
-x^{2\beta}\Delta_\g=-\Delta_b-x^{2\beta}\Delta_{\g_{\mathcal T}}.
\]
On $\TT$-invariant functions, the second term vanishes, and we thus obtain the $b$-elliptic operator $-\Delta_b$ on the orbit space. On general functions, however, the rescaled operator is not $b$-elliptic: its limiting principal symbol vanishes on nonzero covectors tangent to $\mathcal T$. Moreover, the coefficient $x^{2\beta}$ is not smooth at $x=0$ unless $2\beta$ is an integer.
\end{example}

\subsection{Weighted H\"older spaces}

\phantomsection\label{notloc:b-distance}

\subsubsection{Definitions and basic properties}

All norms use the fixed cylinder metric $\g_b$ and are independent of $\beta$. Write $d_b$ for the distance function associated to $\g_b$.

\begin{definition}\label{def:little-spaces}
Fix $0<\alpha<1$. For $m\in\mathbb{N}_0$ and $V^\nu=(x\partial_x)^{\nu_1}\partial_\vartheta^{\nu_2}$, set
\begin{align*}
 [a]_{\alpha;b}
 &=\sup_{p\ne p'}\frac{|a(p)-a(p')|}{d_b(p,p')^\alpha},\\
 \|a\|_{\mathcal C_b^{m,\alpha}}
 &=\sum_{|\nu|\leq m}\|V^\nu a\|_\infty
   +\sum_{|\nu|=m}[V^\nu a]_{\alpha;b},
 \qquad p,p'\in\Omega^\circ.
\end{align*}
Define
\begin{eqnarray*}
 \mathcal C_b^{m,\alpha}(\Omega)
 &=& \{a\in\mathcal{C}^m(\Omega\setminus F_c):\|a\|_{\mathcal C_b^{m,\alpha}}<\infty\},
\end{eqnarray*}
and let $c_b^{m,\alpha}(\Omega)$ be the closure of $\mathcal{C}^{\infty}(\Omega)$ in this norm. The big space imposes no ordinary $\mathcal C^m$ regularity at $F_c$; the smooth functions defining the little space extend through $x=0$. For $t\in\mathbb{R}$, multiplication by $x^t$ defines the weighted spaces $x^t\mathcal C_b^{m,\alpha}$ and $x^tc_b^{m,\alpha}$, with norm $\|u\|_{x^t\mathcal C_b^{m,\alpha}}=\|x^{-t}u\|_{\mathcal C_b^{m,\alpha}}$ and its restriction to the little space.
\end{definition}

\begin{lemma}
    The $b$-H\"older space $\mathcal C_b^{m,\alpha}(\Omega)$ and the little $b$-H\"older space $c_b^{m,\alpha}(\Omega)$ are Banach spaces with respect to the norm $\| \cdot \|_{\mathcal{C}_b^{m,\alpha}}$.
\end{lemma} \begin{proof}
    Set $Z=[0,\infty)\times S^1$, with $\tau=\log(\varepsilon/x)$. Then $\g_b=\dd\tau^2+\dd\vartheta^2$ and $x\partial_x=-\partial_\tau$. Let $\mathcal C^{m,\alpha}(Z)$ denote the bounded ordinary H\"older space for this product metric. The map
    \begin{eqnarray*}
        \mathcal C_b^{m,\alpha}(\Omega)\ni a(x,\vartheta)
        &\longmapsto& a(\varepsilon e^{-\tau},\vartheta)\in\mathcal C^{m,\alpha}(Z)
    \end{eqnarray*}
    is an isometry onto the Banach space $\mathcal C^{m,\alpha}(Z)$, so $\mathcal C_b^{m,\alpha}(\Omega)$ is Banach.
     Finally, the space $c_b^{m,\alpha}$ is  a closed subspace of the Banach space $\mathcal{C}_b^{m,\alpha}(\Omega)$.
\end{proof}

   

\begin{lemma}\label{lem:little-cutoffs}
Restriction to $F_c$ extends uniquely to a bounded linear operator $$\text{Tr} : c_b^{m,\alpha}(\Omega) \longrightarrow c^{m,\alpha}(S^1)$$ whose kernel is the closure in $c_b^{m,\alpha}(\Omega)$ of the smooth functions vanishing near $F_c$.
\end{lemma} \begin{proof}
For smooth $a$, continuity at $x=0$ implies $\|\partial_{\vartheta}^j a(0,\cdot)\|_{\infty} \leq \|\partial_{\vartheta}^j a\|_{\infty}$, for $0\leq j\leq m$. Moreover, $[\partial_{\vartheta}^m a(0,\cdot)]_{\alpha} \leq [\partial_{\vartheta}^m a]_{\alpha;b}$. Hence, the restriction map is bounded and extends uniquely to $\text{Tr}:c_b^{m,\alpha}(\Omega)\to c^{m,\alpha}(S^1)$. Suppose $\text{Tr}(a)=0$, and choose smooth functions $a_j$ such that $\|a_j-a\|_{c_b^{m,\alpha}}\to0$. Set $b_j(x,\vartheta):=a_j(x,\vartheta)-a_j(0,\vartheta)$. Extending functions on $S^1$ independently of $x$ is bounded in these norms. Hence,
\begin{eqnarray*}
    \|b_j-a\|_{c_b^{m,\alpha}} &\leq& \|a_j-a\|_{c_b^{m,\alpha}} +C\|\text{Tr}(a_j)\|_{c^{m,\alpha}(S^1)} \,\longrightarrow\,0.
\end{eqnarray*}
Since $b_j(0,\vartheta)=0$, we can write $b_j=xh_j$, where $h_j$ is smooth. Choose $\chi\in\mathcal C^\infty([0,\infty))$ equal to zero on $[0,1]$ and one on $[2,\infty)$, and set $\chi_\eta(x):=\chi(x/\eta)$. On the support of $1-\chi_\eta$, we have $x\leq2\eta$, and the $b$ derivatives of $\chi_\eta$ are bounded independently of $\eta$. The Leibniz rule therefore bounds every $b$ derivative of $(1-\chi_\eta)xh_j$ through order $m+1$ by $C_j\eta$. The additional derivative controls the H\"older seminorm, giving
\begin{eqnarray*}
    \|(1-\chi_\eta)b_j\|_{c_b^{m,\alpha}}
    &\leq& C_j'\eta.
\end{eqnarray*}
Choose $\eta_j\downarrow0$ such that $C_j'\eta_j\leq1/j$. Then $\chi_{\eta_j}b_j$ is smooth, vanishes near $F_c$, and
\begin{eqnarray*}
    \|\chi_{\eta_j}b_j-a\|_{c_b^{m,\alpha}}
    &\leq& \frac{1}{j}+\|b_j-a\|_{c_b^{m,\alpha}}
    \,\longrightarrow\,0.
\end{eqnarray*}
Conversely, every smooth function vanishing near $F_c$ has zero trace, and continuity of $\text{Tr}$ gives the same conclusion for its closure.
\end{proof}

\subsubsection{Decay and smooth dependence on the cone angle}
 
Write $c_{\mathrm{loc}}^{m,\alpha}(U)$ for the space of functions whose restriction to each compact smooth subdomain $K\subset U$ belongs to $c^{m,\alpha}(K)$.

\begin{lemma}\label{lem:decay-implies-little}
Let $\Omega'=[0,\varepsilon']\times S^1$, where $0<\varepsilon'<\varepsilon$. If $u\in c_{\mathrm{loc}}^{m,\alpha}(\Omega\setminus F_c) \cap x^{t+q}\mathcal C_b^{m,\alpha}(\Omega)$ for some $t\in\mathbb{R}$ and $q>0$, then $u\in x^tc_b^{m,\alpha}(\Omega')$.
\end{lemma}

\begin{proof}
Choose $\chi\in\mathcal{C}^{\infty}([0,\infty))$ equal to zero on $[0,1]$ and to one on $[2,\infty)$, and set $\chi_\eta(x)=\chi(x/\eta)$, where $0<\eta<\varepsilon'/2$. Set $\rho(s)=(1-\chi(s))s^q$. Since $(s\partial_s)^j\rho(s)=q^js^q$ for $s\leq1$ and $\rho(s)=0$ for $s\geq2$, we have
$\rho\in\mathcal C_b^{m,\alpha}((0,\infty)\times S^1)$. Writing $s=x/\eta$, we obtain
\begin{equation*}
 \|(1-\chi_\eta)x^q\|_{\mathcal C_b^{m,\alpha}(\Omega')} \, \leq \, \eta^q \|\rho\|_{\mathcal C_b^{m,\alpha}((0,\infty)\times S^1)} \, \leq \,  C\eta^q.
\end{equation*}
Here, $x\partial_x=s\partial_s$ and $\mathrm dx/x=\mathrm ds/s$. Leibniz's rule and $[fg]_{\alpha;b}\leq\|f\|_\infty[g]_{\alpha;b} +\|g\|_\infty[f]_{\alpha;b}$ give
\begin{align*}
 \|(1-\chi_\eta)x^{-t}u\|_{\mathcal C_b^{m,\alpha}(\Omega')}
 &\leq C\|(1-\chi_\eta)x^q\|_{\mathcal C_b^{m,\alpha}(\Omega')}
          \|x^{-t-q}u\|_{\mathcal C_b^{m,\alpha}(\Omega')}\\
 &\leq C\eta^q
          \|u\|_{x^{t+q}\mathcal C_b^{m,\alpha}(\Omega)}.
\end{align*}

For fixed $\eta$, set $K_\eta=[\eta/2,\varepsilon']\times S^1$. Since $x^{-t}u\in c^{m,\alpha}(K_\eta)$, choose $v_j\in\mathcal{C}^{\infty}(K_\eta)$ with
\begin{equation*}
 \|v_j-x^{-t}u\|_{\mathcal C^{m,\alpha}(K_\eta)} \longrightarrow0.
\end{equation*}
Define $a_j=\chi_\eta v_j$ on $K_\eta\cap\Omega'$ and extend $a_j$ by zero for $x<\eta/2$. Then $a_j\in\mathcal{C}^{\infty}(\Omega')$, and
\begin{equation*}
 \|a_j-\chi_\eta x^{-t}u\|_{\mathcal C_b^{m,\alpha}(\Omega')} \leq C_\eta \|v_j-x^{-t}u\|_{\mathcal C^{m,\alpha}(K_\eta)} \longrightarrow0.
\end{equation*}
Hence, $\chi_\eta x^{-t}u\in c_b^{m,\alpha}(\Omega')$. Since $C\eta^q\|u\|_{x^{t+q}\mathcal C_b^{m,\alpha}(\Omega)} \to0$, we have $\chi_\eta x^{-t}u\to x^{-t}u$ in $\mathcal C_b^{m,\alpha}(\Omega')$. Closedness of $c_b^{m,\alpha}(\Omega')$ gives $x^{-t}u\in c_b^{m,\alpha}(\Omega')$, so $u\in x^tc_b^{m,\alpha}(\Omega')$.
\end{proof}

\begin{corollary}[Smooth dependence on the cone angle]\label{cor:4.13}
    Fix $\delta>0$, $m \in \mathbb{N}_0$, and $0<\alpha<1$. The map $\left( \frac{\delta}{2}, \infty \right) \to x^{\delta} c_b^{m,\alpha}(\Omega)$, $\beta \longmapsto x^{2\beta}$, is smooth, with \begin{eqnarray*}
        \frac{\dd^k}{\dd \beta^k} x^{2\beta} &=&2^k x^{2\beta} (\log x)^k.
    \end{eqnarray*}
\end{corollary} \begin{proof}
    On a compact interval $J \Subset \left( \frac{\delta}{2}, \infty \right)$, choose $0 < q < \min_{\beta \in J}(2\beta-\delta)$. Set $\widetilde\Omega=[0,\widetilde\varepsilon]\times S^1$, with $\widetilde\varepsilon>\varepsilon$. For every $j\in\mathbb{N}_0$, the functions $x^{2\beta-\delta}(\log x)^j$ are smooth for $x>0$ and have uniformly bounded $x^q\mathcal C_b^{m,\alpha}(\widetilde\Omega)$ norms. Applying Lemma~\ref{lem:decay-implies-little} with $t=0$ and $\Omega$ as the smaller collar gives $x^{2\beta-\delta}(\log x)^j\in c_b^{m,\alpha}(\Omega)$. Set $p_j(\beta) = 2^j x^{2\beta}(\log x)^j$. The $b$-derivatives of $x^{-\delta}p_j(\beta)$ through order $m+1$ are uniformly bounded for $\beta \in J$. Hence,
\begin{eqnarray*}
    \sup_{\beta \in J} \| p_j(\beta) \|_{x^\delta c_b^{m,\alpha}(\Omega)} &\leq& C_j.
\end{eqnarray*}
Taylor's formula, applied pointwise, gives
\begin{eqnarray*}
    p_j(\beta+h) - p_j(\beta) - hp_{j+1}(\beta)
    &=& h^2 \int_0^1 (1-s)p_{j+2}(\beta+sh) \, \dd s
\end{eqnarray*}
whenever $\beta,\beta+h \in J$. Taking the weighted $b$-Hölder norm yields
\begin{eqnarray*}
    \left\|
        \frac{p_j(\beta+h)-p_j(\beta)}{h} - p_{j+1}(\beta)
    \right\|_{x^\delta c_b^{m,\alpha}(\Omega)}
    &\leq& \frac{C_{j+2}}{2}|h|.
\end{eqnarray*}
Thus $p_j' = p_{j+1}$ in $x^\delta c_b^{m,\alpha}(\Omega)$ on the interior of $J$, for every $j \in \mathbb{N}_0$. Since $J$ was arbitrary, $\beta \mapsto x^{2\beta}$ is smooth, with the stated derivatives.
\end{proof}

\subsubsection{Elliptic estimates and compactness}

The following weighted elliptic estimate improves decay near the conical orbit up to the first positive indicial root $1$, and isolates the first nontrivial Fourier modes. It will be used both to establish curvature regularity and to obtain uniform estimates in the compactness argument.

\begin{proposition}\label{prop:fourier-estimate}
Let $\Omega'=[0,\varepsilon']\times S^1\subset\Omega$, with $0<\varepsilon'<\varepsilon$. Fix $k\in\mathbb{N}_0$, $0<\alpha<1$, $s,t>0$, and $0<a\leq1$ with $a<t$. If $H\in x^sc_b^{k+2,\alpha}(\Omega)$ and $\Delta_bH\in x^tc_b^{k,\alpha}(\Omega)$, then
\begin{equation*}
 \|H\|_{x^ac_b^{k+2,\alpha}(\Omega')}
 \leq C\bigl(\|H\|_{\mathcal C^0(\Omega)}
       +\|\Delta_bH\|_{x^tc_b^{k,\alpha}(\Omega)}\bigr).
\end{equation*}
If $t>1$, there is a unique $\ell\in\operatorname{span}\{\cos\vartheta,\sin\vartheta\}$ such that, for every $1<b<\min\{2,t\}$,
\begin{equation*}
 \|\ell\|_{\mathcal C^{k+2,\alpha}(S^1)}
 +\|H-x\ell\|_{x^bc_b^{k+2,\alpha}(\Omega')}
 \leq C\bigl(\|H\|_{\mathcal C^0(\Omega)}
       +\|\Delta_bH\|_{x^tc_b^{k,\alpha}(\Omega)}\bigr).
\end{equation*}
The constants depend only on $k,\alpha,t,\varepsilon,\varepsilon'$ and, respectively, $a$ or $b$.
\end{proposition}

\begin{proof}
Write $z=xe^{i\vartheta}$, $F=\Delta_bH$, and
\[
M=\|H\|_{\mathcal C^0(\Omega)}+\|F\|_{x^tc_b^{k,\alpha}(\Omega)}.
\]
Since $\Delta_b=x^2\Delta_{\mathbb R^2}$, we have $\Delta_{\mathbb R^2}H=x^{-2}F$ on the punctured disk. Fix $\varepsilon'<R<\varepsilon$, and let $u$ be the Dirichlet Green potential on $D_R=\{|z|<R\}$ of the function $x^{-2}F$. Since $|x^{-2}F|\leq Mx^{t-2}$ and $t>0$, the potential extends continuously across the origin, with $\|u\|_{\mathcal C^0(D_R)}\leq CM$. The function $h=H-u$ is bounded and harmonic on the punctured disk, hence extends harmonically across the origin, with $\|h\|_{\mathcal C^0(D_R)}\leq M$. Moreover, $H(0)=0$, so $h(0)=-u(0)$.

Elementary estimates for the logarithmic Green kernel give, for every $0<q<\min\{1,t\}$,
\[
|u(z)-u(0)|\leq C_qM|z|^q.
\]
If $t>1$, the potential is differentiable at the origin, and for every $1<q<\min\{2,t\}$,
\[
|\nabla u(0)|\leq CM,\qquad
|u(z)-u(0)-\langle\nabla u(0),z\rangle|
\leq C_qM|z|^q.
\]
Indeed, splitting the Green integral into $|y|\leq2|z|$ and $|y|>2|z|$, the near contribution is $O(|z|^t)$, while the far contribution is controlled by
\[
C|z|\int_{2|z|}^R r^{t-2}\,\dd r,
\]
or, after subtracting the linear term, by
\[
C|z|^2\int_{2|z|}^R r^{t-3}\,\dd r.
\]
These estimates give the asserted bounds, including the logarithmic borderline cases by taking $q$ strictly below the corresponding exponent.

Since $h$ is harmonic, $h(z)-h(0)=\langle\nabla h(0),z\rangle+O(M|z|^2)$. Consequently, if $t\leq1$, then
\[
|H(z)|\leq C_qM|z|^q,\qquad 0<q<t.
\]
If $t>1$, there is an $\ell\in\operatorname{span}\{\cos\vartheta,\sin\vartheta\}$, with $\|\ell\|_{\mathcal C^{k+2,\alpha}}\leq CM$, such that
\[
|H(x,\vartheta)-x\ell(\vartheta)|\leq C_qMx^q,
\qquad 1<q<\min\{2,t\}.
\]

On dyadic annuli $x\asymp r$, the operator $\Delta_b$ becomes the uniformly elliptic cylinder Laplacian after rescaling. Interior Schauder estimates, applied to $H$ or $H-x\ell$ as appropriate, therefore upgrade these bounds to
\[
\|H\|_{x^q\mathcal C_b^{k+2,\alpha}(\Omega')}\leq C_qM
\quad (t\leq1,\ 0<q<t),
\]
and
\[
\|H-x\ell\|_{x^q\mathcal C_b^{k+2,\alpha}(\Omega')}\leq C_qM
\quad (t>1,\ 1<q<\min\{2,t\}).
\]
Choosing $q>a$, respectively $q>b$, and applying Lemma~\ref{lem:decay-implies-little} gives the asserted estimates in the little H\"older spaces. When $t>1$, the first estimate also follows from $H=x\ell+O_b(x^q)$, since $x\ell\in x^ac_b^{k+2,\alpha}$ for $a\leq1$. Finally, the expansion with $q>1$ shows that $\ell=\lim_{x\to0}x^{-1}H(x,\cdot)$, proving uniqueness.
\end{proof}

\begin{lemma}[Cylinder inverse]\label{lem:cylinder-inverse}
Let $N\in\{1,2\}$, $s\in\mathbb R\setminus N\mathbb Z$, $m\in\mathbb N_0$, and $0<\alpha<1$. On the flat cylinder $Z_N=\mathbb R_T\times(\mathbb R/(2\pi/N)\mathbb Z)_\vartheta$, the operator
\[
 L_s=\partial_T+i\partial_\vartheta+s
\]
is an isomorphism between the bounded complex-valued H\"older spaces $\mathcal C^{m+1,\alpha}(Z_N)$ and $\mathcal C^{m,\alpha}(Z_N)$, with
\begin{equation}\label{eq:reconstruction-cylinder-estimate}
 \|L_s^{-1}B\|_{\mathcal C^{m+1,\alpha}(Z_N)}
 \leq C\|B\|_{\mathcal C^{m,\alpha}(Z_N)}.
\end{equation}
If $B$ belongs to $c_b^{m,\alpha}$ near $x=0$, where $x=e^T$, then $L_s^{-1}B$ belongs to $c_b^{m+1,\alpha}$ on a smaller collar, and
\[
 \operatorname{Tr}(L_s^{-1}B)
 =(s+i\partial_\vartheta)^{-1}\operatorname{Tr}B.
\]
\end{lemma}
\begin{proof}
Put $j=\lfloor s/N\rfloor+1$. The kernel
\begin{equation}\label{eq:reconstruction-cylinder-kernel}
 K_{s,N}(T,\vartheta)
 =\frac{N}{2\pi}
   \frac{e^{(Nj-s)T+iNj\vartheta}}
        {e^{N(T+i\vartheta)}-1}
\end{equation}
has Fourier multiplier $(i\xi+s-n)^{-1}$, $n\in N\mathbb Z$, and hence is a fundamental solution for $L_s$. Its singularity is $1/(2\pi(T+i\vartheta))$, and it decays exponentially at both ends. Thus $K_{s,N}\in L^1(Z_N)$, and convolution gives a bounded solution with $\|L_s^{-1}B\|_\infty\leq C\|B\|_\infty$. The local Schauder estimate~\cite[Lemma~2.2.4]{Biquard2008}, applied on fixed-radius disks and translated along the cylinder, gives \eqref{eq:reconstruction-cylinder-estimate}. Every homogeneous Fourier mode is a multiple of $e^{(n-s)T}$, which is unbounded at one end; this proves uniqueness.

For the little-space assertion, choose a smooth cutoff $\eta$ equal to one near the left end and supported inside the collar on which $B$ is little. Set
\[
 B_0=\operatorname{Tr}B,\qquad
 A_0=(s+i\partial_\vartheta)^{-1}B_0.
\]
The periodic first-order operator is invertible because $s\notin N\mathbb Z$, and $A_0\in c^{m+1,\alpha}$. The datum
\[
 \eta B-L_s(\eta A_0)
\]
has zero trace and vanishes near the right end. Lemma~\ref{lem:little-cutoffs}, followed by smooth approximation on compact cylinders, approximates it in $\mathcal C^{m,\alpha}(Z_N)$ by smooth compactly supported data. Their inverses are smooth and decay exponentially, with all derivatives, at the left end; by Lemma~\ref{lem:decay-implies-little}, they belong there to the zero-trace little space. Estimate~\eqref{eq:reconstruction-cylinder-estimate} passes this conclusion to the limit. Finally, the kernel formula shows that $L_s^{-1}((1-\eta)B)$ is smooth and exponentially decaying on a smaller left collar. Adding this term and $\eta A_0$ proves the assertion and the trace formula.
\end{proof}

\begin{lemma}\label{lem:compactness}
If $\delta'>\delta$, $m>k+2$, and $\alpha'>\alpha$, then the following inclusion is compact: $$x^{\delta'}\mathcal C_b^{m,\alpha'}(\Omega)\hookrightarrow x^\delta c_b^{k+2,\alpha}(\Omega) \, .$$
\end{lemma}
\begin{proof}
Choose $\chi\in\mathcal{C}^{\infty}([0,\infty))$ with $\chi=0$ on $[0,1]$ and $\chi=1$ on $[2,\infty)$. For $0<\eta<\varepsilon/2$, set $\chi_\eta(x)=\chi(x/\eta)$ and $T_\eta H=\chi_\eta H$. For fixed $\eta$, the $b$-H\"older norms are equivalent to ordinary Hölder norms on $[\eta,\varepsilon]\times S^1$. The compact inclusion $\mathcal C^{m,\alpha'}\hookrightarrow c^{k+2,\alpha}$ on this cylinder makes $T_\eta:x^{\delta'}\mathcal C_b^{m,\alpha'}(\Omega)\to x^\delta c_b^{k+2,\alpha}(\Omega)$ compact. The Leibniz rule and the scaling of $\chi_\eta$ give
\begin{equation*}
 \|H-T_\eta H\|_{x^\delta\mathcal C_b^{k+2,\alpha}(\Omega)}
 \leq C\eta^{\delta'-\delta}
       \|H\|_{x^{\delta'}\mathcal C_b^{m,\alpha'}(\Omega)}.
\end{equation*}
Since $x^\delta c_b^{k+2,\alpha}(\Omega)$ is closed in $x^\delta\mathcal C_b^{k+2,\alpha}(\Omega)$, we have $H\in x^\delta c_b^{k+2,\alpha}(\Omega)$. The same inequality shows that $T_\eta$ converges in operator norm to the inclusion, which is therefore compact.
\end{proof}

\subsection{The space of conical metrics}\label{sec:space_conical}
Choose a smooth $\Gamma$-invariant metric $h_{1/2}^\circ$ on $\Sigma\setminus\{c\}$, conformal to the round metric $\g_{S^2}$ on the sphere $S_{\Sigma}^2$, equal to $|z|^{-1}|\dd z|^2$ on $\Omega^\circ$ and to $\g_{S^2}$ near $\partial\Sigma$. For $\beta\in I$, set $h_\beta^\circ=x^{2\beta-1}h_{1/2}^\circ$. On $\Omega^\circ$,
\begin{equation}\label{eq:cone-model}
 h_\beta^\circ=x^{2\beta}\g_b
 =\dd r_\beta^2+\beta^2r_\beta^2\dd\vartheta^2,
 \qquad r_\beta=\frac{x^\beta}{\beta}.
\end{equation}
We also write $x$ for $\pi^*x$. Near $\mathcal T$, the $c_b^{m,\alpha}$-norm of a tensor is the sum of the norms of its coefficients in the cone coframe $(\dd r_\beta,\beta r_\beta\dd\vartheta,\dd\vartheta_1,\dd\vartheta_2)$. Set $K=\{p\in S^4:x(p)\geq\varepsilon/2\}$.

\begin{definition}
\phantomsection\label{notloc:tensor-space}
Let $\mathcal S_{\beta}^{m,\alpha}$ consist of the $\mathsf G$-invariant tensors $H\in c_{\mathrm{loc}}^{m,\alpha}(S^4\setminus\mathcal T;\operatorname{Sym}^2T^*S^4)$ which, on principal orbits, have the form $$H=a \pi^*h_{\beta}^{\circ} +b_1\dd\vartheta_1^2+b_2\dd\vartheta_2^2,$$ with $a-a_0,b_1-b_0,b_2-b_0\in x^\delta c_b^{m,\alpha}(\Omega)$ for constants $a_0,b_0\in\mathbb{R}$. 

For $\beta\in I$ and $t\in\mathbb{R}$, tensors $H = \pi^{\ast} H_{\Sigma} + b_1 \dd \vartheta_1^2 + b_2 \dd \vartheta_2^2$ with $H_\Sigma-a_0h_\beta^\circ,b_1-b_0,b_2-b_0\in x^tc_b^{m,\alpha}(\Omega)$ carry the norm
\begin{align}\label{eq:conical-tensor-norm}
 \|H\|_{\beta,t}^{m,\alpha}
 &=|a_0|+|b_0|
   +\|x^{-t}(H_\Sigma-a_0h_\beta^\circ)\|_{c_b^{m,\alpha}(\Omega)}\nonumber\\
 &\quad+\sum_{j=1}^2\|x^{-t}(b_j-b_0)\|_{c_b^{m,\alpha}(\Omega)}
   +\|H\|_{c^{m,\alpha}(K)}.
\end{align}
The constants $a_0,b_0$ are unique for $t>0$ and are set to zero for $t\leq0$. 
\end{definition}

\begin{lemma}
    For $m\in\mathbb{N}_0$ and $\beta\in I$, the space $(\mathcal S_\beta^{m,\alpha},\|\cdot\|_{\beta,\delta}^{m,\alpha})$ is Banach.
\end{lemma} \begin{proof}
The map
\begin{eqnarray*}
 H&\longmapsto&\bigl(a_0,b_0,x^{-\delta}(a-a_0),x^{-\delta}(b_1-b_0),x^{-\delta}(b_2-b_0),H|_K\bigr)
\end{eqnarray*}
embeds $\mathcal S_\beta^{m,\alpha}$ into the Banach space $\mathbb{R}^2\times(c_b^{m,\alpha}(\Omega))^3\times c^{m,\alpha}(K;\operatorname{Sym}^2T^*S^4)$, with equivalent norms. Agreement on the overlap, the prescribed tensor form, and $\mathsf G$-invariance are closed conditions in this product. Its image is therefore closed, so $\mathcal S_\beta^{m,\alpha}$ is Banach.
\end{proof}

To construct conical toric Einstein metrics on $S^4$ with cone angle $2\pi\beta$ along $\mathcal T$ we use metrics satisfying the following diagonal Ansatz:
\begin{definition}\label{def:X_beta}
Define
\phantomsection\label{notloc:metric-space}
 $\mathcal X_\beta =\{\g\in\mathcal S_\beta^{k+2,\alpha}:\g \, > 0 \text{ and}\ a_0,b_0>0\}$.
\end{definition}

The map $\mathcal I_\beta : \mathcal S_\beta^{m,\alpha} \to \mathcal S_{1/2}^{m,\alpha}$
\begin{equation}\label{eq:angle-identification}
 \mathcal I_\beta H=\pi^*(x^{1-2\beta}H_\Sigma)
       +b_1\dd\vartheta_1^2+b_2\dd\vartheta_2^2
\end{equation}
is a linear isomorphism and  preserves positivity: $x=1$ near $\partial\Sigma$ and $x^{1-2\beta}h_\beta^\circ=h_{1/2}^\circ$. 
Endow
\[
    \mathcal X_I:=\bigsqcup_{\beta\in I}\{\beta\}\times\mathcal X_\beta, \qquad I = [\beta_0, 1/2],
\] 
with the topology for which $(\beta,\g)\mapsto(\beta,\mathcal I_\beta\g)$ is a homeomorphism onto $I\times\mathcal X_{1/2}$. Notice that, under this identification, $\mathcal X_I$ is relatively open in $I\times \mathcal S_{1/2}^{k+2,\alpha}$.

\subsection{The scalar Einstein equations}\label{sec:reduced_equations}

\phantomsection\label{notloc:scalar-potentials}
For $\g\in\mathcal X_\beta$, put $v=\frac12\log a$, $f=\frac12\log(b_1b_2)$, and $u=\frac12\log(b_1/b_2)$ on $\Sigma^\circ\setminus\{c\}$. Then
\begin{equation}\label{eq:conical-ansatz}
 \g=\pi^*(e^{2v}h_\beta^\circ)
       +e^{f+u}\dd\vartheta_1^2+e^{f-u}\dd\vartheta_2^2,
\end{equation}
with $f-f_0,u,v-v_0\in x^\delta c_b^{k+2,\alpha}(\Omega)$, where $f_0=\log b_0$ and $v_0=\frac12\log a_0$.

The following is included here for the reader's convenience. Necessity follows immediately from Proposition~\ref{prop:quotient-equations} and Remark~\ref{rem:diagonal-reduction}, whereas sufficiency will actually be proved in a slightly more general setting in $\S$\ref{sec:fredholm} (Proposition \ref{prop:scalar-zero-equivalence}).

\begin{proposition}[Reduced Einstein equations]\label{prop:reduced-equations}
Given $\g\in\mathcal X_\beta$ and $\hat\g=e^{2v}h_\beta^\circ$, the Einstein equation $\operatorname{Ric}(\g)=\lambda\g$ is equivalent to the  system
\begin{align}
 \Delta_{\hat\g}f+|\dd f|_{\hat\g}^2+2\lambda&=0,
 \label{eqn:red-einstein-1}\\
 \Delta_{\hat\g}u+\langle\dd f,\dd u\rangle_{\hat\g}&=0,
 \label{eqn:red-einstein-2}\\
 4K_{\hat\g}+|\dd f|_{\hat\g}^2-|\dd u|_{\hat\g}^2&=0
 \label{eqn:red-einstein-3}
\end{align}
on $\Sigma^\circ\setminus\{c\}$.
\end{proposition}

Near $c$, conformal change and $\Delta_b\log x=0$ give $\Delta_{\hat\g}=e^{-2v}x^{-2\beta}\Delta_b$, and $K_{\hat\g}=-e^{-2v}x^{-2\beta}\Delta_bv$. The scalar equations therefore become
\begin{equation}\label{eq:scalar-b-system}
 \begin{aligned}
 \Delta_bf&=-|\dd f|_b^2-2\lambda e^{2v}x^{2\beta},\\
 \Delta_bu&=-\langle\dd f,\dd u\rangle_b,\\
 \Delta_bv&=\tfrac14\bigl(|\dd f|_b^2-|\dd u|_b^2\bigr).
 \end{aligned}
\end{equation}

\subsection{Regularity and branched covers}

We now address two regularity issues: smoothness of solutions in the standard smooth structure away from the conical torus, and improved asymptotic regularity near the torus. These then imply removability of the singularity on branched covers.

\begin{lemma}\label{lem:ambient-einstein-regularity}
An Einstein metric $\g \in \mathcal X_\beta$ is smooth with respect to the standard smooth structure of $S^4\setminus\mathcal T$.
\end{lemma}

\begin{proof}
By DeTurck--Kazdan~\cite[Theorem~5.2]{DeTurckKazdan1981}, $\g$ is smooth in harmonic coordinates. Let $\mathscr A_\g$ denote the resulting smooth structure, which is $\mathcal C^1$-compatible with the standard one. We show that these two smooth structures agree.

Consider the involutions $\kappa_j$ given by complex conjugation in the $j$th coordinate of $\mathbb C^2$. Since $\g$ is torus-invariant and diagonal, they are isometries. Thus the $\TT$-action and the $\kappa_j$ are smooth in $\mathscr A_\g$ by Myers--Steenrod~\cite[Chapter~II, Theorem~1.2]{KobayashiGroups}. Their common fixed surface is the meridional sphere $S_\Sigma^2$, away from its intersection with $\mathcal T$, and is smooth in both structures. Its induced metric is conformal to the prescribed smooth metric $h_\beta^\circ$, which is itself conformal to the round metric. Hence the identity between the two smooth structures on this surface is a $\mathcal C^1$ conformal map, and therefore smooth by regularity of conformal maps between Riemann surfaces.

On principal orbits, smoothness of the ambient identification follows immediately from equivariance. Near a singular orbit, choose equivariant slice coordinates in both structures, adapted to the common fixed surface. The transition map restricts to a smooth map on the signed section coordinates and commutes with the reflections $r_j\mapsto-r_j$. At a vertex, its components therefore have the form
\[
r_j\longmapsto r_j A_j(r_1^2,r_2^2),\qquad j=1,2,
\]
with $A_j$ smooth. Equivariance then extends this map to the ambient complex slice coordinates as
\[
z_j\longmapsto z_j A_j(|z_1|^2,|z_2|^2),
\]
which is smooth. The same argument, with one normal coordinate and an additional tangential coordinate, applies along an open edge. Thus the two ambient smooth structures agree everywhere on $S^4\setminus\mathcal T$, and $\g$ is smooth in the standard structure.
\end{proof}


\begin{lemma}\label{lem:cone-regularity}
If $\g\in\mathcal X_\beta$ is Einstein, then $f-f_0,u,v-v_0\in r_\beta^2\mathcal C_b^{m,\alpha}$ near $c$, for every $m\in\mathbb N_0$ (see $\S$\ref{sec:reduced_equations}). In particular, its curvature is bounded near $\mathcal T$, and $K_{\hat\g}=O(r_\beta^2)$.
\end{lemma}
\begin{proof}
Lemma~\ref{lem:ambient-einstein-regularity} implies smoothness off $\mathcal T$. Rescaled interior Schauder estimates for \eqref{eq:scalar-b-system} \cite[Theorem~6.2]{GilbargTrudinger} imply that $f-f_0,u,v-v_0 \in x^{\delta} \mathcal{C}_b^{m,\alpha}$ for all $m$. To improve the decay without losing the endpoint $2\beta=1$, remove the quadratic term in the first equation by setting $\rho=e^f$:
\begin{equation*}
 \Delta_b\rho=-2\lambda e^{2v}x^{2\beta}\rho
 =-2\lambda e^{2v_0}\rho_0x^{2\beta}+ R, \qquad R \in x^{2\beta+\delta}\mathcal{C}_b^{m,\alpha},
 \qquad \rho_0=e^{f_0}.
\end{equation*}
Since $\Delta_bx^{2\beta}=4\beta^2x^{2\beta}$, the Laplacian of $\rho-\rho_0+\lambda e^{2v_0}\rho_0x^{2\beta}/(2\beta^2)$ is contained in $x^{2\beta+\delta} \mathcal{C}_b^{m,\alpha}$. Set $H=\rho-\rho_0+\lambda e^{2v_0}\rho_0x^{2\beta}/(2\beta^2)$. For every $m\in\mathbb{N}_0$, Lemma~\ref{lem:decay-implies-little} gives
\begin{equation*}
 H\in x^{\delta/2}c_b^{m+2,\alpha},
 \qquad \Delta_bH\in x^{2\beta+\delta/2}c_b^{m,\alpha}
\end{equation*}
on a smaller collar. Proposition~\ref{prop:fourier-estimate}, with $(s,a,t)=(\delta/2,2\beta,2\beta+\delta/2)$, therefore gives $f-f_0\in x^{2\beta}c_b^{m+2,\alpha}$.  Now $\Delta_bu \in x^{2\beta+\delta} \mathcal{C}_b^{m,\alpha}$, so the same weights improve $u$ to order $x^{2\beta}$. The last equation then gives $\Delta_bv \in x^{4\beta} \mathcal{C}_b^{m,\alpha}$; use $(s,a,t)=(\delta/2,2\beta,3\beta)$ to obtain the same decay for $v-v_0$. Thus $f - f_0, u, v-v_0 \in r_{\beta}^2 \mathcal{C}_b^{m,\alpha}$ for every $m$, including $\beta=\frac{1}{2}$. Their cone-frame first derivatives are $O(r_\beta)$ and their second derivatives are bounded. Equation~\eqref{eqn:red-einstein-3} gives $K_{\hat\g}=O(r_\beta^2)$. The curvature formulas for the diagonal warped product give
\begin{equation*}
 |\operatorname{Rm}(\g)|_\g
 \leq C\bigl(|K_{\hat\g}|+|\hat\nabla^2f|_{\hat\g}
              +|\hat\nabla^2u|_{\hat\g}
              +|\dd f|_{\hat\g}^2+|\dd u|_{\hat\g}^2\bigr),
\end{equation*}
so the preceding derivative bounds give bounded curvature.
\end{proof}

The following result is closely related to the general orbifold regularity theorem of LeBrun \cite[Proposition~1.1]{LeBrun2015EdgesOrbifolds}, which asserts that Einstein edge-cone metrics with cone angle $2\pi/d$ are smooth orbifold Einstein metrics.

\begin{proposition}\label{prop:removability}
If $d\in\mathbb{N}$, $1/d\in I$, and $\g\in\mathcal X_{1/d}$ is Einstein, then $p_{d-1}^*\g$ extends uniquely to a smooth positive-definite $\TT$-invariant Einstein metric on $\widetilde X_{d-1} \simeq \sharp_{d-1} (S^2 \times S^2)$.
\end{proposition}

\begin{proof}
Choose the local coordinates in the branched-cover construction so that $z=w^d$, where $z$ is the conformal coordinate fixed in Section~5.1. Since
\[
|z|^{2/d-2}|\dd z|^2=d^2|\dd w|^2,
\]
we have
\[
p_{d-1}^*\g=d^2e^{2v}|\dd w|^2
+e^{f+u}\dd\vartheta_1^2+e^{f-u}\dd\vartheta_2^2,
\]
where the pullbacks of $f,u,v$ are denoted by the same symbols. Lemma~\ref{lem:cone-regularity}, together with $r_{1/d}=d|w|$, shows that these functions differ from their constant limits by elements of $|w|^2\mathcal C_b^{m,\alpha}$ for every $m$. In particular, they extend across $w=0$ as $\mathcal C^{1,1}$ functions, and the metric extends continuously and positive-definitely.

Under $z=w^d$, the scalar system~\eqref{eq:scalar-b-system} becomes
\[
\begin{aligned}
\Delta_w f&=-|\nabla_w f|^2-2\lambda d^2e^{2v},\\
\Delta_w u&=-\langle\nabla_w f,\nabla_w u\rangle,\\
\Delta_w v&=\tfrac14\bigl(|\nabla_w f|^2-|\nabla_w u|^2\bigr).
\end{aligned}
\]
The first derivatives of $f,u,v$ are $O(|w|)$, so integration by parts on $\{|w|>\eta\}$ produces boundary terms of order $O(\eta^2)$. Letting $\eta\to0$, the extensions satisfy the system weakly across $w=0$. Its right-hand sides are locally Lipschitz, so standard interior Schauder regularity gives $f,u,v\in\mathcal C_{\mathrm{loc}}^{2,\alpha}$ for every $\alpha<1$. Bootstrapping yields smoothness.

Away from the ramification torus, smoothness follows from Lemma~\ref{lem:ambient-einstein-regularity}, since $p_{d-1}$ is a local diffeomorphism. The extended metric is Einstein by continuity, and is $\TT$-invariant by construction. Uniqueness follows from density of the complement of the ramification torus.
\end{proof}

\section{Fredholm theory}\label{sec:fredholm}

The aim of this section is to formulate the Einstein equation as a smooth Fredholm problem of index zero, suitable for the degree argument. We work with the three scalar equations \eqref{eqn:red-einstein-1}--\eqref{eqn:red-einstein-3} from $\S$\ref{sec:conical}. The smoothness conditions at the singular orbits (not the interior cone singularity) represent a problem, as these are infinitely many constrains to our function spaces. To overcome this issue, in $\S$\ref{sec:scalar_eqn_function_spaces} we  parametrise metrics by triples $(p,q,v)$ in an enlarged Banach space $\widetilde{\mathcal M}$ which temporarily allows conical singularities along the axes in $\partial \Sigma$, and formulate the Einstein equations as a smooth map between Banach spaces
\[
\mathcal R_\beta :\widetilde{\mathcal M}\longrightarrow\mathcal F.
\]
There is a codimension-one subspace \(\mathcal M\subset\widetilde{\mathcal M}\) such that every solution of the Einstein equation in \(\mathcal M\) is smooth along \(\partial\Sigma\), while retaining the prescribed conical singularity in the interior. 
In $\S$\ref{sec:reg_comp}, we compare the metric regularity with that of the scalar functions. Using the contracted Bianchi identity, in $\S$\ref{sec:bianchi} we reconstruct the missing components of the Einstein equation, and identify a natural codimension-one subspace $\mathcal Z_{\beta,\g} \subset \mathcal F$. $\S$\ref{sec:index} establishes the Fredholm and index zero properties for $\mathcal R_\beta$ in the ambient spaces, and finally in $\S$\ref{sec:norm_Fredholm} we prove those properties for $\mathcal R_\beta$ viewed as a section of the Banach bundle with fibres $\mathcal Z_{\beta,\g}$ over $I\times \mathcal M$.



We keep the notation $k\geq2$, $0<\alpha<1$, $I=[\beta_0,1/2]$, and $0<\delta<\min\{1,2\beta_0\}$ from $\S$\ref{sec:conical}.

\subsection{The scalar equations and function spaces}\label{sec:scalar_eqn_function_spaces}

In this section we set up the Einstein equation as a smooth map between fixed Banach spaces, in preparation for computing its Fredholm index. As mentioned above, the main difference from $\S$\ref{sec:conical} is that we need to enlarge the space of metrics slightly, so that the three scalar equations can be treated as an elliptic system without imposing boundary conditions along $\partial \Sigma$. Geometrically, this amounts to allowing the metric to have conical singularities along the edges of $\Sigma$. 

Let $S:=S_\Sigma^2$ be the global section of Section~4, with round metric $h_0$. Write $\eta_j$ for reflection in $r_j=0$, and retain the involutions $\sigma,\tau$ defined there. Set $p_\pm:=(0,0,\pm1)$. The four lifts of the conical point form
\[
 C:=\{(\pm2^{-1/2},\pm2^{-1/2},0)\}.
\]
\phantomsection\label{notloc:scalar-spaces}
Let $\mathcal H_\delta^{m,\alpha}$ consist of real-valued functions $a$ invariant under $\eta_1,\eta_2,\tau$ and satisfying
\[
 a\in c_{\mathrm{loc}}^{m,\alpha}(S\setminus C),
 \qquad a\in x^\delta c_b^{m,\alpha}\quad\text{near each point of }C.
\]
Use the $c^{m,\alpha}$ norm on a fixed compact complement and the collar norms of Section~5. Fix a smooth cutoff $\zeta$, invariant under $\eta_1,\eta_2,\sigma,\tau$ and equal to one near $C$, and set
\[
 \mathcal U^{m,\alpha}
   :=\mathbb R\zeta\oplus\mathcal H_\delta^{m,\alpha}
\]
with the direct-sum norm. Different choices of $\zeta$ give equivalent norms.

On the interior of $\Sigma$, set
\[
 f_*:=\log(r_1r_2),\qquad u_*:=\log(r_1/r_2),
 \qquad f=f_*+p,\quad u=u_*+q.
\]
On the reflected quadrants, use $\dd f_*=\dd r_1/r_1+\dd r_2/r_2$ and $\dd u_*=\dd r_1/r_1-\dd r_2/r_2$.
\phantomsection\label{notloc:scalar-domain}
Define the following enlarged space of scalar metric coefficients:
\[
 \widetilde{\mathcal M}
 :=\{(p,q,v)\in(\mathcal U^{k+2,\alpha})^3:
       p\circ\sigma=p,\ q\circ\sigma=-q,\ v\circ\sigma=v\}.
\]
The condition on $q$ forces its constant term at $C$ to vanish. For each $\beta\in I$, these functions determine the metric on principal orbits
\begin{equation}\label{eq:scalar-enlarged-metric}
 \g=e^{2v}h_\beta^\circ
       +r_1^2e^{p+q}\dd\vartheta_1^2
       +r_2^2e^{p-q}\dd\vartheta_2^2.
\end{equation}
We allow the normal and angular coefficients to differ along the axes, and define the codimension-one subspace
\begin{equation}\label{eq:scalar-normalization}
 \mathcal M:=\{(p,q,v)\in\widetilde{\mathcal M}:
                      p(p_+)-2v(p_+)=0\}.
\end{equation}
Since $q(p_+)=0$ and $h_\beta^\circ=h_0$ near the axes, this sets both limiting axis angles at $p_+$ to $2\pi$.

Write $h_{1/2}^\circ=b h_0$ and set
\[
 w:=x/b,
 \qquad m=m_{\beta,v}:=e^{2v}x^{2\beta}.
\]
Then $w=1$ near the axes and, with $\hat\g=e^{2v}h_\beta^\circ$,
\begin{equation}\label{eq:scalar-fixed-weight}
 x^{2\beta}\Delta_{h_\beta^\circ}=w\Delta_{h_0},
 \qquad
 m=w\frac{\hat\g}{h_0}.
\end{equation}
Near $c$, writing $h_0=\nu(z)|\dd z|^2$ gives $w=\nu(z)x^2$ and $w\Delta_{h_0}=\Delta_b$.

Denote the left-hand sides of \eqref{eqn:red-einstein-1}--\eqref{eqn:red-einstein-3} by
\begin{equation}\label{eqn:Fi}
\begin{gathered}
 F_1 =\Delta_{\hat\g}f+|\dd f|_{\hat\g}^2+2\lambda,\\ 
 F_2 =\Delta_{\hat\g}u+\langle\dd f,\dd u\rangle_{\hat\g},\\
 F_3 =4K_{\hat\g}+|\dd f|_{\hat\g}^2-|\dd u|_{\hat\g}^2. 
\end{gathered}
\end{equation}
Set
\begin{equation}\label{eqn:F_target}
 \mathcal F
 :=\{(R_1,R_2,R_3)\in(\mathcal H_\delta^{k,\alpha})^3:
       R_1\circ\sigma=R_1,\ R_2\circ\sigma=-R_2,
       \ R_3\circ\sigma=R_3\}.
\end{equation}

The Einstein operator in this augmented domain is then given by
\begin{equation}\label{eq:scalar-residual-map}
\mathcal R_\beta : \widetilde{\mathcal M} \to \mathcal F, \qquad
 \mathcal R_\beta(p,q,v)
   :=(mF_1,mF_2,-mF_3/4).
\end{equation}
Its domain and ambient target are independent of $\beta$. Near $c$,
\begin{equation}\label{eq:scalar-residual-end}
 \mathcal R_\beta=
 \left(\Delta_bf+|\dd f|_b^2+2\lambda m,
       \Delta_bu+\langle\dd f,\dd u\rangle_b,
       \Delta_bv-\tfrac14(|\dd f|_b^2-|\dd u|_b^2)\right).
\end{equation}

\begin{lemma}\label{lem:scalar-map-smooth}
The map
\[
\mathcal R:I\times\widetilde{\mathcal M}\longrightarrow\mathcal F,\qquad (\beta,p,q,v)\longmapsto\mathcal R_\beta(p,q,v)
\]
is smooth.
\end{lemma}

\begin{proof}
Near $c$, use the multiplication estimates of Section~5 and $2\beta>\delta$: every $\beta$ derivative of $x^{2\beta}$ belongs to $x^\delta c_b^{k,\alpha}$. Along an axis, the terms of order $r_j^{-2}$ cancel. The remaining quotients have the form $r_j^{-1}\partial_{r_j}a$ with $a$ even, and
\begin{equation}\label{eq:scalar-even-division}
\frac{\partial_ra(r,y)}r=\int_0^1\partial_r^2a(tr,y)\,\dd t.
\end{equation}
This defines a bounded linear map from even $c^{k+2,\alpha}$ functions to $c^{k,\alpha}$ functions. Apply the same formula successively in $r_1,r_2$ at a vertex. The strict inequality $2\beta_0>\delta$ also gives a smooth extension to an open interval containing $I$.
\end{proof}

\subsection{Regularity comparison}\label{sec:reg_comp}

We now show that, at least when the smoothness conditions at $\partial \Sigma$ are satisfied, regularity of the three scalar coefficients is equivalent to regularity of the corresponding metric tensor. In $\S$\ref{sec:norm_Fredholm} we will subsequently show that solutions indeed satisfy the smoothness conditions, and thus belong to the original geometric space $\mathcal X_\beta$ defined in $\S$\ref{sec:conical}.

\begin{lemma}\label{lem:scalar-axis-regularity}\label{lem:even-angular-extension}
Let $m\in\mathbb N_0$ and $0<\alpha<1$. Near an axis or vertex, a tensor
\[
 T=A h_0+r_1^2B_1\dd\vartheta_1^2+r_2^2B_2\dd\vartheta_2^2
\]
belongs to ambient $c^{m,\alpha}$ if and only if $A,B_1,B_2\in c^{m,\alpha}$ on the signed section, are even in each normal variable, and satisfy
\begin{equation}\label{eq:axis-conicity-traces}
 B_j=A\quad\text{on }r_j=0
\end{equation}
for each axis present. The scalar and ambient norms are equivalent on fixed smaller neighborhoods. More generally, radial substitution $a(r,y)\mapsto a(|x|,y)$, $x\in\mathbb R^d$, $d\geq2$, is bounded on even $c^{m,\alpha}$ functions, including all tangential derivatives. These assertions also hold for the big H\"older spaces.
\end{lemma}
\begin{proof}
Set $\widehat F(x,y):=F(|x|,y)x_ix_j/|x|^2$, where $F\in\mathcal C^{m,\alpha}$ is even and $F(0,y)=0$. The total-degree-$m$ Taylor polynomial at $(0,y_0)$ consists of terms $r^{2j}(y-y_0)^\gamma$, with $j\geq1$; their extensions are the Cartesian polynomials $x_ix_j|x|^{2j-2}(y-y_0)^\gamma$. Taylor's remainder estimate and the chain rule therefore give, jointly in all variables,
\[
 |D^\ell(\widehat F-\widehat P_{y_0})(x,y_0)|
 \leq C\|F\|_{\mathcal C^{m,\alpha}}
                  |x|^{m-|\ell|+\alpha},\qquad |\ell|\leq m.
\]
Thus the derivatives extend to the indicated polynomial jets on the axis. If two points are separated by at least a fixed fraction of the larger normal radius, compare their top derivatives through the axis jets, whose order-$m$ coefficients are $\alpha$-H\"older in $y$. Otherwise subtract the Taylor polynomial at the first point's axis projection and rescale both normal and tangential variables by its normal radius. On this fixed annulus the remainder has bounded $\mathcal C^{m,\alpha}$ norm, while the top derivative of the extended polynomial is constant. This proves the bounded angular extension without derivative loss. Allowing $j=0$ proves the scalar assertion.

For little spaces, take smooth even approximants and replace each angular coefficient $F_n$ by $F_n(r,y)-F_n(0,y)$. Trace and constant extension are bounded; these approximants converge and their extensions are smooth. The preceding estimate gives convergence in $\mathcal C^{m,\alpha}$.

At a vertex use the round-sphere chart $(x,z)\in\mathbb R^2\times \mathbb R^2$, with $r_1=|x|$, $r_2=|z|$, and denote its round metric by $g_0$. Put $\Pi(x)=|\dd x|^2-(x\cdot\dd x)^2/|x|^2$. Then
\begin{equation}\label{eq:axis-tensor-reconstruction}
 T=A g_0+(B_1-A)\Pi(x)+(B_2-A)\Pi(z).
\end{equation}
Apply the angular extension in one normal variable and the scalar extension in the other, treating all remaining coordinates as parameters. This proves the forward bound, uniformly in $|x|/|z|$. Conversely, restriction to $(x,z)=((r,0),(s,0))$ gives
\[
 B_1=T_{x_2x_2},\qquad B_2=T_{z_2z_2},\qquad
 A=T_{x_1x_1}/(g_0)_{x_1x_1}.
\]
The denominator is smooth and positive, so these formulas give the reverse bound. Rotations give evenness and equality of radial and angular components on an axis, proving \eqref{eq:axis-conicity-traces}.

Along an open axis only one angular extension is needed. Recover its coefficient from a normal Cartesian component, the other from the nonvanishing circle component divided by $r_j^2$, and $A$ from a tangential component of $T|_{S}=Ah_0$. All denominators are smooth and positive. These restriction operations also preserve little spaces.
\end{proof}

For metrics, $A=e^{2v}$, $B_1=e^{p+q}$, and $B_2=e^{p-q}$. Logarithms and exponentials are locally Lipschitz in $c^{m,\alpha}$ on bounded sets of positive coefficients bounded away from zero. Lemma~\ref{lem:scalar-axis-regularity} therefore gives the same local norm bounds for $(p,q,v)$ and the ambient metric. The axis conditions are $p+q=2v$ on $r_1=0$ and $p-q=2v$ on $r_2=0$. Since $h_\beta^\circ=h_0$ near the axes, the bounds are independent of $\beta$.

\subsection{The Bianchi identity and the restricted target space}\label{sec:bianchi}

The three scalar Einstein equations determine all components of the Einstein tensor except its trace-free horizontal part. The contracted Bianchi identity allows us to reconstruct this missing component from the residuals $F_i$'s from \eqref{eqn:Fi}.Its necessary vanishing at the vertices imposes a single compatibility condition, defining the natural target space $\mathcal Z_{\beta,\g}$ in \eqref{eq:scalar-corner-functional} for our Fredholm problem.

Let $\mathcal G$ be the finite group generated by $\eta_1,\eta_2,\sigma,\tau$, and define its character $\epsilon$ by
\[
 \epsilon(\eta_1)=\epsilon(\eta_2)=-1,
 \qquad \epsilon(\sigma)=\epsilon(\tau)=1.
\]
Let $\mathcal Q_\delta^{k,\alpha}$ consist of trace-free covariant two-tensors $Q$ satisfying $\gamma^*Q=\epsilon(\gamma)Q$ for $\gamma\in\mathcal G$, with
\[
 Q\in c_{\mathrm{loc}}^{k,\alpha}
       (S\setminus C;\operatorname{Sym}^2T^*S)
\]
whose coefficients in smooth conformal coordinates near $C$ belong to $x^{\delta-2}c_b^{k,\alpha}$. Define $\mathcal A_\delta^{k-1,\alpha}$ similarly for one-forms, with coefficients in $x^{\delta-3}c_b^{k-1,\alpha}$ near $C$.

\begin{lemma}\label{lem:scalar-reconstruction}
The operator
\[
 D=\operatorname{div}_{h_0}:
 \mathcal Q_\delta^{k,\alpha}
 \longrightarrow\mathcal A_\delta^{k-1,\alpha}
\]
is a bounded isomorphism.
\end{lemma}
\begin{proof}
In a conformal coordinate $z=y_1+iy_2$, write $h_0=e^{2\omega}|\dd z|^2$ and $Q=\operatorname{Re}(a\,\dd z^2)$. Then
\begin{equation}\label{eq:reconstruction-local-divergence}
 DQ=2e^{-2\omega}\operatorname{Re}
       ((\partial_{\bar z}a)\dd z).
\end{equation}
This formula gives boundedness in the stated weighted norms.

We first prove injectivity. At each point of $C$, the stabilizer contains a half-turn $z\mapsto-z$ with character $1$, so $a$ is even. If $DQ=0$, then $a$ is holomorphic off $C$ and satisfies $a=O(x^{\delta-2})$. Since $0<\delta<1$, its only possible poles are simple, and evenness excludes them. Thus $Q$ extends to a holomorphic quadratic differential on $S$. Such a differential vanishes: its coefficient in a stereographic coordinate is entire and is $O(|z|^{-4})$ at infinity.

For surjectivity, let $\upsilon\in\mathcal A_\delta^{k-1,\alpha}$. Near a point of $C$, equation $DQ=\upsilon$ becomes
\[
 \partial_{\bar z}a=G,\qquad
 G=\tfrac12e^{2\omega}(\upsilon_1-i\upsilon_2).
\]
The half-turn makes $G$ odd. With $s=\delta-2$, $T=\log x$, and $a=x^sA$, this is
\begin{equation}\label{eq:reconstruction-cylinder-equation}
 (\partial_T+i\partial_\vartheta+s)A
 =2e^{-i\vartheta}x^{1-s}G.
\end{equation}
Both sides are $\pi$-periodic, and the right-hand side belongs to $c_b^{k-1,\alpha}$. Since $s\in(-2,-1)$, Lemma~\ref{lem:cylinder-inverse} with $N=2$ applies. Cut off the datum, extend it by zero to $Z_2$, solve, and cut off the resulting tensor in the conformal disk. Doing this at all points of $C$ and averaging with character $\epsilon$ gives $Q_e\in\mathcal Q_\delta^{k,\alpha}$ such that $\upsilon-DQ_e$ extends to a $c^{k-1,\alpha}$ one-form on $S$ and vanishes near $C$.

On the unpunctured sphere, the kernel of $D$ is zero by the same holomorphic argument. Its adjoint kernel consists of conformal Killing fields. Anti-invariance under both axis reflections excludes the gradient fields and leaves only rotation about the $t$-axis; this rotation is odd under $\sigma$, so it does not have character $\epsilon$. The elliptic Fredholm alternative~\cite[equation~(2.2.5)]{Biquard2008} therefore solves
\[
 DQ_0=\upsilon-DQ_e
\]
in character $\epsilon$, with a bounded inverse on ordinary H\"older spaces and, by smooth approximation, on little H\"older spaces. Since the datum vanishes near $C$, $Q_0$ is smooth there. Its normalized coefficients $x^{2-\delta}Q_0$ have positive decay and belong to the little space by Lemma~\ref{lem:decay-implies-little}. Hence $Q_e+Q_0\in\mathcal Q_\delta^{k,\alpha}$ solves $DQ=\upsilon$. All these constructions are bounded in the indicated norms, proving the assertion.
\end{proof}

On the reflected sphere, define the signed orbit-area density
\[
 \varrho:=r_1r_2e^p.
\]
For $R=(R_1,R_2,R_3)\in\mathcal F$, set
\begin{equation}\label{eq:scalar-reconstruction-forcing}
 J_{\beta,\g}(R)
  :=-\frac{\varrho}{2w}
     \left(\dd R_1-R_1\dd\log m+R_2\dd u-2R_3\dd f\right),
 \qquad Q_R:=D^{-1}J_{\beta,\g}(R).
\end{equation}
The one-form $J_{\beta,\g}(R)$ has character $\epsilon$ and belongs to $\mathcal A_\delta^{k-1,\alpha}$. Near $C$, its coefficient weight is $x^{\delta-3}$, since
\[
 w^{-1}=O(x^{-2}),\qquad
 \dd\log m=2\dd v+2\beta\,\dd x/x.
\]
At the axes, $\varrho\,\dd f$ and $\varrho\,\dd u$ extend as little H\"older one-forms in smooth coordinates. Thus $R\mapsto Q_R$ is bounded and depends smoothly on $(\beta,p,q,v)$.

Choose a fixed conformal coordinate $(s,t)$ near $p_+$ with axes $s=0,t=0$ and $\sigma(s,t)=(t,s)$. For a tensor $Q$ with character $\epsilon$, either axis reflection makes its diagonal coefficients vanish at $p_+$. Its mixed coefficient is even under both reflections and need not vanish. Thus $Q(p_+)$ has one possible real component. Since $\tau$ determines $Q(p_-)$ from $Q(p_+)$, vanishing at the two vertices imposes one condition.

Define
\begin{equation}\label{eq:scalar-corner-functional}
 \kappa_{\beta,\g}(R):=(Q_R)_{st}(p_+),
 \qquad \mathcal Z_{\beta,\g}:=\ker\kappa_{\beta,\g}.
\end{equation}
Rescaling the coordinate multiplies $\kappa$ by a nonzero constant and leaves its kernel unchanged.

\begin{proposition}\label{prop:scalar-target}
The functional $\kappa_{\beta,\g}$ is nonzero for every $(\beta,p,q,v)\in I\times\widetilde{\mathcal M}$. Its kernels form a smooth codimension-one subbundle of the trivial bundle with fiber $\mathcal F$, and
\[
 \mathcal R_\beta(p,q,v)\in\mathcal Z_{\beta,\g}.
\]
\end{proposition}
\begin{proof}
Let $E:=\operatorname{Ric}(\g)-\lambda\g$ and let $H$ be the trace-free part of its horizontal block. The vertical eigenvalues and horizontal trace are
\[
 E_1=-\tfrac12(F_1+F_2),\qquad E_2=-\tfrac12(F_1-F_2),\qquad
 \operatorname{tr}_{\hat\g}E_H=\tfrac12F_3-F_1.
\]
The contracted Bianchi identity gives
\begin{equation}\label{eq:scalar-bianchi}
 \operatorname{div}_{\hat\g}(\varrho H)
 =-\frac{\varrho}{2}
       (\dd F_1+F_2\dd u+\tfrac12F_3\dd f).
\end{equation}
Since $\operatorname{div}_{\hat\g}Q=(\hat\g/h_0)^{-1}DQ$ for a trace-free covariant tensor, \eqref{eq:scalar-fixed-weight} gives $D(\varrho H)=J_{\beta,\g}(\mathcal R_\beta(\g))$.

Near $p_+$, write
\[
 \hat\g=e^{2\omega}(\dd s^2+\dd t^2),\qquad
 \ell_1=se^a,\quad\ell_2=te^b,\qquad \varrho=st e^{a+b},
\]
with $a,b,\omega\in c^{k+2,\alpha}$ even in both variables. The horizontal Ricci formula for a diagonal metric (see \cite[\S~4.2, equation~(4.6)]{Lott2010}) gives 
\[
 E_{st}=\frac{\omega_t-a_t}{s}+\frac{\omega_s-b_s}{t}+B,
 \qquad B\in c^{k,\alpha},
\]
where $B$ is odd in both variables. Its diagonal coefficients belong to $c^{k,\alpha}$ by \eqref{eq:scalar-even-division}. Consequently
\begin{equation}\label{eq:scalar-corner-vanishing}
 (\varrho H)_{ss}=st\,C_{ss},\qquad
 (\varrho H)_{st}=e^{a+b}
       \{t(\omega_t-a_t)+s(\omega_s-b_s)+stB\},
 \qquad C_{ss}\in c^{k,\alpha}.
\end{equation}
Both coefficients extend in $c^{k,\alpha}$ and vanish at the vertex, without requiring equality of the normal and angular coefficients. The same formula gives extension along an open axis. Near $C$, their weight is $x^{\delta-2}$, so $\varrho H\in\mathcal Q_\delta^{k,\alpha}$. Lemma~\ref{lem:scalar-reconstruction} identifies it with $Q_{\mathcal R_\beta(\g)}$, proving $\kappa_{\beta,\g}(\mathcal R_\beta(\g))=0$.

To prove nonvanishing, take $Q=2\psi(s^2+t^2)\dd s\dd t$, where $\psi=1$ near the vertex and has small compact support. Extend with character $\epsilon$ and by $\tau$ near $p_-$. Then $Q_{st}(p_+)=1$ and $DQ$ is supported in annuli about the vertices. Set $R_1=0$ and solve
\[
 R_2\dd u-2R_3\dd f=-2w\varrho^{-1}DQ.
\]
Multiplying the $\dd s$ equation by $s$ and the $\dd t$ equation by $t$ gives an even $c^{k+1,\alpha}$ coefficient matrix whose value at the vertex is $\left(\begin{smallmatrix}1&-2\\-1&-2\end{smallmatrix}\right)$. It is invertible on a sufficiently small disk. The right-hand side extends evenly across both axes: the $\dd s$ coefficient of $DQ$ is divisible by $t$, and its $\dd t$ coefficient by $s$. Thus $R_2,R_3\in c^{k+1,\alpha}$ are supported in the annuli. Uniqueness gives their required exchange symmetries. Extension by zero and by $\tau$ gives $R\in\mathcal F$ with $J(R)=DQ$. Hence $Q_R=Q$ and $\kappa_{\beta,\g}(R)=1$.

Smooth dependence and non-vanishing give local smooth projections: if $\kappa(W)\ne0$, use $R\mapsto R-\kappa(R)W/\kappa(W)$.
\end{proof}

\subsection{Fredholm and index in the ambient space}\label{sec:index}

In order to compute the index, we compare the linearisation of the three scalar Einstein equations with invariant Laplacians on the round spheres $S^2$, $S^4$, and $S^6$, whose orbit-volume densities are $1$, $r_1r_2$, and $(r_1r_2)^2$, respectively.

\begin{lemma}\label{lem:scalar-reference-index}
For $a\in\{0,1,2\}$, set
\[
 P_a^*=w\bigl(\Delta_{h_0}
                 +a\langle\dd f_*,\dd\,\cdot\,\rangle_{h_0}\bigr).
\]
The restriction
\[
 P_a^*:
 \{u\in\mathcal U^{k+2,\alpha}:u\circ\sigma=u\}
 \longrightarrow
 \{F\in\mathcal H_\delta^{k,\alpha}:F\circ\sigma=F\}
\]
is Fredholm of index zero. Its kernel consists of the constants, and its range is the kernel of
\begin{equation}\label{eq:scalar-mean-functional}
 L_a(F)=\int_\Sigma \frac{F}{w}(r_1r_2)^a\,\dd v_{h_0}.
\end{equation}
The restriction
\[
 P_a^*:
 \{u\in\mathcal U^{k+2,\alpha}:u\circ\sigma=-u\}
 \longrightarrow
 \{F\in\mathcal H_\delta^{k,\alpha}:F\circ\sigma=-F\}
\]
is an isomorphism.
\end{lemma}
\begin{proof}
For $a=0,1,2$, consider $M_a=S^{2a+2}\subset\mathbb R^{a+1}\times\mathbb R^{a+1}\times\mathbb R$. The $O(a+1)\times O(a+1)$ quotient is $(\Sigma,h_0)$, with orbit-volume density proportional to $(r_1r_2)^a$. Consequently
\[
 \Delta_{M_a}u=\Delta_{h_0}u
                  +a\langle\dd f_*,\dd u\rangle_{h_0}.
\]
For $a=0$ these are the two reflections on $S^2$. The scalar assertion of Lemma~\ref{lem:even-angular-extension} identifies the even little H\"older spaces at the axes with the invariant spaces on $M_a$.

Near the codimension-two preimage $T_a$ of $c$,
\[
 P_a^*=\Delta_b+a\{c_1\partial_T+c_2\partial_\vartheta\},
 \qquad (c_1,c_2)=(\partial_Tf_*,\partial_\vartheta f_*).
\]
Half-turn invariance gives $c_j=O_b(x^2)$ with all derivatives. On $Z_1=\mathbb R_T\times(\mathbb R/2\pi\mathbb Z)_\vartheta$,
\[
 x^{-\delta}\Delta_bx^\delta
 =(\partial_T-i\partial_\vartheta+\delta)
  (\partial_T+i\partial_\vartheta+\delta).
\]
The inverse from Lemma~\ref{lem:cylinder-inverse}, and its reflection in $\vartheta$, therefore give an inverse $G_\delta$ gaining two derivatives, since $0<\delta<1$ excludes all integer roots. It preserves the little spaces at the left end and has trace operator $(\delta^2+\partial_\vartheta^2)^{-1}$.

Cut the drift coefficients off in $x<\varepsilon$ and extend them by zero to the right end. Their conjugated perturbation $E_\varepsilon$ satisfies
\[
 \|E_\varepsilon\|_{
 \mathcal C^{k+2,\alpha}(Z_1)\to\mathcal C^{k,\alpha}(Z_1)}
       \leq C\varepsilon^2.
\]
Here the cutoffs are dilates of one fixed cutoff, so their $b$ derivatives are uniformly bounded. For sufficiently small $\varepsilon$, $(1+G_\delta E_\varepsilon)^{-1}G_\delta$ exists by its Neumann series and preserves the little spaces. Apply it to cut-off $x^{-\delta}F$, multiply the solution by $x^\delta$, and cut off again. Averaging with the prescribed symmetries gives $u_e\in\mathcal H_\delta^{k+2,\alpha}$ such that
\[
 H:=(F-P_a^*u_e)/w
\]
extends to an invariant $c^{k,\alpha}$ function on $M_a$ and vanishes near $T_a$. All cutoffs equal one sufficiently near $c$.

The flux of $u_e$ across a tube of radius $\eta$ is $O(\eta^\delta)$, since its normal derivative is $O(\eta^{\delta-1})$ and the boundary measure is $O(\eta)$. Thus $L_a(F)=0$ makes $H$ mean zero. Ordinary Poisson theory on $M_a$ gives a mean-zero solution $\Delta_{M_a}u_0=H$ with the same symmetries. It is smooth near $T_a$; invariance gives one common limiting constant, zero when $F\circ\sigma=-F$. Its remaining terms are $O(x)$, hence lie in $x^\delta c_b^{k+2,\alpha}$. Therefore $u_e+u_0$ belongs to the required domain and solves the equation, with a bounded norm.

The integral $L_a$ is bounded because its local density is $O(x^{\delta-1})\,\dd x\,\dd\vartheta$. It is nonzero on $\sigma$-invariant data and vanishes on $\sigma$-anti-invariant data. The same flux argument shows that it annihilates the range. Finally, a kernel element is bounded; integration against a test function gives boundary errors $O(\eta^\delta)+O(\eta)$. It is weakly harmonic across $T_a$, hence smooth and constant on $M_a$. Only the invariant restriction admits constants. These kernel and range descriptions give both assertions.
\end{proof}

\begin{proposition}\label{prop:scalar-ambient-index}
At every point of $\widetilde{\mathcal M}$, the derivative
\[
 D\mathcal R_\beta:
       \widetilde{\mathcal M}\longrightarrow\mathcal F
\]
is Fredholm of index zero.
\end{proposition}
\begin{proof}
For $a=0,1,2$, set
\[
 P_a=w\bigl(\Delta_{h_0}
                    +a\langle\dd f,\dd\,\cdot\,\rangle_{h_0}\bigr).
\]
The difference $P_a-P_a^*=aw\langle\dd p,\dd\,\cdot\,\rangle$ is compact on the spaces of Lemma~\ref{lem:scalar-reference-index}. Away from $C$, it is first order with $c^{k+1,\alpha}$ coefficients; compactness follows from $c^{k+1,\alpha}\hookrightarrow c^{k,\alpha}$. Near $C$, it is $a\langle\dd_b p,\dd_b\,\cdot\,\rangle_b$, where $\dd_b p\in x^\delta c_b^{k+1,\alpha}$. It annihilates the constant term at the end. On the decaying part of the domain, its restriction to $x<\varepsilon$ has operator norm $O(\varepsilon^\delta)$. It is therefore a norm limit of compact cutoff operators. Each $P_a$ has index zero on either symmetry subspace.

Write $(\phi,\psi,\chi)=(\dot p,\dot q,\dot v)$. Differentiating \eqref{eq:scalar-residual-map} gives
\begin{equation}\label{eq:scalar-triangular-linearization}
 D\mathcal R_\beta
 \begin{pmatrix}\phi\\\psi\\\chi\end{pmatrix}
 =\begin{pmatrix}
 P_2\phi+4\lambda m\chi\\
 P_1\psi+w\langle\dd u,\dd\phi\rangle_{h_0}\\
 P_0\chi-\frac12w\langle\dd f,\dd\phi\rangle_{h_0}
              +\frac12w\langle\dd u,\dd\psi\rangle_{h_0}
 \end{pmatrix}.
\end{equation}
The off-diagonal first-order terms are bounded. At the axes, this follows from \eqref{eq:scalar-even-division} for the factors $r_j^{-1}\partial_{r_j}$; near $C$, it follows from the weighted product estimates.

Omitting $4\lambda m\chi$ gives a lower triangular operator. Its diagonal blocks $P_2,P_1,P_0$ are Fredholm of index zero. A triangular parametrix gives index equal to their sum, hence zero. Multiplication by $4\lambda m$ is compact from the $\sigma$-invariant subspace of $\mathcal U^{k+2,\alpha}$ to that of $\mathcal H_\delta^{k,\alpha}$. Its restriction to the constant term at the end has finite rank and belongs to the target because $2\beta>\delta$. On the decaying part, it gains weight $2\beta$. Its restriction to $x<\varepsilon$ tends to zero in operator norm, and its compactly supported cutoffs are compact. Adding this term therefore preserves the index.
\end{proof}

\subsection{The Einstein section and its properties}\label{sec:norm_Fredholm}

The next result shows that every zero of $\mathcal R_\beta$ in $\mathcal M$ determines an Einstein metric that extends smoothly across the singular orbits $\pi^{-1}(\partial\Sigma)$.

\begin{proposition}\label{prop:scalar-zero-equivalence}
For a metric $\g$ represented by a point of $\mathcal M$, $\mathcal R_\beta(\g)=0$ if and only if $\g\in\mathcal X_\beta$ and $\operatorname{Ric}(\g)=\lambda\g$. At a zero, the $\ker D \mathcal R_\beta$ identifies with the infinitesimal Einstein deformations in $T_\g\mathcal X_\beta$.
\end{proposition}
\begin{proof}
If $\mathcal R_\beta(\g)=0$, then \eqref{eq:scalar-reconstruction-forcing} and Lemma~\ref{lem:scalar-reconstruction} give $\varrho H=0$. The vertical equations and the horizontal trace equation also hold, so $\g$ is Einstein on principal orbits.

At an open axis $s=0$, equation \eqref{eq:scalar-corner-vanishing} gives
\[
 (\varrho H)_{st}|_{s=0}
       =(\varrho/s)|_{s=0}\,\partial_t(\omega-a).
\]
Thus $a-\omega$ is constant along the axis, whose cone angle is $2\pi e^{a-\omega}$. The same holds on the other axis. Exchange symmetry equates these constants, and \eqref{eq:scalar-normalization} sets their common value to zero at $p_+$. Both axes have angle $2\pi$. Lemma~\ref{lem:scalar-axis-regularity}, with $m=k+2$, gives $c^{k+2,\alpha}$ regularity of the ambient metric at the axes and vertices. Its expansion at the conical point is the one prescribed in $\mathcal X_\beta$. The regularity results of Section~5 apply. Conversely, an Einstein metric in $\mathcal X_\beta$ has the required scalar symmetries, satisfies the normalization by smoothness at $p_+$, and has $\mathcal R_\beta(\g)=0$.

At an Einstein metric, suppose $D\mathcal R_\beta(\phi,\psi,\chi)=0$. Differentiating $D(\varrho H)=J(\mathcal R_\beta)$ and applying Lemma~\ref{lem:scalar-reconstruction} gives $\varrho\dot H=0$. Differentiating the vertical and trace identities gives the remaining linearized Einstein equations. The differentiated axis identity makes the variations of $a-\omega$ constant along the axes. Symmetry and the differentiated normalization set them to zero. The linear correspondence and norm bounds in Lemma~\ref{lem:scalar-axis-regularity} therefore place the metric variation in $T_\g\mathcal X_\beta$. The converse follows by differentiating the scalar expressions.
\end{proof}

\begin{theorem}\label{thm:scalar-index-zero}
The maps $\mathcal R_\beta$ define a smooth section of the Banach bundle with fibers $\mathcal Z_{\beta,\g}$ over $I\times\mathcal M$. At every zero, its vertical derivative
\[
 D\mathcal R_\beta:T_\g\mathcal M
                       \longrightarrow\mathcal Z_{\beta,\g}
\]
is Fredholm of index zero. Its zero set identifies with the conical Einstein metrics in $\mathcal X_I$, with the topology from $\S$\ref{sec:space_conical}.
\end{theorem}
\begin{proof}
Smoothness and the bundle assertion follow from Proposition~\ref{prop:scalar-target}. At a zero, differentiating $\kappa_{\beta,\g}(\mathcal R_\beta(\g))=0$ gives
\[
 \kappa_{\beta,\g}\circ D\mathcal R_\beta=0.
\]
By Proposition~\ref{prop:scalar-ambient-index}, the derivative on $\widetilde{\mathcal M}$ with target $\mathcal F$ has index zero. Replacing the target by the codimension-one subspace $\mathcal Z_{\beta,\g}$ increases the index by one. Restricting the domain by \eqref{eq:scalar-normalization} decreases it by one. Hence
\[
 \operatorname{ind}
   (D\mathcal R_\beta:T_\g\mathcal M\to\mathcal Z_{\beta,\g})
       =0+1-1=0.
\]
The change $p\mapsto p+C$ rescales both orbit lengths by $e^{C/2}$ and leaves the three scalar expressions unchanged. Its derivative is transverse to \eqref{eq:scalar-normalization}.

Proposition~\ref{prop:scalar-zero-equivalence} identifies the zeros. Near the axes and vertices, Lemma~\ref{lem:scalar-axis-regularity} and the logarithmic change of variables give locally bounded maps between the scalar and ambient norms. On the remaining principal region, restriction and reconstruction are smooth operations. Near $C$, logarithms and exponentials of the positive coefficients are continuous in the weighted little norms of Section~5. These bounds hold under \eqref{eq:angle-identification}, which keeps the scalar coefficients fixed. The identification is therefore a homeomorphism, including for the families over $I$.
\end{proof}

\begin{remark}\label{rem:scalar-degree-projections}
Let $K\subset\mathcal M$ be compact. Since $\kappa_{\beta,\g}\ne0$, compactness of $I\times K$ gives finitely many $W_j\in\mathcal F$ such that $\sum_j\kappa_{\beta,\g}(W_j)^2>0$ on $I\times K$. Compactness of $I$ gives an open neighborhood $U\subset\mathcal M$ of $K$ on which this inequality holds for every $\beta\in I$. On $I\times U$, set
\[
 Z_{\beta,\g}
   =\frac{\sum_j\kappa_{\beta,\g}(W_j)W_j}
          {\sum_j\kappa_{\beta,\g}(W_j)^2},
 \qquad
 \Pi_{\beta,\g}R=R-\kappa_{\beta,\g}(R)Z_{\beta,\g}.
\]
Then $\Pi_{\beta,\g}$ is a smooth projection onto $\mathcal Z_{\beta,\g}$. Take $K$ to be the projection onto $\mathcal M$ of the total zero set, compact by Section~7 and the identification above. Proposition~\ref{prop:localized-degree} applies with $E=\mathcal M$, $F=\mathcal F$, $s_\beta=\mathcal R_\beta$, and $P_{\beta,\g}=\Pi_{\beta,\g}$.
\end{remark}

\section{Compactness}\label{sec:compactness}

The main goal of this section is to establish the compactness of the set of conical toric Einstein metrics as the cone angles $2\pi \beta$ vary within a compact  interval: $\beta \in I \Subset (0,1/2]$. This is one of the key requirements in order to define a homotopy-invariant mod-two degree (see Proposition \ref{prop:localized-degree}).

To that end, one first establishes geometric convergence of a sequence, keeping track of the conical singularities. This is done by obtaining curvature bounds via a blow-up analysis. Once these bounds are available, volume non-collapsing  follows immediately from the energy identity  \eqref{eqn:L2-bound}. One then has to translate this geometric convergence into convergence in the prescribed coordinates and weighted topology as per $\S\S$\ref{sec:conical}--\ref{sec:fredholm}.

\subsection{Geometric  a priori  estimates}

\begin{lemma}\label{lem:global-bounds}
Let $\beta\in I$ and let $\g\in\mathcal X_\beta$ be Einstein. Its metric completion has diameter at most $\pi\sqrt{3/\lambda}$,  $\operatorname{Vol}(\g)\leq24\pi^2/\lambda^2$, and
\begin{eqnarray}
\label{eqn:L2-bound}
 \int_{S^4\setminus\mathcal T}|{\rm Rm}(\g)|_\g^2\dd v_\g &=&64\pi^2.
\end{eqnarray}
\end{lemma}
\begin{proof}
By Lemma~\ref{lem:cone-regularity}, the coefficients of $\g-\g_0$ in the cone frames lie in $r^2\mathcal C_b^{m,\alpha}$ for every $m$, where $\g_0$ is the product-cone model. Since $2\pi\beta\leq2\pi$, the metric completion is an $\text{RCD}(\lambda,4)$ space \cite[Theorem~A]{BKMR21}.  Bonnet--Myers gives the diameter bound. Bishop comparison centered at a regular point $p$, where $\operatorname{Vol}_\g(B_\g(p,r))/r^4\to\pi^2/2$ as $r\to0$, gives the volume bound \cite[Corollary~2.6 and Theorem~2.3]{Sturm06}. In smooth normal disc coordinates $z=xe^{i\vartheta}$, define
\begin{equation*}
 \Phi(z,y):= \bigl((\beta e^{-v_0})^{1/\beta}|z|^{1/\beta-1}z,y\bigr),  \qquad y\in\mathcal T,
\end{equation*}
with $\Phi(0,y)=(0,y)$. Extend $\Phi$ to a homeomorphism of $S^4$, a diffeomorphism off $\mathcal T$, equal to the identity outside a collar. Then $(e^{v_0}x^\beta/\beta)\circ\Phi=|z|$ near $\mathcal T$. Apply Chern--Gauss--Bonnet to $\Phi^*\g$. By Lemma~\ref{lem:cone-regularity}, $\Phi^*\g$ is an edge-cone metric in the sense of \cite[\S1]{AtiyahLeBrun}. Since $\Phi$ is a diffeomorphism off $\mathcal T$ and a homeomorphism of $(S^4,\mathcal T)$, the curvature integral and both Euler characteristics are unchanged. Since $\chi(\mathcal T)=0$, the correction $(1-\beta)\chi(\mathcal T)$ vanishes. The Einstein equation gives $|\text{Ric}|_\g^2=4\lambda^2$ and $\text{Scal}=4\lambda$, hence $-4|\text{Ric}|_\g^2+\text{Scal}^2=0$. Thus \cite[Theorem~2.1]{AtiyahLeBrun}, written using the full tensor norm, gives \eqref{eqn:L2-bound}.
\end{proof}

\begin{lemma}\label{lem:quotient-curvature}
The quotient metric $\hat\g$ has Gaussian curvature $K<\lambda/3$ on $\Sigma^\circ\setminus\{c\}$.
\end{lemma}
\begin{proof}
Set $\mu:=K-\lambda/3$ and $F:=\rho^{3/2}\mu$. Equations~\eqref{eqn:edge-slice-expansion}--\eqref{eqn:vertex-slice-expansion} show that $F$ extends continuously to $\partial\Sigma$ with value zero. Lemma~\ref{lem:cone-regularity} gives $K(c)=0$, so $F$ also extends continuously across $c$, with $F(c)=-\lambda\rho(c)^{3/2}/3<0$. A positive maximum of $F$ would occur in $\Sigma^\circ\setminus\{c\}$, contradicting \eqref{eqn:weyl-smooth-scalar}. Hence $F\leq0$. The strong maximum principle gives either $F<0$ or $F\equiv0$ on $\Sigma^\circ\setminus\{c\}$. Continuity at $c$ and $F(c)<0$ exclude the second alternative.
\end{proof}

\begin{lemma}\label{lem:density-maximum}
The orbit-area density $\rho$ attains its maximum at $c$, and $\dd\rho\neq0$ on $\Sigma^\circ\setminus\{c\}$.
\end{lemma}
\begin{proof}
Along every geodesic starting at $c$, $\rho$ decreases strictly until the geodesic reaches $\partial\Sigma$. By \eqref{eqn:general-horizontal} and Lemma~\ref{lem:quotient-curvature}, $\nabla^2f\leq-(2\lambda/3)\hat\g$. Let $\gamma:[0,l]\to\Sigma$ be a minimizing geodesic from $c$ to a point in $\Sigma^\circ\setminus\{c\}$, parametrized by arc length. It cannot meet $\partial\Sigma$, since the edges are totally geodesic and the corners have angle $\pi/2$. By Lemma~\ref{lem:cone-regularity}, $(f\circ\gamma)'(0)=0$. Integrating $(f\circ\gamma)''\leq-2\lambda/3$ gives
\begin{equation*}
\dd f(\dot\gamma(l))\leq-2\lambda l/3<0,\qquad f(\gamma(l))\leq f(c)-\lambda l^2/3<f(c).
\end{equation*}
Since $\rho=e^f$ in the interior and $\rho=0$ on $\partial\Sigma$, these inequalities prove both assertions.
\end{proof}

\begin{lemma}\label{lem:density-isoperimetry}
Let $\g\in\mathcal X_\beta$ be Einstein and let $0<s<\rho(c)$. Set $A(s):=\operatorname{Area}_{\hat\g}\{\rho>s\}$ and $L(s):=\operatorname{Length}_{\hat\g}\{\rho=s\}$. Then $\{\rho>s\}$ contains $c$, is homeomorphic to an open disk, and satisfies
\begin{eqnarray}
\label{eq:density-isoperimetric}
 L(s)^2&\geq&4\pi\beta A(s)-\frac{\lambda}{3}A(s)^2.
\end{eqnarray}
\end{lemma}
\begin{proof}
Set $E:=\{\rho>s\}$. Since $s>0$ and $\rho=0$ on $\partial\Sigma$, $\overline E$ is a compact subset of $\Sigma^\circ$. Each component of $E$ contains an interior maximum of $\rho$, hence contains $c$ by Lemma~\ref{lem:density-maximum}. Any component of $\Sigma\setminus E$ disjoint from $\partial\Sigma$ would contain an interior minimum away from $c$, contradicting the same lemma. Thus $E$ is a disk. Since $\dd\rho\neq0$ on $\partial E$, its boundary is a smooth simple closed curve. Choose a domain $\Omega$ with smooth boundary and $E\Subset\Omega\Subset\Sigma^\circ$, and a conformal coordinate $z$ on $\Omega$ with $z(c)=0$. Near $c$, write $\hat\g=|z|^{2\beta-2}e^{2v}|\dd z|^2$.  By Lemma~\ref{lem:cone-regularity}, $v-v_0\in |z|^{2\beta}\mathcal C_b^{m,\alpha}$ for every $m$. Hence
\begin{equation*}
 \left|\int_{|z|=\eta}\partial_\nu v\,\dd\ell\right|  \leq C\eta^{2\beta}\longrightarrow0,
\end{equation*}
where $\nu$ and $\dd\ell$ are the Euclidean unit normal and arc length. The curvature measure is therefore $K\,\dd v_{\hat\g}+2\pi(1-\beta)\delta_c$. The curvature measure has finite total variation on $\Omega$. Its only atom is at $c$, with mass $2\pi(1-\beta)<2\pi$. Let $\kappa$ denote the curvature measure above. Writing $\hat\g=e^{2u}|\dd z|^2$, we have
\begin{equation*}
 u(z)=h(z)+\frac{1}{2\pi}\int_\Omega
 \log\frac{1}{|z-\zeta|}\,\dd\kappa(\zeta),
\end{equation*}
where $h$ is harmonic. In the conventions of~\cite[Definition~4.1]{BartolucciCastorina2019}, the conformal exponent and representing measure are $2u$ and $2\kappa$. The bound $\kappa(\{c\})=2\pi(1-\beta)<2\pi$ is precisely the no-cusp condition in~\cite[Remark~4.3]{BartolucciCastorina2019}. By Lemma~\ref{lem:quotient-curvature}, $K\leq\lambda/3$, so
\begin{equation*}
 \sup_{\substack{B\subset E\\ B\ \mathrm{Borel}}} \left(\kappa(B)-\frac{\lambda}{3} \operatorname{Area}_{\hat\g}(B)\right)  =2\pi(1-\beta).
\end{equation*}
The supremum is attained at $B=\{c\}$. Applying \cite[Theorem~4.7, equation~(4.9)]{BartolucciCastorina2019} with $K_0=\lambda/3$ gives \eqref{eq:density-isoperimetric}.
\end{proof}

\begin{lemma}\label{lem:density}
Orient $E_1,E_2$ by the boundary orientation of $\Sigma$. Then
\begin{eqnarray}
 P=P(\g)&:=&\frac{\lambda\text{Vol}(\g)}{2\pi^2}
 =-2\int_{E_1} \star_{\hat\g}\dd\rho
 =-2\int_{E_2} \star_{\hat\g}\dd\rho,
 \label{eq:density-boundary}\\
 \sup_\Sigma\rho&\leq&\frac{\sqrt3}{4\pi\beta}P
 \leq\frac{\sqrt3}{4\pi\beta_0}P.
 \label{eq:density-bound}
\end{eqnarray}
\end{lemma}
\begin{proof}
Integrating~\eqref{eqn:general-trace} over $\Sigma$ and using
\begin{equation}
\label{eq:orbit-volume}
    \operatorname{Vol}(\g)
    =4\pi^2\int_\Sigma\rho\,\dd v_{\hat\g}
\end{equation}
(see~\eqref{eq:conical-ansatz}) gives
\begin{equation*}
 -\int_{\partial\Sigma}\star_{\hat\g}\dd\rho =2\lambda\int_\Sigma\rho\,\dd v_{\hat\g}=P.
\end{equation*}
By Lemma~\ref{lem:cone-regularity}, $|\dd\rho|_{\hat\g}\leq Cr$ near $c$, so $\int_{\{r=\eta\}} \star_{\hat\g}\dd\rho\to0$ as $\eta\to0$. The symmetry exchanging $E_1$ and $E_2$ gives $-\int_{E_i} \star_{\hat\g}\dd\rho=P/2$ for $i=1,2$, proving \eqref{eq:density-boundary}. Set $M:=\rho(c)$ and let $A,L$ be as in Lemma~\ref{lem:density-isoperimetry}. Define
\begin{equation*}
H(A(s)):=\int_{\{\rho>s\}}\rho\,\dd v_{\hat\g}, \qquad 0<s<M.
\end{equation*}
The coarea formula gives $A'(s)<0$ and $\frac{\dd}{\dd s}H(A(s))=sA'(s)$. Hence $H'(A(s))=s$ and $H''(A(s))=1/A'(s)<0$. Set $a:=A(s)$. Stokes' theorem gives $\int_{\{\rho=s\}}|\dd\rho|_{\hat\g}\,\dd\ell_{\hat\g} =2\lambda H(a)$. Cauchy--Schwarz and Lemma~\ref{lem:density-isoperimetry} then give
\begin{eqnarray*}
 4\pi\beta a-\frac{\lambda}{3}a^2 \leq L(s)^2 \leq2\lambda H(a)(-A'(s)) =-\frac{2\lambda H(a)}{H''(a)}.
\end{eqnarray*}
Where $4\pi\beta a-\lambda a^2/3>0$, this gives
\begin{equation*}
 H''(a)+\frac{2\lambda H(a)}  {4\pi\beta a-\lambda a^2/3}\geq0.
\end{equation*}
The extension to $a=0$ satisfies $H(0)=0$ and $H'(0)=M$. The function
\begin{eqnarray*}
 H_0(a)&=&\frac{12\pi\beta M}{\lambda}y(1-y)(1-2y),  \qquad y=\frac{\lambda a}{12\pi\beta},
\end{eqnarray*}
satisfies $H_0''+2\lambda H_0/(4\pi\beta a-\lambda a^2/3)=0$, with $H_0(0)=0$ and $H_0'(0)=M$. The derivative $H_0'$ first vanishes at $a_*:=(2\pi\beta/\lambda)(3-\sqrt3)$. Set $W:=H'H_0-HH_0'$. For $0<a<\min\{a_*,A(0)\}$, we have $H_0>0$ and $W'\geq0$. Since $W\to0$ and $H/H_0\to1$ as $a\to0$, it follows that $H\geq H_0$ and $H'\geq(H/H_0)H_0'\geq H_0'>0$. If $A(0)<a_*$, taking $a\uparrow A(0)$ would contradict $H'(A(0))=0$. Thus $A(0)\geq a_*$. Since $H$ is increasing,
\begin{equation*}
 \frac{P}{2\lambda} =H(A(0))\geq H(a_*)\geq H_0(a_*)  =\frac{2\pi\beta M}{\sqrt3\lambda}.
\end{equation*}
This gives \eqref{eq:density-bound}.
\end{proof}

\subsection{Ricci-flat bubbles}

This subsection studies smooth Ricci-flat limits arising when the conical orbit escapes to infinity under curvature rescaling. We first prove rigidity results and then analyze the quotient topology. Proposition~\ref{prop:polarbubbles} shows that a nonflat limit contains an entire compact edge, including both endpoints and its primitive stabilizer. The bounded-distance cone case is treated in \S\ref{sec:blow-up}.

\begin{lemma}\label{lem:properness-static-field}
A complete Ricci-flat four-manifold is flat if it admits a nowhere-zero Killing field whose orthogonal distribution is integrable.
\end{lemma}
\begin{proof}
Let $Z$ be the nowhere-zero Killing field. The second fundamental form of each leaf of $Z^\perp$ is symmetric and, by the Killing equation, skew-symmetric. It therefore vanishes, so each leaf is totally geodesic. An ambient geodesic initially tangent to a leaf remains in it, since $\g(Z,\dot\gamma)$ is constant. Completeness of $\g$ therefore implies completeness of each leaf. Extend local coordinates on a leaf along the flow of $Z$, with parameter $t$. Then $\g=g_3+N^2\dd t^2$, where $g_3$ is the induced leaf metric and $N:=|Z|_\g>0$. Since $Z$ is Killing, $g_3$ and $N$ are independent of $t$. Ricci flatness gives
\begin{equation*}
 \Delta_{g_3}N=0,\qquad  N\text{Ric}_{g_3}=\nabla^2_{g_3}N.
\end{equation*}
Anderson's static-vacuum estimate~\cite[Theorem~1.1]{Anderson00}, together with completeness of the leaf and positivity of $N$, gives $\text{Ric}_{g_3}=0$ and $\dd N=0$. Since a three-dimensional Ricci-flat metric is flat and $N$ is constant, $\g=g_3+N^2\dd t^2$ is flat.
\end{proof}

\begin{lemma}\label{lem:properness-marked-action}
Let $(\g_i,p_i)\to(X,\g_\infty,p)$ smoothly, with $\g_i$ in the diagonal family and uniformly bounded curvature. Suppose $X$ is complete, nonflat, and Ricci-flat, and $d_{\g_i}(p_i,\mathcal T)\to\infty$. Then each fixed circle generator vanishes within a uniformly bounded distance of $p_i$. Along a subsequence, the $T^2$-actions converge smoothly, with the same stabilizers and edge labels.
\end{lemma}
\begin{proof}
Write $K_{i,1},K_{i,2}$ for the fixed circle generators. Suppose, along a subsequence, that $d_{\g_i}(p_i,\{K_{i,a}=0\})\to\infty$ for some $a\in\{1,2\}$. On a unit ball at distance at least one from $\{K_{i,a}=0\}$, $N_{i,a}:=|K_{i,a}|_{\g_i}>0$. Along unit-speed geodesics orthogonal to $K_{i,a}$, the curvature bound gives $|N_{i,a}''|\leq C N_{i,a}$. At the center of the ball, comparison with $u''=Cu$ in both directions gives the first estimate; the diagonal form gives the second:
\begin{eqnarray}
\label{eq:positive-jacobi-bound}
|\dd\log N_{i,a}|&\leq&\sqrt C\coth\sqrt C, \qquad |\nabla K_{i,a}|^2=2|\dd N_{i,a}|^2.
\end{eqnarray}
These bounds and $|\nabla^2K_{i,a}|\leq C|\text{Rm}(\g_i)|\,|K_{i,a}|$ give, along a subsequence, smooth convergence of $K_{i,a}/N_{i,a}(p_i)$ on compact sets to a nowhere-zero Killing field with integrable orthogonal distribution. Lemma~\ref{lem:properness-static-field} then forces $\g_\infty$ to be flat, a contradiction. Each generator has a zero within a uniformly bounded distance of $p_i$. Its circle rotates the normal plane with period $2\pi$, so $|\nabla K_{i,a}|^2=2$ at its zeros. Integrating $|\nabla^2K_{i,a}|\leq C|K_{i,a}|$ along geodesics from these zeros to $p_i$ gives $c\leq |K_{i,a}(p_i)|+|\nabla K_{i,a}(p_i)|\leq C$, where $c,C>0$ are independent of $i$. At a limiting zero, the differential of the time-$t$ flow rotates the normal plane through angle $t$. Thus any positive period is a multiple of $2\pi$, so $2\pi$ is the minimal period. Let $\Phi_i$ be the convergence embeddings on relatively compact invariant domains. Following the center-of-mass construction of~\cite{GroveKarcher1973}, consider the points $a\Phi_i(a^{-1}x)$, $a\in T^2$. For large $i$, they lie in a small convex normal ball about $\Phi_i(x)$. Their center of mass for Haar probability measure defines a smooth map $F_i$. Haar invariance and the isometric action give $F_i(ax)=aF_i(x)$. Since $F_i$ is $\mathcal C^1$-close to $\Phi_i$, it is an embedding on a smaller invariant domain. Injectivity gives $\text{Stab}_{T^2}(x)=\text{Stab}_{T^2}(F_i(x))$. Thus the limiting action is effective and has the same edge labels.
\end{proof}

For related uniqueness with prescribed asymptotic control of the harmonic map, see Li--Sun~\cite[Proposition~4.38]{LS25}. The following argument obtains global half-plane coordinates from completeness and the sign of the quotient curvature.

\begin{lemma}\label{lem:properness-one-corner}
Let $(X,\g)$ be a complete diagonal Ricci-flat four-manifold with an effective $T^2$-action and no nontrivial finite stabilizers. Suppose the quotient of the free orbits is simply connected and $\partial(X/T^2)$ consists of two smooth edges meeting at one smooth corner. Let $h$ be the quotient metric. If $K_h\leq0$, then $\g$ is flat.
\end{lemma}
\begin{proof}
Set $R:=\rho$. Since $\lambda=0$ and $K_h\leq0$, equations~\eqref{eqn:general-trace} and~\eqref{eqn:general-horizontal} give $\Delta_hR=0$ and $\nabla_h^2\log R\leq0$. If $\dd R(q)=0$, then along each geodesic $\gamma$ starting at $q$, $(\log R\circ\gamma)'(0)=0$ and $(\log R\circ\gamma)''\leq0$. Thus $R\leq R(q)$ near $q$. The maximum principle makes $R$ constant, contradicting $R>0$ in the interior and $R=0$ on the boundary. Hence $\dd R\neq0$. Since $\dd(\star_h\dd R)=(\Delta_hR)\dd v_h=0$, the form $\star_h\dd R$ is closed. Simply connectedness makes it exact; choose $z$ with $\dd z=-\star_h\dd R$. In the local coordinates $(R,z)$, write $h=e^{2\nu}(\dd R^2+\dd z^2)$. Equation~\eqref{eqn:general-trace}, Remark~\ref{rem:diagonal-reduction}, and $K_h\leq0$ give the first estimate below; the horizontal equation gives the remaining identities:
\begin{eqnarray*}
u_R^2+u_z^2&\leq&R^{-2},\\
\nu_R&=&\frac R4(u_R^2-u_z^2)-\frac1{4R},\qquad \nu_z=\frac R2u_Ru_z.
\end{eqnarray*}
Set $F:=(R,z)$. Let $\gamma:[0,1]\to\mathbb R_{>0}\times\mathbb R$ be piecewise smooth and put $r_0:=\min_\gamma R>0$. If $\tilde\gamma$ is a lift of $\gamma$ defined on $[0,T)$, then for $0\leq t<T$,
\begin{eqnarray*}
|\nu(\tilde\gamma(t))-\nu(\tilde\gamma(0))|
&\leq&\frac{\sqrt5}{4r_0}
\operatorname{Length}_{\mathbb R^2}(\gamma).
\end{eqnarray*}
The upper bound for $\nu$ gives finite $h$-length for $\tilde\gamma$, so completeness of the quotient gives a limit $q$ as $t\uparrow T$. Since $R(q)=R(\gamma(T))>0$ and there are no exceptional orbits, $q$ is principal. Hence $F(q)=\gamma(T)$, and the local inverse of $F$ continues the lift. Every lift therefore extends to its endpoint. Lifting segments from a point of $F$'s image shows that $F$ is onto. For a ball $B\Subset\mathbb R_{>0}\times\mathbb R$ centered at $y$, lifting radial segments from each point of $F^{-1}(y)$ defines inverse branches over $B$. Uniqueness of lifts makes their images disjoint; lifting the reversed segments shows that they exhaust $F^{-1}(B)$. Thus $F$ is a covering of the simply connected half-plane $\mathbb R_{>0}\times\mathbb R$, and hence a diffeomorphism. For fixed $z$, integrating $\nu_R\geq-1/(2R)$ gives $e^\nu\leq CR^{-1/2}$, so $\int_0^{R_0}e^{\nu(R,z)}\dd R\leq2C\sqrt{R_0}<\infty$. Near an open edge, \eqref{eqn:edge-slice-expansion} gives $R=rD$ with $D(0,s)>0$. The closed form $\star_h\dd R$ extends smoothly to the edge, so its primitive $z$ does too. The extended functions $(R,z)$ are local coordinates there:
\begin{eqnarray*}
\left|\det\frac{\partial(R,z)}{\partial(r,s)}\right|_{r=0} &=&D(0,s)^2>0.
\end{eqnarray*}
Normalize $z$ to vanish at the corner. With the chosen sign of $\dd z$, the function $z+iR$ is holomorphic. After reflection across both edges, \eqref{eqn:vertex-slice-expansion} gives it a double zero at the corner. Choose $x+iy:=\sqrt{2(z+iR)}$ with $x,y>0$ in the interior; this defines smooth coordinates at the corner. Then $R=xy$ and $z=(x^2-y^2)/2$, with the edges on the positive coordinate axes. For the flat metric, $u_0=\log(x/y)$ and $e^{2\nu_0}=(2\sqrt{R^2+z^2})^{-1}$. The difference $w:=u-u_0$ satisfies $w_{RR}+R^{-1}w_R+w_{zz}=0$. Thus $(\xi,\eta,z)\mapsto w(\sqrt{\xi^2+\eta^2},z)$ is harmonic on $\mathbb R^3$ away from the $z$-axis. Equations~\eqref{eqn:edge-slice-expansion} and~\eqref{eqn:vertex-slice-expansion} show that $w$ extends smoothly across the $z$-axis away from the origin and is bounded at the origin. Hence $w$ extends harmonically to $\mathbb R^3$. Since $|\nabla w|\leq2/R$, the mean-value property applied to $\partial_z w$ gives
\begin{eqnarray*}
|\partial_zw(p)|\leq\frac2{|B_L|} \int_{B_L(p)}R^{-1}\leq\frac CL\longrightarrow0.
\end{eqnarray*}
Thus $w_z=0$ and $(Rw_R)_R=0$. Smoothness at $R=0$ gives $w=B$. The formulas for $\nu_R,\nu_z$ give $\nu=\nu_0+C$. In $(x,y)$ coordinates, $h=e^{2C}(\dd x^2+\dd y^2)$, $G_{11}=e^Bx^2$ and $G_{22}=e^{-B}y^2$. Smoothness at $x=0$ and $y=0$, with minimal circle periods $2\pi$, gives $e^B=e^{2C}=e^{-B}$, so $B=C=0$.
\end{proof}

\begin{proposition}\label{prop:polarbubbles}
Let $(\g_i,p_i)$ in the diagonal family, with fixed circle generators, converge smoothly to a complete nonflat Ricci-flat limit $(X,\g_\infty,p)$. Suppose $\g_\infty$ has bounded curvature, $K_h\leq0$, and $d_{\g_i}(p_i,\mathcal T)\to\infty$. Then $X/T^2$ contains a compact edge that, for large $i$, the embeddings of Lemma~\ref{lem:properness-marked-action} map onto an entire original edge, including both endpoints and with the same stabilizer circle.
\end{proposition}
\begin{proof}
By Lemma~\ref{lem:properness-marked-action}, each compact subset of $\Sigma$ embeds in the original quotient disk, preserving stabilizers and edge labels. Thus $\Sigma$ has genus zero, no nontrivial finite stabilizers, and at most two corners. Moreover, $\Sigma$ is noncompact. Otherwise $X$ would be compact, and Bochner's formula would make both generators parallel. Their independence on the free orbits gives a local product $\mathbb R^2\times N$, where $N$ is a Ricci-flat surface. Since $N$ is flat, $X$ would be flat, a contradiction. Also, $\Sigma$ has one end. Otherwise $X$ would have more than one end, since $\TT$ is compact. The splitting theorem~\cite{CheegerGromoll71} would give $X\cong\mathbb R\times Y$ isometrically, with $Y$ a Ricci-flat three-manifold. Then $Y$, and hence $X$, would be flat, a contradiction.

A circular boundary component $C\subset\partial\Sigma$ without corners would embed in an original edge interval, which is impossible. If $C$ had one corner $q$, then $C\setminus\{q\}$ would be connected with constant stabilizer, so the stabilizers on both sides of $q$ would agree, contrary to~\eqref{eqn:vertex-slice-expansion}. Thus $C$ has two corners, joined by a compact edge $E$. Otherwise every boundary component is a line. Since $\Sigma$ has genus zero and one end, each simple closed curve in $\Sigma^\circ$ bounds a compact region. That region contains no boundary line and is therefore a disk. Hence $\pi_1(\Sigma^\circ)=0$. An arc meeting $\partial\Sigma$ only at endpoints on two distinct lines would separate $\Sigma$ into two noncompact components, contradicting its single end. Thus there is at most one boundary line. If $\partial\Sigma=\varnothing$ or the line had no corners, one fixed circle generator would be nowhere zero, making $X$ flat by Lemma~\ref{lem:properness-static-field}. Lemma~\ref{lem:properness-one-corner} excludes one corner. Hence the boundary line has two corners, joined by a compact edge $E$.

Let $E_i$ be the original edge with the same label as $E$. For large $i$, Lemma~\ref{lem:properness-marked-action} gives a quotient embedding $\bar F_i$ with $\bar F_i(E)\subset E_i$ and $\bar F_i(\partial E)=\partial E_i$. Since $\bar F_i(E)$ is connected, $\bar F_i(E)=E_i$.
\end{proof}

\begin{figure}[tbp]
\centering
\includegraphics[width=\textwidth]{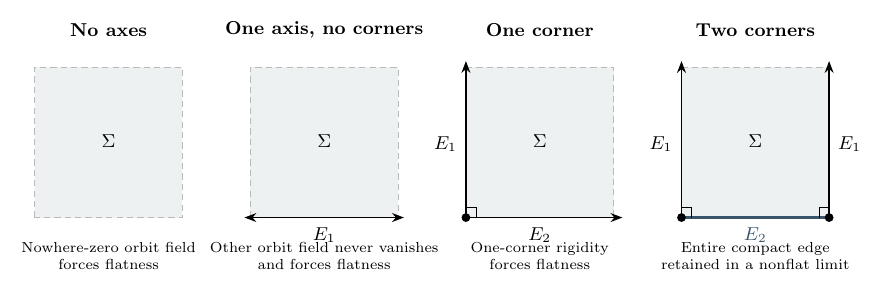}
\caption{Representative quotient boundary patterns in Proposition~\ref{prop:polarbubbles}. The first two are excluded by Lemma~\ref{lem:properness-static-field}, and the third by Lemma~\ref{lem:properness-one-corner}. The fourth illustrates the compact-edge conclusion, with the same primitive label on the two noncompact sides. These are topological schematics, not an exhaustive classification of quotient surfaces.}
\label{fig:section7-bubble-patterns}
\end{figure}

\subsection{Blow-up limits and non-collapse}\label{sec:blow-up}

We construct partially unwrapped covers on which the blow-up sequence has a complete Ricci-flat limit. The orbit directions and the reflected quotient require separate noncollapse arguments. We then show that the cone escapes and exclude the remaining deck identifications. Proposition~\ref{prop:no-collapsed-smooth-limit} gives an injectivity-radius lower bound for the rescaled metrics downstairs.

\begin{notation}[Blow-up sequence]\label{not:rescaled-covers}
Let $\g_i\in\mathcal X_{\beta_i}$ be Einstein, with $\beta_i\in I$. Choose $s_i>0$ tending to zero and regular points $p_i$. Set
\begin{eqnarray}
\label{eq:rescaled-metrics}
 \widetilde\g_i&=&s_i^{-2}\g_i,
 \qquad h_i=s_i^{-2}\hat\g_i,\nonumber\\
 \widetilde\rho_i&=&s_i^{-2}\rho_i,
 \qquad \varepsilon_i=\lambda s_i^2.
\end{eqnarray}
Then $\text{Ric}(\widetilde\g_i)=\varepsilon_i\widetilde\g_i$. Assume that
\begin{eqnarray}
\label{eq:curvature-normalization}
 |\text{Rm}(\widetilde\g_i)|_{\widetilde\g_i}&\leq&2,
 \qquad
 |\text{Rm}(\widetilde\g_i)|_{\widetilde\g_i}(p_i)=1.
\end{eqnarray}
\end{notation}

The normalization gives $|\mathrm{Rm}(\g_i)|_{\g_i}(p_i)=s_i^{-2}\to\infty$. Although the rescaled curvatures are bounded, their injectivity radii may tend to zero.

\subsubsection{Local covers and vertical non-collapse}

Short periodic orbit directions can cause collapse under a curvature bound. We unwrap the directions whose axes escape and retain the primitive $2\pi$ periods of those whose axes remain nearby: a circle shrinking at an included axis already bounds a smooth normal disk. The radii below grow slowly enough that the positive curvature allowed by $K_{h_i}\leq\varepsilon_i/3$ contributes vanishing total curvature to the comparison disks in Lemma~\ref{lem:reflected-quotient-injectivity}.

\begin{definition}[Partially unwrapped covers]\label{def:rescaled-covers}
For a sequence as in Notation~\ref{not:rescaled-covers}, put $b_i:=\pi(p_i)$ and write $c_i$ for the conical point. All distances and balls below are taken in $h_i$. Pass to a subsequence so that the distances from $b_i$ to $c_i$ and to each axis are either bounded or tend to infinity. Choose $R_i\to\infty$ so that every escaping axis or cone lies outside $B_{h_i}(b_i,2R_i)$ and
\[
\varepsilon_iV_{2,-2}(4R_i+1)\longrightarrow0,
\qquad V_{2,-2}(r):=\pi\bigl(\cosh(\sqrt2r)-1\bigr).
\]
Here $V_{2,-2}(r)$ is the area of a radius-$r$ ball in the hyperbolic plane of curvature $-2$.

Choose a connected relatively open domain $\mathcal Q_i\subset\Sigma$ with piecewise smooth boundary such that
\[
B_{h_i}(b_i,R_i)\subset\mathcal Q_i\subset B_{h_i}(b_i,2R_i).
\]
In a fixed planar realization of $\Sigma$, let $\operatorname{Fill}(\mathcal Q_i)$ be obtained by filling the bounded complementary components. Set
\[
\mathcal Q_i^+:=
\begin{cases}
\mathcal Q_i, & d_{h_i}(b_i,c_i)\to\infty,\\
\operatorname{Fill}(\mathcal Q_i)\cup\{c_i\}, & \sup_i d_{h_i}(b_i,c_i)<\infty,
\end{cases}
\qquad
\mathcal Q_i^*:=\mathcal Q_i^+\cap\bigl(\Sigma^\circ\setminus\{c_i\}\bigr).
\]
Put $U_i:=\pi^{-1}(\mathcal Q_i^+)$. Its \emph{partially unwrapped cover} $\widehat U_i$ is the metric completion of the universal cover of the underlying manifold $U_i$, including the axes and the conical orbit when present. Choose a lift $\widehat p_i$ of $p_i$ and also write $\widetilde\g_i$ for the lifted metric.
\end{definition}

Filling bounded planar components adds no new axis, since a complementary component meeting $\partial\Sigma$ connects to the exterior of the planar disk. When the cone remains nearby, $\mathcal Q_i^+$ is simply connected. A primitive orbit circle retains its $2\pi$ period precisely when its axis meets $\mathcal Q_i^+$; otherwise it lifts to a line. Thus each connected lift of a principal orbit is
\[
\begin{cases}
\mathbb R^2, & \text{both axes escape},\\
S^1\times\mathbb R, & \text{exactly one axis remains nearby},\\
\TT, & \text{both axes remain nearby}.
\end{cases}
\]
Including the conical orbit makes its transverse meridian contractible, so the lifted metric keeps cone angle $2\pi\beta_i$. It is smooth across the retained axes. On smaller domains, restrict this covering and take the relevant connected component.

By the \emph{reflected quotient} we mean the connected lift through $\widehat p_i$ of a meridional section through $p_i$, with its induced metric. It is smooth across the axes and at the vertices by Lemma~\ref{lem:slice-regularity}. When the cone escapes, each principal chamber is simply connected: pushing off the axes identifies $\pi_1(\mathcal Q_i^*)$ with $\pi_1(\mathcal Q_i)$, and the section makes its map into $\pi_1(U_i)$ injective. Its lift to the universal cover is therefore simply connected.

\begin{figure}[tbp]
\centering
\includegraphics[width=.94\textwidth]{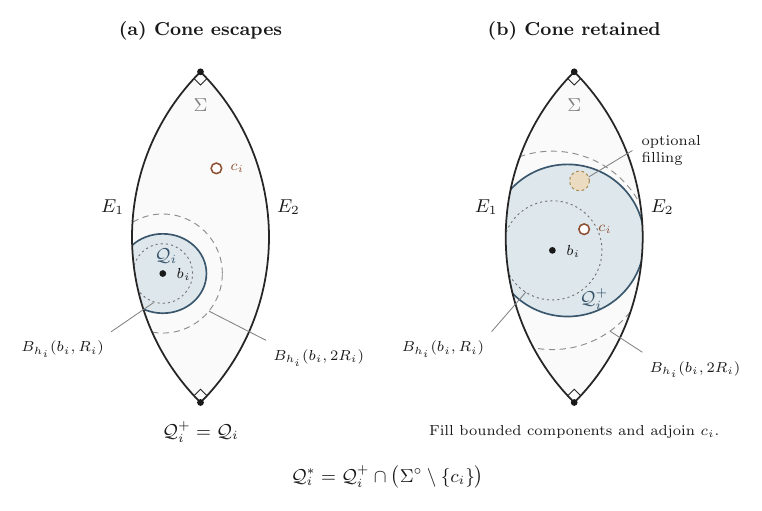}
\caption{The domains of Definition~\ref{def:rescaled-covers}, drawn in a right-angled bigon. Dotted and dashed curves indicate the balls of radii $R_i$ and $2R_i$. On the left the cone escapes and only $E_1$ meets the domain. On the right the cone lies in the inner ball; the ochre patch illustrates an optional bounded complementary component filled to form $\mathcal Q_i^+$. Both axes meet the domain in this example. The principal part $\mathcal Q_i^*$ excludes the axes and $c_i$. The balls and domains are schematic.}
\label{fig:section7-partial-unwrapping}
\end{figure}

The next result shows that a retained primitive circle can become short only by approaching its axis.

\begin{lemma}\label{lem:short-axis}
Let $K_i$ be a fixed primitive circle generator of a diagonal metric $\g_i$ with $|\mathrm{Rm}|\leq2$, and suppose $Z_i:=\{K_i=0\}$ is smooth. For a fixed $D>0$, assume $B_{\g_i}(p_i,2D+1)$ is smooth and $d_{\g_i}(p_i,Z_i)\leq D$. If $|K_i|_{\g_i}(p_i)\to0$, then $d_{\g_i}(p_i,Z_i)\to0$. If both fixed generators $K_{i,1},K_{i,2}$ satisfy these hypotheses and their lengths at $p_i$ tend to zero, then $d_{\g_i}(p_i,\{K_{i,1}=K_{i,2}=0\})\to0$.
\end{lemma}
\begin{proof}
The time-$t$ flow of $K_i$ rotates the normal plane to $Z_i$ through angle $t$, so $|\nabla K_i|^2=2$ on $Z_i$. Integrating $|\nabla^2K_i|\leq C|K_i|$ along curves of length at most $L$ starting on $Z_i$ gives
\begin{equation*}
 \sqrt2e^{-CL}
 \leq (|K_i|^2+|\nabla K_i|^2)^{1/2}
 \leq \sqrt2e^{CL},
\end{equation*}
with $C$ independent of $i$. Choose a unit-speed minimizing geodesic $\gamma_i$ from $Z_i$ to $p_i$. Set $l_i:=d_{\g_i}(p_i,Z_i)$ and $j_i(t):=|K_i(\gamma_i(t))|_{\g_i}$. The diagonal form gives
\begin{equation}
\label{eq:primitive-jacobi-data}
 j_i''+\mathrm{Sec}_{\g_i}(\dot\gamma_i,K_i)j_i=0,
 \qquad j_i(0)=0,\qquad j_i'(0)=1.
\end{equation}
Suppose $l_i\geq d>0$ and pass to a subsequence with $l_i\to l\in[d,D]$. Extend $\gamma_i$ beyond $p_i$ by $d/2$. A zero on this extension would lie within distance $d/2<l_i$ of $p_i$, contradicting the definition of $l_i$. Thus $j_i>0$ on $(0,l_i+d/2]$. The bound $|j_i''|\leq2j_i$ gives, after a further subsequence, $\mathcal C^1$ convergence on $[0,l+d/4]$ to a nonnegative function $j$. Since $j_i(l_i)\to0$, we have $j(l)=j'(l)=0$. The limit satisfies $|j''|\leq2j$ almost everywhere, so Gronwall's inequality gives $j\equiv0$, contradicting $j'(0)=1$. Choose $q_i\in Z_{i,1}$ with $d(p_i,q_i)\to0$. The uniform gradient bound gives $|K_{i,2}|(q_i)\to0$. Commutation makes $K_{i,2}$ tangent to $Z_{i,1}$ and, together with orthogonality, gives $\nabla_{K_{i,1}}K_{i,2}=0$. Differentiating on $Z_{i,1}$, where $\nabla K_{i,1}$ rotates the normal plane, gives $\nabla_\nu K_{i,2}=0$ for every normal vector $\nu$. The preceding lower bound therefore gives
\begin{equation*}
 |\nabla^{Z_{i,1}}K_{i,2}|(q_i)
 =|\nabla K_{i,2}|(q_i)\geq c>0
\end{equation*}
for all sufficiently large $i$, with $c$ independent of $i$. Set $A_i:=\nabla^{Z_{i,1}}K_{i,2}(q_i)$. Since $A_i$ is skew on a two-dimensional plane, $\|A_i^{-1}\|\leq C$. In exponential coordinates on $Z_{i,1}$ and a radial parallel frame, write $K_{i,2}$ as $F_i(y)$. Jacobi comparison and the Killing identity give
\begin{equation*}
 \sup_{|y|\leq r}\|DF_i(y)-A_i\|\leq Cr,
\end{equation*}
uniformly in $i$. For fixed sufficiently small $r$, the quantitative inverse function theorem gives $F_i(y_i)=0$ with $|y_i|\leq C|F_i(0)|=C|K_{i,2}|(q_i)$. The exponential disk need only be immersed. Thus $z_i:=\exp_{q_i}(y_i)$ is a common zero, and $d(p_i,z_i)\leq d(p_i,q_i)+C|K_{i,2}|(q_i)\to0$.
\end{proof}

In the next lemma, we show that the partially unwrapped covers cannot collapse in the orbit directions, assuming the base is non-collapsed. The latter  will be verified in Lemma \ref{lem:reflected-quotient-injectivity}.

\begin{lemma}\label{lem:retained-period-volume}
Consider the covers of Definition~\ref{def:rescaled-covers}, and suppose that $d_{\tilde \g_i}(p_i, \mathcal T) \to \infty$. If the reflected quotient has injectivity radius bounded below on fixed balls and each retained primitive circle has an axis at uniformly bounded distance from the basepoint, then unit balls centered in any fixed pointed ball have volume bounded below.
\end{lemma}
\begin{proof}
Write $N_a=|K_a|$ and suppose $\operatorname{Vol}B(x_i,1)\to0$ for centers in a fixed pointed ball. The transport estimate in Lemma~\ref{lem:short-axis} uniformly bounds $N_a$ and $\dd N_a$ for each circle retained in the cover. After passage to a subsequence, each $N_a(x_i)$ tends to zero or stays bounded below. If $N_a(x_i)\to0$, Lemma~\ref{lem:short-axis} gives $d(x_i,Z_a)\to0$; if both $N_a(x_i)$ tend to zero, $x_i$ approaches a common zero of the generators. A retained circle that is not short has uniformly comparable length on a smaller fixed ball. Normalize each unwrapped direction at the center; \eqref{eq:positive-jacobi-bound} gives the same comparability there.

Near a point where neither generator vanishes, use a quotient rectangle and two orbit intervals. Near one axis, use a polar disk and two intervals; near a common zero, use two polar disks. The reflected-quotient injectivity makes the corresponding quotient normal or Fermi neighborhood embedded at a fixed radius. Jacobi comparison controls its metric, and \eqref{eq:primitive-jacobi-data} gives $N_a(r)=r+O(r^3)$ along each collapsing radius. Choose every remaining orbit interval to have a fixed length in the rescaled metric. The lower bound for $N_a$ allows this interval to fit within a retained period; a direction lifted to a line has no period restriction. Definition~\ref{def:rescaled-covers} leaves no further orbit identifications and retains the full $2\pi$ period of each collapsed circle. Integrating $N_1N_2\,\dd v_h\,\dd\vartheta_1\,\dd\vartheta_2$, the volume form in~\eqref{eq:orbit-volume}, gives a positive lower bound in every model. Choose their radii uniformly small enough that these regions lie in $B(x_i,1)$. Their positive volume lower bound contradicts $\operatorname{Vol}B(x_i,1)\to0$.
\end{proof}

\subsubsection{Conical regularity}

Let $\g_i$ have the form \eqref{eq:conical-ansatz}, with uniformly bounded Einstein constants $\varepsilon_i$ and $\beta_i\in I$. On a fixed coordinate disk $D$ centered at the conical orbit, rescale the transverse and local orbit coordinates to normalize their leading coefficients to one. Record the original primitive generators and their periods separately. The leading normal metric is $x^{2\beta_i-2}(\dd x^2+x^2\dd\vartheta^2)$. Set $W_i:=(f_i,u_i,v_i)$.

The next lemma improves weighted regularity and excludes curvature concentration at the cone when the limit is flat. At $\beta=1/2$, the first angular Fourier mode has bounded cone Hessian, so its coefficient must also be shown to tend to zero.

\begin{lemma}\label{lem:uniform-cone-estimates}
Suppose $W_i$ is uniformly bounded in $x^\delta\mathcal C_b^{2,\alpha}(D)$. On a fixed smaller disk $D'$, we have
\[
\sup_i\|W_i\|_{x^{\delta'}\mathcal C_b^{m,\alpha'}(D')}<\infty
\]
for every $m$, $0<\alpha'<1$, and $0<\delta'<\min\{1,2\beta_0\}$. If also $\beta_i\to\beta_\infty$, $\varepsilon_i\to0$, and the normalized coefficients converge locally in $\mathcal C^2$ off $x=0$ to those of a flat product cone, then $|\text{Rm}(\g_i)|_{\g_i}(q_i)\to 0$ for every sequence of regular points with $x(q_i)\to0$.
\end{lemma}
\begin{proof}
These local rescalings give $W_i\to0$ as $x\to0$. We record the original periods of the primitive circle generators separately. Setting $x=as$ maps $a<x<2a$ to $1<s<2$ and changes the principal part of \eqref{eq:scalar-b-system} to $(s\partial_s)^2+\partial_\vartheta^2$. Schauder estimates~\cite[Theorem~6.2]{GilbargTrudinger}, applied to the equations and their derivatives, give
\begin{equation*}
 \sup_i\|W_i\|_{x^\delta\mathcal C_b^{m,\alpha}(D')}<\infty
\end{equation*}
for every $m$, on a fixed smaller disk $D'$. Set $\rho_i:=e^{f_i}$. With Einstein constant $\varepsilon_i$, the first equation of \eqref{eq:scalar-b-system} becomes
\begin{equation*}
 \Delta_b\rho_i =-2\varepsilon_i x^{2\beta_i}e^{2v_i}\rho_i.
\end{equation*}
Fix $0<\delta'<\min\{1,2\beta_0\}$. Proposition~\ref{prop:fourier-estimate} gives $\rho_i-1,f_i\in x^{\delta'}\mathcal C_b^{m,\alpha}(D')$, with uniform bounds. Since $\langle\dd f_i,\dd u_i\rangle_b \in x^{\delta+\delta'}\mathcal C_b^{m,\alpha}(D')$ with uniformly bounded norm, the same proposition applied to the $u_i$ equation gives $\sup_i\|u_i\|_{x^{\delta'}\mathcal C_b^{m,\alpha}(D')}<\infty$. Hence $\tfrac14(|\dd f_i|_b^2-|\dd u_i|_b^2) \in x^{2\delta'}\mathcal C_b^{m,\alpha}(D')$ with uniform bounds. Proposition~\ref{prop:fourier-estimate} then gives $v_i\in x^{\delta'}\mathcal C_b^{m,\alpha}(D')$, also with uniform bounds. Choose $\delta'<\eta<\min\{1,2\beta_0\}$ and repeat these estimates at weight $\eta$. The inclusion $x^\eta\mathcal C_b^{m+1,\alpha}(D') \subset x^{\delta'}c_b^{m,\alpha'}(D')$ (Lemma~\ref{lem:decay-implies-little})  gives
\begin{equation*}
\sup_i\|W_i\|_{x^{\delta'}c_b^{m,\alpha'}(D')}<\infty, \qquad 0<\alpha'<1.
\end{equation*}

We now prove the curvature convergence assertion. Choose $\beta_-<\beta_\infty<a_0<\min\{1,2\beta_-\}$. Discard finitely many indices so that $\beta_-\leq\beta_i\leq(\beta_\infty+a_0)/2$. Applying the first part to this tail with lower angle bound $\beta_-$ gives $\sup_i\|W_i\|_{x^{a_0}\mathcal C_b^{m,\alpha'}(D')}<\infty$ for every $m$. Thus, for $0<x<1$,
\begin{equation*}
|\dd W_i|_{\g_{\beta_i}}\leq Cx^{a_0-\beta_i}  \leq Cx^{(a_0-\beta_\infty)/2}\longrightarrow0
\end{equation*}
uniformly in $i$ as $x\to0$. Set $r_i:=x^{\beta_i}/\beta_i$ and $H_i:=W_i+\tfrac12\varepsilon_i(r_i^2,0,0)$. Since $\Delta_b r_i^2=4x^{2\beta_i}$, \eqref{eq:scalar-b-system} gives
\begin{eqnarray*}
 \Delta_b(H_i)_f
 &=&\Delta_b f_i+2\varepsilon_i x^{2\beta_i}\\
 &=&-|\dd f_i|_{\g_b}^2
     -2\varepsilon_i x^{2\beta_i}(e^{2v_i}-1).
\end{eqnarray*}
Together with the other two equations of \eqref{eq:scalar-b-system}, this gives $\sup_i\|\Delta_bH_i\|_{x^{2a_0}\mathcal C_b^{m,\alpha'}(D')}<\infty$ for every $m$. Set $N_{i,\pm}:=e^{(f_i\pm u_i)/2}$. For a unit horizontal vector $\xi$, the plane spanned by $\xi$ and the corresponding circle direction has sectional curvature $-N_{i,\pm}^{-1}\nabla_{\hat\g_i}^2N_{i,\pm}(\xi,\xi)$. The plane tangent to the torus orbit has sectional curvature $-\langle\dd\log N_{i,+},\dd\log N_{i,-}\rangle_{\hat\g_i}$. Equation~\eqref{eqn:red-einstein-3} gives $4K_{\hat\g_i}=|\dd u_i|_{\hat\g_i}^2-|\dd f_i|_{\hat\g_i}^2$. Combining this with the preceding curvature formulas and the Hessian formula for $\hat\g_i=e^{2v_i}\g_{\beta_i}$ gives
\begin{eqnarray*}
 |\text{Rm}(\g_i)|_{\g_i}
 &\leq&C\bigl(
 |\nabla_{\g_{\beta_i}}^2 f_i|_{\g_{\beta_i}}
 +|\nabla_{\g_{\beta_i}}^2 u_i|_{\g_{\beta_i}}
 +|\dd W_i|_{\g_{\beta_i}}^2\bigr),
\end{eqnarray*}
where $C$ is independent of $i$. If $\beta_\infty<1/2$, choose $2\beta_\infty<a<\min\{1,2a_0\}$. Proposition~\ref{prop:fourier-estimate} gives $\sup_i\|H_i\|_{x^a\mathcal C_b^{m,\alpha'}(D')}<\infty$ for every $m$. Hence
\begin{equation*}
 |\nabla_{\g_{\beta_i}}^2f_i|_{\g_{\beta_i}}
 +|\nabla_{\g_{\beta_i}}^2u_i|_{\g_{\beta_i}}
 \leq C\bigl(|\varepsilon_i|+x^{a-2\beta_i}\bigr).
\end{equation*}
Since $\varepsilon_i\to0$ and $a>2\beta_\infty$, the right-hand side tends to zero at $q_i$. If $\beta_\infty=1/2$, choose $1<b<\min\{2,2a_0\}$. Proposition~\ref{prop:fourier-estimate} gives
\begin{eqnarray*}
 f_i+\tfrac12\varepsilon_i r_i^2
 &=&x\ell_{i,f}+R_{i,f},\\
 u_i&=&x\ell_{i,u}+R_{i,u},
\end{eqnarray*}
where $\ell_{i,f},\ell_{i,u}\in \text{span}\{\cos\vartheta,\sin\vartheta\}$ and $\sup_i\|(R_{i,f},R_{i,u})\|_{x^b\mathcal C_b^2(D')}<\infty$. Flatness of the limit and the mixed sectional curvature formula give $\nabla_{\hat\g_\infty}^2N_{\infty,\pm}=0$. Since $|\dd N_{\infty,\pm}|_{\hat\g_\infty}\to0$ as $x\to0$, their vanishing Hessians and the normalization give $N_{\infty,\pm}\equiv1$. Thus $f_\infty=u_\infty=0$. For each fixed $x>0$, the expansions give
\begin{equation*}
 \limsup_{i\to\infty}
 \bigl(\|\ell_{i,f}\|_{L^\infty(S^1)}
       +\|\ell_{i,u}\|_{L^\infty(S^1)}\bigr)
 \leq Cx^{b-1}.
\end{equation*}
Letting $x\downarrow0$ gives $\ell_{i,f},\ell_{i,u}\to0$. For $\ell\in\text{span}\{\cos\vartheta,\sin\vartheta\}$, differentiation in the orthonormal frame $(x^{1-\beta}\partial_x,x^{-\beta}\partial_\vartheta)$ gives
\begin{eqnarray*}
 \nabla^2_{\g_\beta}(x\ell)
 &=&(1-\beta)x^{1-2\beta}
 \begin{pmatrix}
  \ell&\partial_\vartheta\ell\\
  \partial_\vartheta\ell&-\ell
 \end{pmatrix}.
\end{eqnarray*}
For $0<x\leq1$, $\beta_i\leq1/2$ gives $x^{1-2\beta_i}\leq1$ and $x^{b-2\beta_i}\leq x^{b-1}$. The preceding estimates therefore give
\begin{eqnarray*}
 |\text{Rm}(\g_i)|_{\g_i}
 &\leq&C\bigl(
 \|\ell_{i,f}\|_{L^\infty(S^1)}
 +\|\ell_{i,u}\|_{L^\infty(S^1)}\\
 &&\qquad+x^{b-1}+|\varepsilon_i|
 +|\dd W_i|_{\g_{\beta_i}}^2\bigr).
\end{eqnarray*}
At $q_i$, every term on the right tends to zero. Hence $|\text{Rm}(\g_i)|_{\g_i}(q_i)\to0$.
\end{proof}

\subsubsection{No horizontal collapse}\label{sec:noHcollapse}

We first establish a weighted isothermal-coordinate estimate, then construct injective polar disks and normalized isothermal charts at the cone. Polar coordinates control metric distances, while isothermal coordinates put the scalar equations in the form~\eqref{eq:scalar-b-system}. We finally exclude short loops in the reflected quotient.

In a complex coordinate $z$, write a two-dimensional metric as $h=a(z)|\dd z+\mu(z)\dd\bar z|^2$, where $a>0$ and $|\mu|<1$. Recall that the function $\mu$ is its \emph{Beltrami coefficient}: it measures the failure of $z$ to be isothermal, and an isothermal coordinate $\zeta$ satisfies $\partial_{\bar z}\zeta=\mu\,\partial_z\zeta$.

\begin{lemma}\label{lem:weighted-isothermal}
Suppose the Beltrami coefficients $\mu_i$ on a disk satisfy $\mu_i=O_b^1(|z|^\sigma)$ uniformly, with $\sigma>0$. For every $0<q<\min\{1,\sigma\}$, a common smaller disk has isothermal coordinates $\zeta_i=z+O_b^{1,\alpha}(|z|^{1+q})$, for every $0<\alpha<1$.
\end{lemma}
\begin{proof}
For compactly supported $\eta$, define the Cauchy transform and the Beurling transform by
\begin{eqnarray*}
 \mathcal C\eta(z)&=&\frac1\pi\int_{\mathbb C}\frac{\eta(\xi)}{z-\xi}\,\dd A(\xi),\\
 S\eta(z)&=&\partial_z\mathcal C\eta(z)
 =-\frac1\pi\operatorname{p.v.}\int_{\mathbb C}\frac{\eta(\xi)}{(z-\xi)^2}\,\dd A(\xi).
\end{eqnarray*}
Cut off $\mu_i$ on a small disk and choose $p>2$ with $2/p<\sigma-q$. Shrink the support until $\|\mu_i\|_\infty\|S\|_{L^p\to L^p}<1/2$. The construction in~\cite[Lemmas~3--5 and Theorem~4]{AhlforsBers} gives $\eta_i=\mu_i(1+S\eta_i)$ and $\phi_i=z+\mathcal C\eta_i$, with uniform $L^p$ bounds; Lemma~8 of the same paper implies that $\phi_i$ is a homeomorphism. Since $\eta_i=\mu_i(1+S\eta_i)$ and $\sigma-2/p>q$, H\"older's inequality makes the following integral convergent:
\begin{eqnarray*}
 A_i=1-\frac1\pi\int_{\mathbb C}\frac{\eta_i(\xi)}{\xi^2}\,\dd A(\xi).
\end{eqnarray*}
Subtracting the constant and linear terms gives
\begin{eqnarray*}
 \phi_i(z)-\phi_i(0)-A_iz
 =\frac{z^2}\pi\int_{\mathbb C}
 \frac{\eta_i(\xi)}{\xi^2(z-\xi)}\,\dd A(\xi).
\end{eqnarray*}
Splitting this integral into $|\xi|\leq2|z|$ and $|\xi|>2|z|$ gives
\begin{equation}
\label{eq:isothermal-first-jet}
 \phi_i(z)-\phi_i(0)=A_iz+O(|z|^{1+q}),
 \qquad |A_i-1|<\tfrac12,
\end{equation}
after a further uniform shrinking. Writing $z=aw$ on $a<|z|<2a$, the coefficients $\mu_i(aw)$ have uniformly bounded Lipschitz norms on $1<|w|<2$. Interior Beltrami estimates~\cite[Chapter~15]{AstalaIwaniecMartin2009} therefore give $D\phi_i-A_i=O_b^{0,\alpha}(|z|^q)$. Taking $\zeta_i=(\phi_i-\phi_i(0))/A_i$ proves the assertion.
\end{proof}

The next lemma gives a uniform cone neighbourhood under a curvature bound.

\begin{lemma}\label{lem:normalized-cone-charts}
Let $\g_i\in\mathcal X_{\beta_i}$ be Einstein, with $\beta_i\in I$, and suppose $\widetilde\g_i=s_i^{-2}\g_i$ have uniformly bounded curvature. Then $(\Sigma,h_i)$ admit injective geodesic polar coordinates of uniform radius about $c_i$, together with normalized isothermal charts satisfying the uniform weighted estimates of Lemma~\ref{lem:uniform-cone-estimates}.
\end{lemma}

\begin{proof}
After a fixed rescaling take the curvature bound to be $2$, and write $h_i$ and $\varepsilon_i>0$ for the quotient metrics and Einstein constants. We first construct the geodesics issuing from $c_i$. On a sector with Cartesian coordinate $y=r_{\beta_i}e^{i\beta_i\vartheta}$, Lemma~\ref{lem:cone-regularity} gives $h_i=|\dd y|^2+O_b(|y|^2)$. Its Cartesian coefficients extend in $\mathcal C^{1,1}$ to the vertex, with Christoffel symbols $O(|y|)$ and bounded first derivatives. The locally Lipschitz geodesic equation therefore gives a unique curve in each initial direction. For $\eta>0$, solve the geodesic equation with $y(\eta)=\eta e^{i\beta_i\vartheta}$ and $y'(\eta)=e^{i\beta_i\vartheta}$. Pass to the limit as $\eta\downarrow0$ in the integral equation and its derivatives with respect to $\vartheta$. Angular differentiation is justified inductively: the $m$th Cartesian derivatives of the Christoffel symbols are $O(|y|^{1-m})$, while each angular variation is $O(t)$, so the differentiated terms remain integrable. The resulting family is smooth for $t>0$. Gauss' lemma gives the angular Jacobi factor $j_i=\beta_it+O(t^3)$ and $(j_i)_t=\beta_i+O(t^2)$.

By Lemma~\ref{lem:quotient-curvature} and \eqref{eqn:general-horizontal}, $\nabla^2f\leq-2\varepsilon_i h_i/3$. Normalize $N_a(c_i)=1$. Mixed curvature gives $|\nabla^2N_a|\leq2N_a$ and $\dd N_a(c_i)=0$, so $N_a(t)\geq\cos(\sqrt2t)$. Every axis is therefore at distance at least $\pi/(2\sqrt2)$ from the cone. The Jacobi equation with $j_i(0)=0$ and $(j_i)_t(0)=\beta_i$ excludes conjugacy before $\pi/\sqrt2$.

Fix $t_0$ below both radii. Along every geodesic leaving $c_i$, we have $f'(0)=0$ and $f''\leq-2\varepsilon_i/3$; integration forbids a return to $c_i$. Thus the radial family extends to $t_0$, since it can meet neither an axis nor the cone. If two radial geodesics first met at a regular point at distance $l<t_0$, the absence of conjugate points would make the radial exponential map a local diffeomorphism near each preimage. The images of the neighborhoods $t<l$ are disjoint, so their inward tangent half-planes are disjoint. In two dimensions this forces exactly two preimages with opposite terminal tangents. Joining the two geodesics gives a geodesic of length $2l$ from $c_i$ to itself, contradicting the preceding concavity argument. The exclusion of a return to $c_i$ uses positivity of $\varepsilon_i$, with no uniform lower bound. The polar disk of radius $t_0$ is therefore injective without a volume estimate.

To obtain isothermal coordinates, write $h_i=\dd t^2+j_i^2\dd\vartheta^2$. Radial integration of the mixed-curvature bound gives $N_a-1=O(t^2)$, with gradient $O(t)$ and Hessian $O(1)$. Equation~\eqref{eqn:red-einstein-3} then gives $K_h=O(t^2)$ and $|\dd K_h|=O(t)$. Since $|\partial_\vartheta|=j_i=O(t)$, it also gives $\partial_\vartheta K_h=O(t^2)$. Differentiating the Jacobi equation with initial data $j_i(0)=0$ and $(j_i)_t(0)=\beta_i$ yields $j_i/(\beta_it)=1+O(t^4)$; the radial and angular derivatives of this ratio are respectively $O(t^3)$ and $O(t^4)$.

Put $z=(\beta_it)^{1/\beta_i}e^{i\vartheta}$ and $x=|z|$. The preceding Jacobi bounds give $\mu_i=O_b^1(x^{4\beta_i})$. Choose $\delta<q<\min\{1,2\beta_0\}$. Lemma~\ref{lem:weighted-isothermal} gives normalized isothermal coordinates $\zeta_i$ on a common disk; write $\zeta=\zeta_i$ in the formulas below. Then $h_i=|\zeta|^{2\beta_i-2}e^{2v_i}|\dd\zeta|^2$ and $v_i=O_b^{0,\alpha}(|\zeta|^q)$. The mixed-curvature bounds give $f_i,u_i=O_b^{1,\alpha}(|\zeta|^{2\beta_i})$. Applying annular Schauder estimates~\cite[Theorem~6.2]{GilbargTrudinger} to~\eqref{eq:scalar-b-system} yields the initial $|\zeta|^q\mathcal C_b^{2,\alpha}$ bounds; Lemma~\ref{lem:uniform-cone-estimates} then supplies the higher bounds.
\end{proof}

\begin{lemma}\label{lem:reflected-quotient-injectivity}
For the covers of Definition~\ref{def:rescaled-covers}, suppose the conical orbit escapes. Then the reflected quotient has injectivity radius bounded below on every fixed pointed ball.
\end{lemma}
\begin{proof}
Use the reflected quotient defined after Definition~\ref{def:rescaled-covers}. Fixed-radius balls in this smooth surface have compact closure and remain away from the artificial boundary. A geodesic segment of length $l<\pi/\sqrt2$ cannot meet the same axis twice unless it lies in that axis. Indeed, the ambient Jacobi field $J=K_a|_\gamma$ vanishes at the two intersections, and in a parallel frame
\[
 (\pi/l)^2\int|J|^2\leq\int|D_tJ|^2
 =\int\langle\mathrm{Rm}(J,\dot\gamma)\dot\gamma,J\rangle
 \leq2\int|J|^2.
\]
Thus $J=0$. A short geodesic loop based off the axes therefore crosses each axis at most once. A crossing changes the corresponding orbit angle by $\pi$, whereas retained periods are $2\pi$; closure forces no crossings. If the basepoint lies on an axis, the same Jacobi argument confines the loop to that axis. Within the reflected section, the axis is a smooth one-dimensional manifold. A geodesic on it cannot reverse direction, so a return to its initial point gives a smoothly closed geodesic. If this geodesic met a corner, start its parametrization there and apply the same inequality to the other generator. Both generators would then vanish along the entire loop, contrary to the isolated nature of their common zeros. The loop therefore stays in an open edge interval, where it cannot close. Every short geodesic loop consequently lies in a single principal chamber, which is simply connected by the observation following Definition~\ref{def:rescaled-covers}.

If the reflected quotient injectivity radius tended to zero, the conjugate-radius bound and the local geodesic-loop alternative~\cite[p.~45]{CheegerGromovTaylor} would produce arbitrarily short loops. If such a loop has self-intersections, choose a subarc of smallest positive length whose endpoints agree. This minimum is positive because a geodesic is locally injective on its compact parameter interval. The subarc has embedded interior and is geodesic away from its new basepoint. We thus obtain a simple geodesic loop of length $L_i<1$, which bounds a disk $F$ in its simply connected chamber. Every point of $F$ can be joined to the artificial boundary by a path of length at most $4R_i$: project to $\mathcal Q_i^+$, follow a minimizing path towards the basepoint, and, if it stays in the domain, continue towards the excluded cone until the first exit. Lift this path from the chosen point. It exits the compact disk $F$ before reaching the artificial boundary, so the inradius of $F$ is at most $D_i=4R_i$.

To estimate the area without assuming injectivity, cover $F$, up to its cut locus, by inward normals from its geodesic boundary and the inward sector at its vertex. Comparison with curvature $-2$ bounds their Jacobians, and the sector angle is at most $2\pi$. Hence
\begin{equation}
\label{eq:loop-disk-area}
 \operatorname{Area}(F)\leq
 \frac{L_i}{\sqrt2}\sinh(\sqrt2D_i)
 +\pi\bigl(\cosh(\sqrt2D_i)-1\bigr)
 \leq C V_{2,-2}(4R_i+1).
\end{equation}
Gauss--Bonnet gives $\int_FK_{h_i} = \pi + \alpha_i \geq\pi$ for a simple geodesic loop, where $\alpha_i \geq 0$ is the interior angle at the loop basepoint. Lemma~\ref{lem:quotient-curvature}, \eqref{eq:loop-disk-area}, and Definition~\ref{def:rescaled-covers} give the contradiction:
\begin{eqnarray*}
\pi\leq\int_FK_{h_i}
\leq\frac{\varepsilon_i}{3}\text{Area}(F)
\leq C\varepsilon_iV_{2,-2}(4R_i+1)\longrightarrow0.
\end{eqnarray*}
\end{proof}

\subsubsection{Compactness of the covers}

All bounds in this part are on fixed pointed balls: if $A>0$, then $R_i>A+1$ for large $i$, so the ball of radius $A$ about $\widehat p_i$ and its unit neighbourhood lie away from the artificial boundary.

\begin{lemma}[Unwrapped covers are non-collapsed]\label{lem:covered-local-bounds}
For every $A>0$, there exists $v_A>0$ such that, for all sufficiently large $i$,
\[
\operatorname{Vol}  B_{\tilde\g_i}(x,1)\geq v_A,
\qquad\text{for every }x\in B_{\tilde\g_i}(\widehat p_i,A)\subset\widehat U_i.
\]
On $B_{\widetilde\g_i}(\widehat p_i,A)$, the lifted metrics admit uniform smooth local coefficient bounds at any fixed positive distance from the conical orbit. When that orbit is present, they also satisfy the weighted bounds of Lemma~\ref{lem:uniform-cone-estimates} in normalized cone charts.
\end{lemma}

\begin{proof}
Suppose first that the cone remains at bounded distance. Lemma~\ref{lem:normalized-cone-charts} supplies its normalized chart. A retained generator has an axis at bounded distance, where $|\nabla K_a|=\sqrt2$. The transport estimate of Lemma~\ref{lem:short-axis} gives $|K_a|(c_i)\geq c>0$, since the first derivatives of the orbit metric coefficients vanish at $c_i$. Normalize only unwrapped directions. The polar disk times a fixed orbit rectangle is then uniformly comparable to a product cone with $\beta\geq\beta_0$. This gives a volume lower bound on a fixed cone neighbourhood. Local volume comparison propagates it along paths of bounded length; portions meeting the cone are covered by the same uniform product-cone charts. Thus, for each $A$,
\begin{eqnarray}
\label{eq:cone-local-volume}
\text{Vol}_{\widetilde\g_i}B(x,a)&\geq&c_Aa^4,
\qquad x\in B(\widehat p_i,A),\quad 0<a<a_A.
\end{eqnarray}

If the cone escapes, Lemma~\ref{lem:reflected-quotient-injectivity} bounds the reflected-quotient injectivity radius below. Definition~\ref{def:rescaled-covers} retains precisely the fixed circles whose axes stay at bounded distance. Lemma~\ref{lem:retained-period-volume} therefore gives noncollapse of the four-dimensional lift.

On each fixed smooth ball, take the comparison radius below its distance to the cone and to the artificial boundary. Noncollapse now permits the local injectivity estimate \cite[Theorem~4.3]{CheegerGromovTaylor}; harmonic coordinates and Einstein regularity then give smooth local bounds \cite[Chapter~11]{petersen2016riemannian}, \cite[Theorem~5.2]{DeTurckKazdan1981}. For retained primitive generators, the transport estimate of Lemma~\ref{lem:short-axis} bounds their first jets. Normalize each unwrapped generator at $\widehat p_i$ and use~\eqref{eq:positive-jacobi-bound} for the same bound. Differentiating the Killing identity then bounds all their derivatives in these charts. The ambient estimates apply also across axes and corners, which are smooth points of the lifted metrics. Lemma~\ref{lem:slice-regularity} and the expansions~\eqref{eqn:edge-slice-expansion}--\eqref{eqn:vertex-slice-expansion} then give uniform bounds for the coefficients in quotient and orbit coordinates. Together with the conical charts above, these are the required bounds.
\end{proof}

\begin{lemma}\label{lem:marked-compactness}
After passing to a subsequence, $(\widehat U_i,\widetilde\g_i,\widehat p_i)$ converges to a complete pointed limit $(X,\g_\infty,p_\infty)$ whose regular part is Ricci flat. The convergence is smooth away from the conical orbit, if present, and holds in the normalized weighted coefficients near that orbit.
\end{lemma}
\begin{proof}
Lemma~\ref{lem:covered-local-bounds} gives noncollapse and uniform local coordinates. Smooth pointed compactness applies away from the cone \cite[Chapter~11]{petersen2016riemannian}. At the cone, Lemma~\ref{lem:uniform-cone-estimates} improves the weight and differentiability; Lemma~\ref{lem:compactness} then gives convergence in the prescribed little H\"older norm. A diagonal subsequence joins these local limits. The artificial boundaries escape to infinity, so the completed limit is complete. It is Ricci flat because $\varepsilon_i\to0$. The limiting orbit fields are Killing and orthogonal. The fixed orbit coordinates retain the circle periods or line parameters chosen in Definition~\ref{def:rescaled-covers}.
\end{proof}

\begin{proposition}[Cone point escapes]\label{prop:noedgebubble}
For a sequence as in Notation~\ref{not:rescaled-covers}, we have $d_{\tilde\g_i}(p_i,\mathcal T)\to\infty$.
\end{proposition}
\begin{proof}
Suppose the conical orbit remains at bounded distance. On $\widehat U_i$, set $\psi_i:=\widetilde\rho_i/\widetilde\rho_i(c_i)$, viewing the quotient functions as pulled back to the cover. Lemma~\ref{lem:density-maximum} and~\eqref{eqn:general-trace} give
\[
0\leq\psi_i\leq1,\qquad \psi_i(c_i)=1,\qquad
\Delta_{h_i}\psi_i+2\varepsilon_i\psi_i=0.
\]
By Lemma~\ref{lem:marked-compactness}, these functions converge to a harmonic function $\psi$ on the limiting quotient, with $\psi(c_\infty)=1$. In an isothermal coordinate, $\psi$ is bounded and Euclidean harmonic off the puncture, so it extends harmonically across it. The strong maximum principle gives $\psi\equiv1$.

Thus $df=0$, and~\eqref{eqn:general-horizontal} reduces to
\[
K_hh=\tfrac12\,du\otimes du.
\]
Comparing ranks yields $K_h=0$ and $du=0$, so the limit is a flat product cone. This contradicts~\eqref{eq:curvature-normalization}, by smooth convergence if the $\hat p_i$'s stay at a uniform positive distance from the cone, and by Lemma~\ref{lem:uniform-cone-estimates} if they approach it.
\end{proof}

\subsubsection{Non-collapsing}

Noncollapse of the covers leaves one possible source of collapse downstairs: deck transformations with displacements tending to zero. The next lemma converts such transformations into a nowhere-zero horizontal Killing field on the limit.

\begin{lemma}\label{lem:free-section-displacements}
Suppose $(\widehat U_i,\widetilde\g_i,\widehat p_i)$ converges smoothly to a connected complete smooth limit $(X,\g_\infty,p)$. Let $\eta_i$ be fixed-point-free isometries of the covering interiors, commuting with the orbit flows and preserving each reflected section and its orientation. If $d_{\widetilde\g_i}(\widehat p_i,\eta_i(\widehat p_i))\to0$, then $X$ admits a nowhere-zero horizontal Killing field whose orthogonal distribution is integrable.
\end{lemma}
\begin{proof}
Move $\widehat p_i$ onto one section by an orbit translation. Since $\eta_i(B(\widehat p_i,R))\subset B(\widehat p_i,R+d(\widehat p_i,\eta_i\widehat p_i))$, the image of each fixed-radius ball lies in a slightly larger fixed-radius ball for all sufficiently large $i$. Smooth compactness of isometries follows from
\begin{equation}
\label{eq:isometry-first-jet}
 \eta_i(\exp_xv)=\exp_{\eta_ix}(\dd\eta_i|_xv).
\end{equation}
Apply it also to the inverses to obtain an onto limiting isometry $\eta$ fixing $p$. The sections converge as exponential images of their tangent planes. On the limiting section, $\dd\eta|_p$ is an orientation-preserving linear isometry. A non-identity rotation would give a nondegenerate fixed point, which persists for $\eta_i$, contrary to freeness. Thus \eqref{eq:isometry-first-jet} makes $\eta$ the identity near $p$ on the section. Commutation with the orbit flows makes it the identity on a principal open set, hence everywhere.

Set $V_i(x)=\exp_x^{-1}(\eta_ix)$ and $a_i=|j^1_{\widehat p_i}V_i|>0$, where the first-jet norm is $(|V_i|^2+|\nabla V_i|^2)^{1/2}$. Then $a_i\to0$. On a normal ball, \eqref{eq:isometry-first-jet} determines $\eta_i$ smoothly from its value and differential at the center. Their deviation from the identity is uniformly comparable to the value and derivative of $V_i$ there. Applying the mean-value theorem to this dependence, then covering $B(\widehat p_i,R)$ by finitely many overlapping normal balls, gives
\begin{equation}
\label{eq:displacement-jet-bound}
 \|V_i\|_{C^j(B(\widehat p_i,R))}\leq C_{R,j}a_i.
\end{equation}
Thus $a_i^{-1}V_i\to Z$ smoothly, with $|j^1_pZ|=1$. Dividing $\eta_i^*\widetilde\g_i-\widetilde\g_i=0$ by $a_i$ proves that $Z$ is Killing. Naturality of the exponential map makes $V_i$ orbit-invariant. The two orbit-coordinate reversals fix $x$ and $\eta_ix$, hence also their unique short geodesic;  $V_i$ is tangent to every reflected section, including those through axes and corners. Hence $Z$ is horizontal and orbit-invariant.

If $Z(q)=0$, its restriction to the reflected section through $q$ is a nontrivial Killing field: otherwise orbit invariance would make $Z$ vanish on a principal open set. Its derivative at $q$ is a nonzero skew map in dimension two, so its zero is nondegenerate. Smooth convergence in~\eqref{eq:displacement-jet-bound} makes that zero persist for $a_i^{-1}V_i$, giving a fixed point of $\eta_i$. Thus $Z$ never vanishes. On the principal set $Z^\flat$ is pulled back from the two-dimensional quotient, so $Z^\flat\wedge\dd Z^\flat=0$; continuity extends this identity everywhere.
\end{proof}

\begin{proposition}\label{prop:no-collapsed-smooth-limit}
For a sequence as in Notation~\ref{not:rescaled-covers}, we have $\inf_i\operatorname{inj}_{\tilde\g_i}(p_i)>0$.
\end{proposition}
\begin{proof}
By Proposition~\ref{prop:noedgebubble} and Lemma~\ref{lem:marked-compactness}, a subsequence of $(\widehat U_i,\widetilde\g_i,\widehat p_i)$ converges to a complete smooth Ricci-flat limit $(X,\g_\infty,p)$ with $|\mathrm{Rm}(\g_\infty)|(p)=1$. If a direction were unwrapped, its axis would escape. Normalize the generator of this line action to have unit length at $p_i$ and apply the estimates in the proof of Lemma~\ref{lem:properness-marked-action}. Those estimates do not require periodicity and give a global nowhere-zero Killing field with integrable orthogonal distribution. Lemma~\ref{lem:properness-static-field} would force flatness. Thus neither orbit direction is unwrapped: both primitive generators retain their $2\pi$ periods on $\widehat U_i$.

Suppose, for a contradiction, that $\operatorname{inj}_{\widetilde\g_i}(p_i)\to0$ along a subsequence. The curvature bound excludes nearby conjugate points, so the local geodesic-loop alternative~\cite[p.~45]{CheegerGromovTaylor} gives loops based at $p_i$ with lengths tending to zero. For large $i$, these loops lie in $U_i$, and their lifts to $\widehat U_i$ have distinct endpoints because the lifted metrics have uniformly positive injectivity radius near $\widehat p_i$. Thus their endpoints are $\widehat p_i$ and $\eta_i(\widehat p_i)$ for nontrivial, fixed-point-free deck transformations satisfying
\[
d_{\widetilde\g_i}\bigl(\widehat p_i,\eta_i(\widehat p_i)\bigr)\to0.
\]
Since both orbit periods are retained, the deck transformations fix the orbit coordinates, commute with the orbit flows, and preserve each reflected section and its orientation. Lemma~\ref{lem:free-section-displacements} therefore gives a nowhere-zero horizontal Killing field on $X$ with integrable orthogonal distribution. Finally, Lemma~\ref{lem:properness-static-field} forces $\g_\infty$ to be flat, contradicting $|\mathrm{Rm}(\g_\infty)|(p)=1$.
\end{proof}

\subsection{Curvature bounds and convergence}

We now exclude curvature blow-up using the compact-edge conclusion and the density estimate. The resulting curvature and volume bounds give geometric compactness; the final argument recovers the prescribed conformal coordinates and weighted topology.

The following lemma gives the final flatness contradiction in Proposition~\ref{prop:curvature-noncollapse}.

\begin{lemma}\label{prop:bounded-density}
Let $(X,\g)$ be a connected complete four-manifold with Killing fields $K_1,K_2$ satisfying
\[
\operatorname{Ric}(\g)=0,\qquad [K_1,K_2]=0,\qquad \g(K_1,K_2)=0.
\]
Set $\rho:=|K_1\wedge K_2|_\g$. Suppose $\{\rho>0\}$ is dense and the distribution $\operatorname{span}\{K_1,K_2\}^{\perp}$ is integrable there. Let $h$ denote the local quotient metric on this set, and assume
\[
\sup_X\rho<\infty,\qquad K_h\leq0.
\]
Then $\g$ is flat.
\end{lemma}
\begin{proof}
On $\{\rho>0\}$, where $\dd\rho\neq0$, follow the unit gradient $e=\nabla\rho/|\dd\rho|$ and let $e^\perp$ be a horizontal unit normal to it. Equations~\eqref{eqn:general-horizontal}--\eqref{eqn:general-trace} and Remark~\ref{rem:diagonal-reduction} give
\begin{eqnarray*}
\nabla_\g^2\rho(e,e)
=\rho\left(-K_h+\tfrac12\dd u(e^\perp)^2\right)\geq0.
\end{eqnarray*}
Thus along its integral curve, $\rho'=|\dd\rho|>0$ and $\rho''\geq0$. The functions $\rho$ and $|\dd\rho|$ stay bounded below by positive constants. If the curve had a finite maximal forward time, its unit speed and completeness would give a limit point in the independent, noncritical locus, where the unit gradient field is smooth; the curve would therefore extend. Thus it exists for all positive time, and $\rho(\gamma(t))\geq\rho(\gamma(0))+t|\dd\rho(\gamma(0))|$, contrary to boundedness. Hence $\dd\rho=0$ on $\{\rho>0\}$. Equation~\eqref{eqn:general-horizontal} reduces to $K_hh=\tfrac12\dd u\otimes\dd u$ and forces $K_h=0$ and $\dd u=0$. By Remark~\ref{rem:diagonal-reduction}, the metric is a flat product there, and density of this set gives flatness everywhere.
\end{proof}

\begin{proposition}\label{prop:curvature-noncollapse}
There exists a constant $C>0$ such that every Einstein metric $\g\in\mathcal X_\beta$, with $\beta\in I$, has $|{\rm Rm}(\g)|_\g\leq C$ on $S^4\setminus\mathcal T$ and ${\rm Vol}(\g)\geq 64 \pi^2 / C^2$.
\end{proposition}
\begin{proof}
Suppose the curvature bound fails. Choose $q_i=|\text{Rm}(\g_i)|(p_i)\geq\tfrac12\sup|\text{Rm}(\g_i)|\to\infty$, and rescale by $\widetilde\g_i=q_i\g_i$ as per Notation \ref{not:rescaled-covers}. Propositions~\ref{prop:noedgebubble} and~\ref{prop:no-collapsed-smooth-limit}, followed by the smooth compactness argument of Lemma~\ref{lem:marked-compactness}, give a complete smooth Ricci-flat limit with $|\text{Rm}|(p_\infty)=1$. Lemma~\ref{lem:quotient-curvature} gives $K_{h_\infty}\leq0$. Proposition~\ref{prop:polarbubbles} gives a compact edge $E_\infty$, identified with an entire original edge $E_i$. Primitive periods and both endpoints are preserved. Set $P_i:=P(\g_i)$, as in Lemma~\ref{lem:density}. Let $b_i$ be the length, in the rescaled metric, of the generator that does not vanish on $E_i$. The expansion~\eqref{eqn:edge-slice-expansion} gives $\partial_\nu\rho=b$. Combining this with~\eqref{eq:density-boundary} and smooth equivariant convergence along the entire edge, with the primitive generators unchanged, gives
\begin{eqnarray*}
\frac{q_iP_i}{2}
=-\int_{E_i}*_{h_i}\dd\widetilde\rho_i
=\int_{E_i}b_i\,\dd s_i
\longrightarrow\int_{E_\infty}b_\infty\,\dd s_\infty\in(0,\infty).
\end{eqnarray*}
Thus $q_iP_i$ is bounded. Lemma~\ref{lem:density} bounds $\widetilde\rho_i=q_i\rho_i$, and hence the limiting density. By Lemma~\ref{prop:bounded-density}, the limit is flat, a contradiction. With $|\text{Rm}|\leq C$, \eqref{eqn:L2-bound} gives $64\pi^2\leq C^2 {\rm Vol}(\g)$.
\end{proof}

\begin{lemma}\label{lem:cone-separation-trace}
There exist constants $d_0>0$ and $0<a_0\leq a_1<\infty$ such that every Einstein metric $\g\in\mathcal X_\beta$, with $\beta\in I$, satisfies
\[
\operatorname{dist}_{\hat\g}(c,\partial\Sigma)\geq d_0,
\qquad a_0\leq\rho(c)\leq a_1.
\]
\end{lemma}
\begin{proof}
Take $C$ from Proposition~\ref{prop:curvature-noncollapse}, $D=\pi\sqrt{3/\lambda}$ from Lemma~\ref{lem:global-bounds} and set $V_{\text{min}} = 64 \pi^2 / C^2$. Set $l:=\operatorname{dist}_{\hat\g}(c,\partial\Sigma)$. Along a shortest quotient geodesic from $c$ to the boundary, choose a generator vanishing at the endpoint. Its length $\phi$ satisfies $\phi''+C\phi\geq0$, $\phi(0)>0$, $\phi'(0)=0$, and $\phi(l)=0$. For $t<\min\{l,\pi/(2\sqrt C)\}$, comparison gives $\phi(t)\geq\phi(0)\cos(\sqrt Ct)>0$. Hence its first zero cannot occur before $d_0=\pi/(2\sqrt C)$, so $l\geq d_0$. Minimizing geodesics from $c$ to interior points avoid the boundary as in Lemma~\ref{lem:density-maximum}. Since $K_{\hat g}\geq-C$, polar comparison, integrated up to the cut time, gives
\begin{eqnarray*}
\operatorname{Area}_{\hat\g}(\Sigma)
\leq2\pi\beta\int_0^D\frac{\sinh(\sqrt Ct)}{\sqrt C}\dd t
\leq\frac{2\pi}{C}(\cosh(\sqrt CD)-1)=:A_0.
\end{eqnarray*}
By Lemma~\ref{lem:density-maximum} and~\eqref{eq:orbit-volume}, $V_{\text{min}}\leq4\pi^2\rho(c)A_0$. This gives $a_0=V_{\text{min}}/(4\pi^2A_0)$. The upper bound follows from Lemmas~\ref{lem:density} and~\ref{lem:global-bounds}.
\end{proof}

The final step recovers the prescribed coordinates from the geometric limit and proves convergence in the topology of $\mathcal X_I$ defined in \S\ref{sec:space_conical}.

\begin{theorem}\label{thm:compact-zero-set}
For every compact interval $I\subset(0,1/2]$, the total zero set
\[
\mathscr Z_I:=\{(\beta,\g)\in\mathcal X_I: {\rm Ric}(\g) = \lambda \g\}
\]
is compact.
\end{theorem}

\begin{proof}
We first obtain geometric convergence, then recover the prescribed conformal coordinates. Take a sequence $(\beta_i,\g_i)\in\mathscr Z_I$ and pass to a subsequence with $\beta_i\to\beta_\infty\in I$. Proposition~\ref{prop:curvature-noncollapse} and Lemma~\ref{lem:global-bounds} bound curvature and diameter above and volume below. Bishop comparison on the complete $\mathrm{RCD}(\lambda,4)$ spaces gives a uniform lower bound for the volume of every fixed-radius ball. At any fixed positive distance from $\mathcal T$, the local Cheeger--Gromov--Taylor estimate~\cite[Theorem~4.3]{CheegerGromovTaylor} therefore bounds the injectivity radius below. In uniformly controlled harmonic charts, elliptic estimates for the Einstein equation~\cite[\S5]{DeTurckKazdan1981} give bounds for all derivatives on smaller charts, and hence smooth subsequential convergence away from the cone~\cite[Chapter~11]{petersen2016riemannian}.

Each of the two fixed circle generators $K_a$ vanishes on its axis with $|\nabla K_a|=\sqrt2$. The estimate $|\nabla^2K_a|\leq C|\mathrm{Rm}|\,|K_a|$ and the diameter bound give uniform bounds for $K_a$ and $\nabla K_a$ by transport from its axis. Apply the averaging construction in the proof of Lemma~\ref{lem:properness-marked-action} also to $\Gamma$ and the coordinate conjugations $\kappa_j$. The resulting convergence is equivariant, preserves the primitive generators and labelled axes, and identifies the meridional sections. Lemma~\ref{lem:cone-separation-trace} separates the cone from the boundary and bounds its orbit-area density above and below, while Lemma~\ref{lem:normalized-cone-charts} supplies normalized isothermal cone charts on a common disk.

We recover the conformal marking on the global section $S=S_\Sigma^2$, using the reflections $\eta_1,\eta_2$ and the four-point set $C$ of \S\ref{sec:scalar_eqn_function_spaces}. Also write $\hat\g_i$ for the reflected quotient metric on $S$. By Lemma~\ref{lem:slice-regularity}, it is smooth across the axes and at the vertices, where four right-angle quadrants meet. The six points $p_\pm$ and $C$ remain distinct. Indeed, a uniform equivariant normal ball about a torus fixed point contains no other fixed point, since its isotropy representation has no fixed vector; the injectivity estimates therefore separate the vertices. The cone-separation bound separates the four lifts of $c$ from the axes and hence from each other.

Since the section is totally geodesic, restrictions of ambient exponential charts give uniform normal charts on $S$ away from the four cone points, including at the axes and vertices. Smooth compactness on the complements of four smaller cone disks, followed by attaching the limiting disks, gives a limiting section diffeomorphic to $S$. Choose a finite conformal atlas consisting of ordinary isothermal disks and normalized cone disks at the four points of $C$. In a cone disk,
\[
\hat\g_i=e^{2v_i}|\zeta|^{2\beta_i-2}|\dd\zeta|^2;
\]
the scalar factor does not affect the conformal structure, which extends across the center in the coordinate $\zeta$. Take disjoint disks about the six marked points and arrange that all overlaps avoid their centers. Smooth convergence on the regular part and the normalized cone estimates give corresponding charts on common smaller domains, with transition maps converging smoothly on the overlaps.

Use the geometric identifications off the marked disks and the exact coordinate identifications on smaller disks about their centers. On the intervening annuli, the maps are $\mathcal C^1$-close and can be joined using cutoff functions in common target coordinate charts. For large $i$ the resulting maps are local diffeomorphisms, hence diffeomorphisms because the surfaces are compact and the target is a sphere. They identify the sections with one fixed sphere and preserve the marked points. In these identifications, the Beltrami coefficients relative to the limiting conformal structure converge uniformly to zero: they vanish on the smaller marked disks, and smooth convergence controls the complement.

Let $c_*\in C$ be the lift in the chamber $r_1,r_2\geq0$. As in Lemma~\ref{lem:fixed-marking}, uniformize the sections by maps $W_i$ satisfying
\[
W_i(p_-)=0,\qquad W_i(p_+)=\infty,\qquad W_i(c_*)=1.
\]
The normalized maps and their inverses converge uniformly to a biholomorphism $W_\infty$ and its inverse~\cite[Theorems~8--9 and Lemma~18]{AhlforsBers}. In the corresponding local holomorphic charts, Cauchy's estimates give convergence of all derivatives on smaller disks. In particular, the limiting chart transitions have nonzero first derivatives.

To recover the marking of Lemma~\ref{lem:fixed-marking}, set $\iota:=\eta_1\eta_2$. This holomorphic involution fixes $p_\pm$ and is nontrivial, so $W_i\circ\iota=-W_i$. Consequently,
\[
Z_i=iW_i^2
\]
descends to $S/\langle\iota\rangle=\widehat\Sigma$ and satisfies $Z_i(p_-)=0$, $Z_i(p_+)=\infty$, and $Z_i(c_+)=i$. Lemma~\ref{lem:fixed-marking} gives $Z_i(c_-)=-i$ and identifies $Z_i$ with the prescribed normalization. Since each $\hat\g_i$ is conformal to the fixed round structure, this recovers exactly the original conformal coordinates, with their labelled reflection axes and $\Gamma$ symmetries. Keeping the orbit angles fixed extends these identifications smoothly across axes and corners, as in Lemma~\ref{lem:slice-regularity} and the proof of Lemma~\ref{lem:ambient-einstein-regularity}.

For the holomorphic transition $z=\phi_i(\zeta)$ from a normalized cone chart to the fixed local chart, we have $\phi_i(0)=0$ and $\phi_i'(0)\to\phi_\infty'(0)\ne0$. Comparing leading normal coefficients gives
\[
1=a_{0,i}|\phi_i'(0)|^{2\beta_i},\qquad
a_{0,i}\longrightarrow|\phi_\infty'(0)|^{-2\beta_\infty}>0.
\]
The common leading orbit coefficient is $b_{0,i}=\rho_i(c)$, so Lemma~\ref{lem:cone-separation-trace} gives a subsequence with a positive limit for $b_{0,i}$. All these coefficients are measured against the original primitive generators.

Choose $\delta<\delta'<\min\{1,2\beta_0\}$, $m>k+2$, and $\alpha<\alpha'<1$. Lemma~\ref{lem:uniform-cone-estimates} bounds the normalized remainders in $x^{\delta'}\mathcal C_b^{m,\alpha'}$. The holomorphic transitions and their inverses converge smoothly on common smaller disks with first derivatives bounded away from zero. Thus the same bounds hold in the fixed chart: the two radial variables are uniformly comparable, and the smooth correction to the leading normal factor is $O(x)$, which is allowed since $\delta'<1$. Since~\eqref{eq:angle-identification} leaves the scalar coefficients unchanged, Lemma~\ref{lem:compactness}, together with smooth convergence on the compact complement, gives
\begin{equation}
\label{eq:properness-topology}
\mathcal I_{\beta_i}\g_i\longrightarrow\mathcal I_{\beta_\infty}\g_\infty
\quad\text{in }\mathcal S_{1/2}^{k+2,\alpha}.
\end{equation}
The geometric limit is positive away from $\mathcal T$, and the limiting leading coefficients are positive. The defining tensor form and symmetries are closed conditions, so $\g_\infty\in\mathcal X_{\beta_\infty}$. The Einstein equation passes to the limit under local $\mathcal C^2$ convergence away from $\mathcal T$, giving $\operatorname{Ric}(\g_\infty)=\lambda\g_\infty$. Hence $\mathscr Z_I$ is compact.
\end{proof}

\section{Proof of Theorem \ref{thm:intro-existence}}
\label{sec:half-angle}

After having established the Fredholm theory and compactness, the only  remaining ingredient for the mod-two degree argument is to show that the degree at $\beta = 1/2$ is non-zero. 

Recall from $\S$\ref{sec:scalar_eqn_function_spaces} that $\mathcal M$ is the fixed Banach space of scalar triples $(p,q,v)$ satisfying the prescribed symmetries and the vertex normalization \eqref{eq:scalar-normalization}. For each $\beta\in I=[\beta_0,1/2]$, the corresponding metric on principal orbits is
\[
 \g=e^{2v}h_\beta^\circ
      +r_1^2e^{p+q}\dd\vartheta_1^2
      +r_2^2e^{p-q}\dd\vartheta_2^2,
\]
as in \eqref{eq:scalar-enlarged-metric}. The (weighted) Einstein operator $\mathcal R_\beta$ is defined in \eqref{eq:scalar-residual-map}. By Theorem~\ref{thm:scalar-index-zero}, the family $\mathcal R_\beta$ defines a smooth section of the Banach bundle over $I\times\mathcal M$ with fibres $\mathcal Z_{\beta,\g}$, whose derivative in the $\mathcal M$-direction is Fredholm of index zero at every solution $\mathcal R_\beta(\g) = 0$. Proposition~\ref{prop:scalar-zero-equivalence} further identifies these zeros with  conical metrics $\g\in\mathcal X_\beta$ satisfying $\operatorname{Ric}(\g)=\lambda\g$.

Theorem~\ref{thm:compact-zero-set} makes the total zero set compact, so Proposition~\ref{prop:localized-degree} and Remark~\ref{rem:scalar-degree-projections} give a well-defined degree $\deg_2(\mathcal R_\beta)$, independent of $\beta\in I$. It therefore remains to compute $\deg_2(\mathcal R_{1/2})$. We show that $\mathcal R_{1/2}$ has a unique zero, represented by a quotient of the equal round product on $S^2\times S^2$
\begin{eqnarray}\label{eq:round-product}
 \g_P&=&\lambda^{-1}\sum_{j=1}^2
 \bigl(\dd s_j^2+{\cos^2}s_j\,\dd\vartheta_j^2\bigr)
\end{eqnarray}
under the branched covering map $S^2 \times S^2 \to S^4$, and that this zero is nondegenerate.

\begin{lemma}\label{lem:half-angle-zero}
The section $\mathcal R_{1/2}$ has exactly one zero in $\mathcal M$. Its corresponding metric $\g_{1/2}\in\mathcal X_{1/2}$ is the quotient of the equal round product on $S^2\times S^2$ by simultaneous reflection across the two equators.
\end{lemma}
\begin{proof}
By Proposition~\ref{prop:scalar-zero-equivalence}, it suffices to prove existence and uniqueness of the Einstein metric in $\mathcal X_{1/2}$. Let $a_j$ reflect the $j$th factor across its equator, and let $S$ exchange the factors. On the quotient by $a_1a_2$, use the coordinates $\zeta=(s_1+is_2)^2$, $\vartheta_1$, and $\vartheta_2$ near the image of the fixed torus. The orbit space of this smooth closed four-manifold is the latitude square modulo $(s_1,s_2)\mapsto(-s_1,-s_2)$: a disk with two corners, primitive edge labels $e_1,e_2$, and no exceptional orbits. The equivariant reconstruction theorem~\cite[\S4.2, p.~546]{OrlikRaymond1} identifies it with the fixed $S^4$ action, sending the central orbit to $\mathcal T$.

Near that orbit,
\[
 \hat\g_P=(4\lambda)^{-1}|\zeta|^{-1}|\dd\zeta|^2,
 \qquad e^f=\lambda^{-1}\cos s_1\cos s_2,
 \qquad u=\log(\cos s_1/\cos s_2).
\]
The normalized conformal identification of Lemma~\ref{lem:fixed-marking} intertwines the descended $S,a_1$ with $\sigma,\tau$. Keeping the orbit angles unchanged extends it smoothly across axes and corners by reflection. At the cone, $\zeta=Az+O(z^2)$ with $A\ne0$, so
\[
 \hat\g_P=\frac{|A|}{4\lambda}|z|^{-1}(1+O_b(|z|))|\dd z|^2,
 \qquad f+\log\lambda,u=O_b(|z|).
\]
Thus the cone angle is $\pi$ and $v-v_0=O_b(x)$, where $v_0=\tfrac12\log(|A|/(4\lambda))$. Since $\delta<1$, Lemma~\ref{lem:decay-implies-little} places these remainders in $x^\delta c_b^{k+2,\alpha}$, so the quotient lies in $\mathcal X_{1/2}$.

Conversely, Proposition~\ref{prop:removability} lifts any Einstein metric in $\mathcal X_{1/2}$ to $S^2\times S^2$, where Theorem~\ref{thm:product-page-uniqueness} identifies it with the equal round product. Conjugate its torus action to the product rotation action. The induced torus automorphism preserves the two primitive stabilizer subgroups, so it is a signed permutation. Product isometries absorb these changes of basis. The deck involution commutes with both rotation circles and fixes the ramification torus; it therefore reflects both factors across their equators and equals $a_1a_2$.

The descended isometry has the form $(q,\vartheta)\mapsto(\phi(q),\vartheta+\psi(q))$. The absence of mixed quotient--orbit terms gives $\dd\psi=0$; a constant torus translation makes $\psi=0$. Compose with $\tau$, if necessary, to fix both corners. Then $\phi$ preserves orientation and each labelled edge. Its holomorphic double fixes the three normalized points, so Lemma~\ref{lem:fixed-marking} gives $\phi=\mathrm{id}$.
\end{proof}

After having established uniqueness, we now prove non-degeneracy, thus allowing us to compute the degree:

\begin{proposition}\label{prop:half-angle-degree}
The derivative
\[
 \mathrm D_{\g_{1/2}}\mathcal R_{1/2}:
 T_{\g_{1/2}}\mathcal M\longrightarrow\mathcal Z_{1/2,\g_{1/2}}
\]
is an isomorphism, and $\deg_2(\mathcal R_{1/2})=1$.
\end{proposition}
\begin{proof}
Let $(\phi,\psi,\chi)$ belong to the kernel. Proposition~\ref{prop:scalar-zero-equivalence} identifies its metric variation with an infinitesimal Einstein deformation, with ambient $c^{k+2,\alpha}$ regularity at the axes and vertices. Denote its pullback to the round product by $H$.

At the branch locus, use $z=w^2$ and put $R=|w|$. Write $f,u$ for the smooth lifted background functions and $\hat\g_P=e^{2V}|\dd w|^2$. The pulled-back scalar variations satisfy
\[
 \begin{split}
 \Delta_w\phi+2\langle\dd f,\dd\phi\rangle
                         +4\lambda e^{2V}\chi&=0,\\
 \Delta_w\psi+\langle\dd f,\dd\psi\rangle
                         +\langle\dd u,\dd\phi\rangle&=0,\\
 \Delta_w\chi-\tfrac12\langle\dd f,\dd\phi\rangle
                         +\tfrac12\langle\dd u,\dd\psi\rangle&=0,
 \end{split}
\]
where the inner products are Euclidean. Their defining asymptotics give
\[
 (\phi,\psi,\chi)=(\phi_0,0,\chi_0)+O(R^{2\delta}),
 \qquad |D_w(\phi,\psi,\chi)|=O(R^{2\delta-1}).
\]
Since $\delta>0$, these functions lie locally in $W^{1,2}$ and their boundary flux on $R=\eta$ is $O(\eta^{2\delta})$. Integration by parts therefore extends the equations weakly across $w=0$. Elliptic regularity makes the variations smooth there. Hence $H$ is a $c^{k+2,\alpha}$ infinitesimal Einstein deformation on the whole product.

\phantomsection\label{notloc:gauge-conventions}
Use $(\delta H)_j=-\nabla^iH_{ij}$ and $\delta^*X=\tfrac12\mathcal L_X\g_P$. The elliptic self-adjoint operator $\delta\delta^*$ has the Killing fields as its kernel. Since $\langle\delta H,K\rangle_{L^2}=\langle H,\delta^*K\rangle_{L^2}=0$ for every Killing field $K$, solve
\[
 \delta\delta^*X=\delta H,
 \qquad X\in c^{k+3,\alpha}.
\]
Average over the torus, the two circle-coordinate conjugations, and the deck involution. They preserve $H$ and $\g_P$, so the averaged field still solves the equation and is horizontal. Then $H_0=H-\delta^*X$ is divergence-free and satisfies the linearized Einstein equation.

Its trace $t$ satisfies $\Delta t+\lambda t=0$. Bochner's identity gives
\[
 \int|\nabla^2t|^2\,\dd V
 =\int(\Delta t)^2\,\dd V-\lambda\int|\dd t|^2\,\dd V=0,
\]
so $t$ is constant and hence zero. The transverse-traceless Einstein equation is elliptic, making $H_0$ smooth. Koiso's classification~\cite[Theorem~5.7, pp.~659--660]{Koiso82} gives no nonzero transverse-traceless infinitesimal Einstein deformations of the equal round product. Thus $H=\delta^*X$.

The common fixed surface $T$ of the two circle-coordinate conjugations is a product of great circles, with a flat metric. The averaged field $X$ is tangent to $T$. Since the horizontal part of $H$ is conformal to the quotient metric, $X|_T$ is a conformal vector field on this flat torus. Its components satisfy the Cauchy--Riemann equations and are therefore constant. In global arclength coordinates on $T$, the deck involution acts by $(s_1,s_2)\mapsto(-s_1,-s_2)$, so invariance forces this constant field to vanish. Every torus orbit meets $T$; torus invariance gives $X=0$, hence $H=0$ and $(\phi,\psi,\chi)=0$.

Theorem~\ref{thm:scalar-index-zero} now makes the derivative an isomorphism. Lemma~\ref{lem:half-angle-zero} and Proposition~\ref{prop:localized-degree} give $\deg_2(\mathcal R_{1/2})=1$.
\end{proof}

Finally, we present the proof of our main theorem:


\begin{proof}[Proof of Theorem~\ref{thm:intro-existence}]
Choose $d=m+1$ and $I=[\beta_0,1/2]$, with $0<\beta_0<1/d$. By Propositions~\ref{prop:scalar-zero-equivalence} and~\ref{prop:removability}, it suffices to find a zero of $\mathcal R_{1/d}$ in $\mathcal M$. We show, more generally, that $\mathcal R_\beta$ has a zero for every $\beta\in I$.

Theorems~\ref{thm:scalar-index-zero} and~\ref{thm:compact-zero-set} give the Fredholm index-zero property and compactness of the total zero set. Applying Proposition~\ref{prop:localized-degree} with the projections of Remark~\ref{rem:scalar-degree-projections}, homotopy invariance gives
\[
 \deg_2(\mathcal R_\beta)
 =\deg_2(\mathcal R_{1/2})=1
 \qquad\text{for every }\beta\in I,
\]
where the last equality follows from Proposition~\ref{prop:half-angle-degree}. Nonzero degree implies that $\mathcal R_\beta$ has a zero. In particular, a zero at $\beta=1/d$ determines an Einstein metric $\g\in\mathcal X_{1/d}$, whose pullback by $p_{d-1}$ extends to a smooth Einstein metric with $\lambda>0$ on $\sharp_m(S^2\times S^2)$, with the effective isometric $\TT$-action constructed in $\S$\ref{sec:branched-cover}.
\end{proof}

\bibliographystyle{alpha}
\begingroup
\renewcommand{\baselinestretch}{1}
\selectfont
\bibliography{references}
\endgroup

\end{document}